\documentclass{amsart}

\usepackage[english]{babel}
\usepackage[utf8]{inputenc}
\usepackage{amsmath}
\usepackage{amsfonts}
\usepackage{mathrsfs}
\usepackage{mathtools}
\mathtoolsset{showonlyrefs}

\usepackage{esint}
\usepackage[colorlinks=true, allcolors=red, linktocpage=true]{hyperref}
\usepackage{amssymb}
\usepackage{graphicx}
\usepackage{commath}
\usepackage{cases}
\usepackage{bm}
\usepackage{dsfont}
\usepackage{xcolor}
\usepackage{verbatim}
\usepackage{todonotes}

\usepackage{amsthm}

\theoremstyle{plain}

\newtheorem{thm}{Theorem}[section]
\newtheorem{prop}[thm]{Proposition}
\newtheorem{lem}[thm]{Lemma}

\newtheorem{cor}[thm]{Corollary}

\numberwithin{equation}{section}

\theoremstyle{definition}

\newtheorem{defn}[thm]{Definition}

\newtheorem{rem}[thm]{Remark}
\newtheorem{notation}[thm]{Notation}
\newtheorem{constants}[thm]{Constants}

\makeatletter
\newcommand*{\rom}[1]{\expandafter\@slowromancap\romannumeral #1@}
\makeatother

\newcommand{\ve}{\varepsilon}

\newcommand{\N}{\mathbb{N}}

\newcommand{\R}{\mathbb{R}}

\newcommand{\Z}{\mathbb{Z}}

\newcommand{\calB}{\mathcal{B}}
\newcommand{\calC}{\mathcal{C}}
\newcommand{\calD}{\mathcal{D}}

\newcommand{\calF}{\mathcal{F}}
\newcommand{\calG}{\mathcal{G}}
\newcommand{\calH}{\mathcal{H}}
\newcommand{\calI}{\mathcal{I}}

\newcommand{\calN}{\mathcal{N}}

\newcommand{\calQ}{\mathcal{Q}}

\newcommand{\calT}{\mathcal{T}}

\def\dist{\mathop\mathrm{dist}} 
\def\supp{\mathop\mathrm{supp}} 
\def\diam{\mathop\mathrm{diam}} 

\newcommand{\Layer}{\mathrm{Layer}}

\newcommand{\Stop}{\mathrm{Stop}}
\newcommand{\Top}{\mathrm{Top}}
\newcommand{\cubes}{\mathcal{D}}

\newcommand{\envelope}{Z}

\newcommand{\upperface}{I}
\newcommand{\upperboundary}{F}
\newcommand{\lowerfacebad}{I}
\newcommand{\lowerboundarybad}{F}
\newcommand{\llowerface}{I}
\newcommand{\llowerboundary}{F}
\newcommand{\planecube}{D}
\newcommand{\graphcube}{\sigma}
\newcommand{\graphboundary}{\zeta}
\newcommand{\planeset}{D}
\newcommand{\graphset}{\sigma}

\newcommand{\PI}{{\rm PI}}
\newcommand{\MPI}{{\rm MPI}}
\newcommand{\Lip}{{\rm Lip}}

\def\Xint#1{\mathchoice
	{\XXint\displaystyle\textstyle{#1}}%
	{\XXint\textstyle\scriptstyle{#1}}%
	{\XXint\scriptstyle\scriptscriptstyle{#1}}%
	{\XXint\scriptscriptstyle\scriptscriptstyle{#1}}%
	\!\int}
\def\XXint#1#2#3{{\setbox0=\hbox{$#1{#2#3}{\int}$ }
		\vcenter{\hbox{$#2#3$ }}\kern-.6\wd0}}

\def\dashint{\Xint-}

\newcommand{\surface}[1]{\Sigma_{#1}}
\newcommand{\surfacelift}[1]{\Sigma_{#1}^{*}}
\newcommand{\innercube}{c_0}

\newcommand{\cLB}{c_6}

\newcommand{\csep}{c_4}

\newcommand{\grid}[1]{G_{#1}}
\newcommand{\dyadic}[1]{J_{#1}}

\newcommand{\stopk}{{k_0}}
\newcommand{\gridk}{{k_1}}
\newcommand{\badk}{{k_2}}
\newcommand{\lambdadyadic}{{\lambda_1}}
\newcommand{\goodA}{M_1}

\newcommand{\Afirst}{A_1}
\newcommand{\Asecond}{A_2}

\newcommand{\Ccomp}{C_5}
\newcommand{\Cprop}{C_7}
\newcommand{\Ccover}{C_8}

\newcommand{\Cdoub}{C_{\mathrm{db}}}
\newcommand{\CPI}{C_{\mathrm{PI}}}
\newcommand{\lambdaPI}{\lambda_{\mathrm{PI}}}

\newcommand{\taustop}{\tau_0}
\newcommand{\alphastop}{\alpha_0}
\newcommand{\vestop}{\ve_0}
\newcommand{\thetastop}{\theta_0}
\newcommand{\Astop}{A_0}

\newcommand{\Mstop}{M_0} 
\newcommand{\Mclose}{M_1}
\newcommand{\Cmax}{C_*}

\newcommand{\Stopprime}{\mathrm{Stop}'}

\newcommand{\ar}{\mathrm{ADR}}

\title[UR and PI]{Any uniformly rectifiable set is contained in a set with Poincar\'{e} inequality}

\author{Matthew Hyde, Michele Villa and Ivan Yuri Violo}

\address{Matthew Hyde, \newline
University of Jyv\"askyl\"a
\\
Department of Mathematics and Statistics
\\
P.O. Box 35 (MaD)
\\
FI-40014 University of Jyv\"askyl\"a, Finland}
\email{\href{mailto:matthew.j.hyde@jyu.fi}{matthew.j.hyde@jyu.fi}}

\address{Michele Villa, \newline
IKERBASQUE, Basque Foundation for Science 
Plaza Euskadi 5, 48009 Bilbao, Spain; Euskal Herriko Unibertsitatea/Universidad del Pa\'{i}s Vasco Enparantza Torres Quevedo Ingeniariaren, 1C, Basurtu-Zorrotza, 48013 Bilbo, Bizkaia, Spain
}
\email{\href{mailto:michele.villa@ehu.eus}{michele.villa@ehu.eus}}

\address{Ivan Yuri Violo, \newline
La Sapienza Università di Roma, Dipartimento di Matematica Guido Castelnuovo,
Piazzale Aldo Moro, 5, 00185 Roma RM, Italy
}
\email{\href{mailto:ivanyuri.violo@uniroma1.it}{ivanyuri.violo@uniroma1.it}}

\thanks{M.H.\ is supported by the Research Council of Finland via the project \textit{Quantitative differentiability and rectifiability in metric spaces}, grant no. 363800.}

\thanks{M.V.\ was supported by the European Union MSCA IF QuReViMal, grant no. 101108515, and by an IKERBASQUE starting grant. M.V. was also supported by the Academy of Finland via the project ``Higher dimensional Analyst's Traveling Salesman theorems and Dorronsoro estimates on non-smooth sets", grant No. 347828.}

\thanks{I.Y.V.\ was supported by the European Union (ERC, Con-Fine, 101078057) and the MIUR Excellence Department Project awarded to the Department of
Mathematics, University of Pisa, CUP I57G22000700001}

\begin{document}

\begin{abstract}
    We show that any uniformly $d$-rectifiable subset of $\R^n$ is contained in an Ahlfors $d$-regular set in $\R^{n+1}$ which supports a $(1,1)$-Poincar\'e inequality. This provides a natural converse to a recent result of Azzam. 
\end{abstract}

	\maketitle

\tableofcontents

\section{Introduction}

In the late 1990s problems connected to quasiconformal and quasisymmetric maps led to the development of a Sobolev spaces theory on metric measure spaces (mms). In their foundational work~\cite{heinonen1998quasiconformal}, Heinonen and Koskela introduced the notion of upper gradient, meant to replace the norm of gradients of smooth functions. To define it, let $(X, \mu)$ be a mms. A nonnegative Borel function $g:X \to [0, \infty]$ is said to be an \textit{upper gradient} of a real valued function $u$ if the inequality
\begin{align}
    |u(x)-u(y)| \leq \int_\gamma g
\end{align}
is satisfied for all rectifiable curves $\gamma$ joining $x$ to $y$ in $X$. The path integral is computed with respect to arc length along $\gamma$. Coupling this notion with that of Poincar\'{e} inequality, Shanmugalingam developed a theory of what she termed \textit{Newton-Sobolev} spaces~\cite{shanmugalingam2000newtonian}, denoted by $N^{1,p}(X)$ (see also \cite{cheeger1999differentiability}).

\begin{defn}
    Let $p\in [1,\infty)$. A metric measure space $(X,\mu)$ is said to support a weak $(1,p)$-Poincar\'{e} inequality (or $(1,p)$-PI) if there exist constants $C>0, \Lambda \geq 1$ such that for all open balls $B$ in $X$ and all pairs of functions $u$ and $g$ defined on $B$, whenever $g$ is an upper gradient of $u$ on $B$ and $u$ is integrable on $B$ the following inequality holds:
    \begin{equation}\label{eq:PI}
        \int_B |u-u_B| \, d\mu \leq C \diam(B) \left( \dashint_{\Lambda B}\int g^p \, d\mu\right)^{1/p}.
    \end{equation}
    Any such $(X,\mu)$ is called a PI space.
\end{defn}
The importance of \eqref{eq:PI} for a first-order Sobolev theory lies in the fact that it ensures the presence of `thick' families of rectifiable curves in $X$. When coupled with the doubling property of $\mu$, this ensures that upper gradients reflect appropriately local changes of the associated functions. 
Indeed, Cheeger proved that if $(X,\mu)$ supports a Poincar\'{e} inequality and is doubling, then it possesses an infinitesimal linear structure sufficiently rich to permit upper gradients of Lipschitz functions to be controlled by their point-wise Lipschitz constant. In other words, Cheeger proved in this setting a version of Rademacher's differentiation theorem about Euclidean Lipschitz functions. Spaces $(X,\mu)$ which support this infinitesimal linear structure needed for Cheeger's Redemacher's result to apply were termed \textit{Lipschitz differentiability spaces}. A converse of Cheeger's theorem was proved in the remarkable works by Bate and Li, and Eriksson-Bique~\cite{bate2018differentiability, eriksson2019characterizing}, where they showed that any Lipschitz differentiability space (with the Radon-Nikodym Property) is PI-rectifiable -- it can be `covered' by a countable number of PI-spaces.

The property of supporting a Poincar\'{e} inequality is quantitative and multiscale. As such, one would expect it to have quantitative and multiscale consequences. Cheeger's theorem, however, is qualitative: it says that Lipschitz functions on $X$ are \textit{infinitesimally} linear in an appropriate generalised sense. In Euclidean setting, there is a quantitative version of Rademacher's theorem stating that, in fact, Lipschitz functions are close to being affine at the typical location and scale -- this is Dorronsoro's theorem~\cite{dorronsoro1985characterization}. This begs the question: does an appropriate version of this theorem hold on PI spaces? And, conversely, do spaces with a Dorronsoro theorem also have PI?

To proceed, recall that a metric measure space $(X,\mu)$ is Ahlfors $d$-regular if $X=\supp(\mu)$ and $\mu$ is Ahlfors regular, meaning that there is $A \geq 1$ such that 
\begin{equation}
    A^{-1} r^d \leq \mu(B(x,r)) \leq Ar^d \,\,\, \mbox{ for all } x \in X, \,\, 0<r<2\diam(X).
\end{equation}
When $(X,\mu)$ is Ahlfors $d$-regular, and $X \subset \R^n$, we now know precisely the condition for which $(X,\mu)$ supports a Dorronsoro type theorem: \textit{uniform $d$-rectifiability} is that condition. This is due to David and Semmes~\cite{david1993analysis} and the proof of the WALA conjecture by Merlo~\cite{merlo2026wala}. 
We recall the definition of uniform rectifiability (introduced in \cite{david1991singular}) below.
\begin{defn}
An Ahlfors $d$-regular set $X \subset \R^n$ is said to be uniformly $d$-rectifiable if there exist constants $0<\theta<1$ and $L \geq 1$ such that, for each ball $B$ centered on $X$, there exists a $L$-Lipschitz map $g: \R^d \supset B(0,r(B)) \to \R^n$ such that 
\begin{equation*}
	\mathcal{H}^d(X \cap B \cap g(B(0,r(B)))) \geq \theta r(B)^d.
\end{equation*}
\end{defn}
Returning to the first question above, Azzam \cite{azzam2021poincare} proved that if $X\subset \R^n$ (equipped with the Euclidean metric) is Ahlfors $d$-regular and supports a $(1,d)$-Poincar\'{e}, then it is uniformly $d$-rectifiable. In this article, we prove an appropriate converse. Namely:

\begin{thm}\label{t:main}
	If $E \subseteq B(0,1) \subseteq \R^n$ is uniformly $d$-rectifiable, then there is an  Ahlfors $d$-regular set $Z \subseteq \R^{n+1}$ with $E \subseteq Z$ such that $(Z,\calH^d|_Z)$ supports a $(1,1)$-Poincar\'e inequality. Moreover, the Ahlfors regularity constants and the constants in the Poincar\'e inequality depending only on $n,d$ and the uniform rectifiability constants for $E$. 
\end{thm}

For Azzam, it was enough to assume a $(1,d)$-PI to guarantee uniform rectifiability. Since we construct a space which supports a $(1,1)$-PI, and since a $(1,p)$-PI is stronger than a $(1,q)$-PI whenever $p < q$, it follows that the space we construct is also $d$-UR. One observes also that Theorem \ref{t:main} and Azzam's result imply that any $d$-ADR $(1,d)$-PI space in $\R^n$ is a subset of a $d$-ADR $(1,1)$-PI space in $\R^{n+1}$.

\subsection{Remarks on the literature}
The literature on Lipschitz differentiability, PI spaces and Sobolev spaces on mms is very broad and we only scratch the surface here. First, we should remark that, aside from those already mentioned, two other definition of Sobolev spaces in metric measure spaces where proposed: that of Hajiłasz \cite{hajlasz1996sobolev} and of Ambrosio-Gigli-Savar\'{e}~\cite{ambrosio2014calculus}. All the definitions coincide for a general mms, except for the one of Hajiłasz which coincides with the others  in the presence of a  PI. Cheeger, after his seminal result~ \cite{cheeger1999differentiability}, proved quantitative analogues \cite{cheeger2012quantitative, cheeger2023quantitative}.
While a notion of uniform rectifiability for general metric measure spaces has been introduced in~\cite{bate2023uniformly} by Bate, Schul and the first author, it is unclear under which conditions Cheeger's quantitative result may imply this notion. For geometric quantified flatness estimates across scales related to uniform rectifiability in metric spaces see also \cite{FasslerViolo2025}. As we mentioned above, if $X\subset \R^n$, quantitative affine approximation of Lipschitz function does imply uniform rectifiability. This is the content of~\cite{merlo2026wala}. In fact, if $X \subset \R^n$ is $d$-$\ar$, one may recover the full Dorronsoro characterisation of Sobolev spaces~\cite{azzam2023quantitative}. See also~\cite{hyde2026quantitative} for an analogous result on RCD spaces. 

In general, the qualitative theory of Lipschitz differentiability is much better understood for general metric measure spaces, thanks to the collective efforts of several mathematicians, see e.g. the work of Bate~\cite{bate2020purely, bate2015structure}, of Bate-Li~\cite{bate2017characterizations}, of Cheeger-Kleiner-Schioppa~\cite{cheeger2015infinitesimal}, of Schioppa~\cite{schioppa2016derivations}, of Keith~\cite{keith2004differentiable, keith2004measurable}, of Alberti-Marchese~\cite{alberti2016differentiability}, and of David, Eriksson-Bique and Schul~\cite{david2025coarse} (in no particular order). 

The theory of Sobolev spaces has much progressed, too: we refer the reader to the monograph \cite{heinonen2015sobolev}. Remarkable, and connected to the progress on Lipschitz differentiability and directions present in $X$ mentioned above, is the construction of Sobolev differential by Eriksson-Bique and Soultanis~\cite{eriksson2024curvewise}.

\subsection*{Outline of the proof}

	Consider $E$ as in the statement of Theorem \ref{t:main}. Since $E$ is uniformly rectifiable, it admits a corona decomposition, i.e., it is well-approximated in \textit{most} locations and scales by a controlled `number' of Lipschitz graphs (and contained in their union up to a set of measure zero). It is well-known that Lipschitz graphs admit $(1,1)$-Poincar\'e inequalities. Thus, the basic idea is build bridges between the graphs from the corona decomposition in such a way that generates a global Poincar\'e inequality. 
	
	There are several problems we need to overcome. First, we only have access to graphs in the \textit{good} scales and locations for $E$ but we need to describe the coarse behaviour of the enveloping set in \textit{all} cubes. This begs the question, what does the set look like in the bad locations and scales and how do you glue this to the good graphs in a way that is consistent with a Poincar\'e inequality? We solve this by considering approximations by $d$-dimensional dyadic skeleta, i.e., $d$-dimensional boundary of dyadic cubes. Instead of using the graphs directly, our approximations in the good regions will be unions of dyadic skeleta which are `aligned' along the original graphs and appear at the scale of the minimal cubes. In the bad locations we have no geometry to deal with, so we just look at a fine-mesh approximation of ambient $n$-dimensional ball. See Section \ref{s:skeleta-per-tree}. The benefit of considering such a approximations is the following: if two dyadic skelata with comparable size have large intersection, then their union supports a Poincar\'e inequality (see Section \ref{s:skeleta-per-cube}). Since the good and bad regions are built using dyadic skeleta, this makes the gluing procedure far easier. One caveat of using skeletal approximations is that it is not immediate that they support a Poincar\'e inequality. Thus, we must check that the skeletal approximations admit a Poincar\'e inequality. We check this condition in Proposition \ref{l:sigma-good}. Note that, since we want to make use of the above good intersection properties of nearby skeleta, we will want to consider \textit{smooth} trees where nearby minimal cubes have comparable sizes. These are defined in Section \ref{s:surfaces}. 
	
	The next problem is to avoid bad intersections of the skeletal approximations and the bridges that we build between each of the skeletal approximations. For this we add an additional dimension and `lift' the approximations in a controlled way which depends on how well it approximates $E$, see \eqref{e:surface-lift}. To connect one skeletal approximation to the next layer of skeletal approximations, we pick a single $d$-dimensional skeleton for each minimal cube and `drag' one of its $(d-1)$-dimensional faces down to the approximation below it (here, \textit{down} and \textit{below} refer the height in this additional dimension). One should have in mind connecting the two sets with a bridge which is bi-Lipschitz to a $d$-dimensional cube, although the reality is a little more complicated that. These bridges are built in Section \ref{s:bridges}. 
	
	Finally, we have to verify a Poincar\'e inequality for the set we construct as above. We start in Section \ref{s:conse-layers} by showing a local Poincar\'e inequality, i.e., we check that a skeletal approximation, together with all its connecting bridges, supports a Poincar\'e inequality. To check a global Poincar\'e inequality we find, for each pair of points, a connecting `Harnack chain' of balls such that consecutive balls lie in one of the spaces that already have Poincar\'e estimates. This gives us a way to estimate the oscillation of functions in consecutive balls and, thus, this type of connectivity with local Poincar\'e estimates is enough to conclude a global Poincar\'e inequality. This is carried out in Section \ref{s:global} and relies on the results in Section \ref{s:PI}.

    \subsection*{Statement on AI} No AI tools were used for this article. 

\section{Preliminaries}

 \subsection{Notation} 
 Throughout the paper $C$ will denote some absolute constant (depending perhaps on $n$) which may change from line to line. If there exists $C \geq 1$ such that $a \leq Cb$, then we will write $a \lesssim b$. If $C$ depends on parameters $(t_1,\dots,t_m)$, then we will write $a \lesssim_{t_1,\dots,t_m} b$. We will write $a \sim b$ if $a \lesssim b$ and $b \lesssim a$. The notation $a \sim_{t_1,\dots,t_m} b$ is defined similarly. 

 For a metric space $(X,d)$, we let $B(x,r)$ denote the open in $X$ centred at $x$ with radius $r$, that is, 
 \begin{align}
     B(x,r) = \{y \in X \colon d(y,x) < r \}. 
 \end{align}
 If $B$ is a ball in $X$, we will often write $x_B$ and $r_B$ to denote its centre and radius, respectively (although such objects are not necessarily uniquely defined). For sets $E,F \subseteq X$, let 
 \begin{align}
     \dist(E,F) &\coloneqq \inf\{d(x,y) \colon x \in E \mbox{ and } y \in F \};\\
     \diam(E) &\coloneqq \sup\{d(x,y) \colon x,y \in E \}. 
 \end{align}
 The \textit{$k$-dimensional Hausdorff measure} of $E$ is defined as 
 \begin{align}
     \calH^k(E) \coloneqq \lim_{\delta \to 0} \left\{ \sum_i \diam(U_i)^k \colon E \subseteq \bigcup_i U \mbox{ and } \diam(U_i) < \delta \right\}.
 \end{align}

		For a set $A \subseteq \R^n \times \{0\}$ and a function $f \colon A \to \R,$ denote \textit{the graph of $f$ over $A$} by $\text{Graph}_A(f)$, that is, 
		\begin{align}\label{e:graph-notation}
        \text{Graph}_A(f) = \{ (x,f(x))  : x \in A\} \subseteq \R^{n+1}.
        \end{align}

\subsection{Dyadic cubes and their skeleta}\label{s:dyadic-cubes}
Let $\Delta$ denote the collection of dyadic cubes on $\R^n.$ For $I \in \Delta$, let $x_I$ denote its \textit{centre} and $\ell(I)$ denote its \textit{side-length}. For $c > 0$, let $cI$ denote the cube concentric to $I$ with side-length equal $\ell(cI) = c\ell(I).$ Let 
\begin{align}\label{e:B_I}
	B_I = B(x_I,\diam(I)).
\end{align}
We let $\partial_d I$ denote its $d$-dimensional skeleton, i.e., if $I = \prod_{j=1}^n I_j$ for some (half open, half closed) intervals $I_j \subseteq \R$, then  
\begin{align}
    \partial_d I = \left\{ (x_1,\dots,x_n) \in \overline{I} \colon x_j \in \partial I_j \mbox{ for all least $n-d$ indices } j \in \{1,\dots,n\}  \right\}
\end{align}

These dyadic skeleta will be the basic building block from which we construct the set $Z$ from Theorem \ref{t:main}. We work with these sets because it is easy to describe how nearby 
skeleta interact with each other. The properties we need are proven below. First, we show that the gluing metric between nearby cubes respect the ambient Euclidean distance (a sort of orthogonality). Second, we show that skeleta corresponding to nested cubes have large intersections. Both of these properties will be useful in showing that unions of such skeleta support a Poincar\'e inequality.

\begin{lem}\label{l:reverse-triangle}
	Let $I,J , K \in \Delta$ be such that $I \subseteq K$, $\partial_d I \cap \partial_d K \neq\emptyset$, $\ell(I) = 2^{-m} \ell(K)$ and $\ell(J) =\ell(K)$. Then, for each $x \in \partial_d I$ and $y \in \partial_d J$, we have 
	\begin{align}\label{e:reverse-triangle}
		|x-y| \sim_m \inf_{z \in \partial_d I \cap \partial_d K} |x-z| + |y-z|. 
	\end{align}
\end{lem}

\begin{proof}
	By scaling and translating, we may assume that $K = [0,1]^n$. Translating and rotating again, we may also assume $I \subseteq [0,1/2]^n$. Let us write $x = (x_1,\dots,x_n)$ and $y = (y_1,\dots,y_n)$. Since $x \in \partial_d I$, there is a set $A_x$ with $|A_x| = n-d$ such that 
	\begin{align}
		x_k \in \{j 2^{-m} \colon j \in \{0,1,\dots,2^{m-1}\}\mbox{ for all } k \in A_x \}. 
	\end{align}
	Similarly, there exists a set $A_y$ with $|A_y| = n-d$ such that 
	\begin{align}
		y_k \in \Z \mbox{ for all } k \in A_y. 
	\end{align}
	
	Suppose to begin with that there exists $k \in A_x$ such that $x_k \neq 0$. In this case, we have $x_k \in [2^{-m},1/2]$. Thus, since $y_k \in \Z$, this implies $|x_k - y_k| \geq 2^{-m}$. By taking an arbitrary point $z \in \partial_d I \cap \partial_d K$, since $x \in K$, we get 
	\begin{align}
		|x-z| + |y-z| \leq 2|x-z| + |x-y| \lesssim_m 2^{-m} + | x-y| \lesssim |x-y|. 
	\end{align}
	Now assume that $x_k = 0$ for all $k \in A_x$. Since $|A_x| = n-d$ and $I \subseteq K$ we have $x \in \partial_d I \cap \partial_d K$. The estimate \eqref{e:reverse-triangle} is now immediate by taking $z = x$ as the infimising point. 
\end{proof}

\begin{lem}\label{l:sub-face}
	Let $I \subseteq J$ and suppose that $x \in \partial_d I \cap \partial_d J$. Then, there are $d$-dimensional faces $F_I \subseteq \partial_d I$ and $F_J \subseteq \partial_d J$ such that $x \in F_I \subseteq F_J$. 
\end{lem}

\begin{proof}
	Let us write $I = \prod_{k=1}^n I_k$ and $J = \prod_{k=1}^n J_k$, where $I_k \subseteq J_k$ are half-open half-closed intervals in $\R$. Let $x = (x_1,\dots,x_k) \in \partial_d I \cap \partial_d J$ and let 
	\begin{align}
		K_I \coloneqq \{k  \colon x_k \in \partial I_k\} \quad \mbox{ and } \quad K_J \coloneqq \{k \colon x_k \in \partial J_k\}.
	\end{align}
	Since $I_k \subseteq J_k$ for each $k$, it is clear that $K_J \subseteq K_I$. Since $x \in \partial_d J$, we have that $|K_J| \geq n-d$. Thus, we can find a set $K$ such that 
	\begin{align}
		|K| = n- d \quad \mbox{ and } \quad K \subseteq K_I \cap K_J. 
	\end{align}
	We define the faces $F_I$ and $F_J$ as 
	\begin{align}
		F_I &\coloneqq \{z \in I \colon z_k = x_k \mbox{ for all } k \in K \}; \\
		F_J &\coloneqq \{z \in J \colon z_k = x_k \mbox{ for all } k \in K\}. 
	\end{align}
	If $z = (z_1,\dots,z_n) \in F_I$, then we have $z_k = x_k$ for all $k \in K$ and $z_k \in I_k \subseteq J_k$ for all $k \in \{1,\dots,n\} \setminus K$. That is, $z \in F_J$, as required.
\end{proof}

\subsection{Christ-David cubes} 

Throughout the paper, it will be convenient for us to work with a version of dyadic cubes tailored to a set $X$. These are the so-called Christ-David cubes, which were first introduced by David \cite{david1988morceaux} and generalized in \cite{christ1990b}. The following formulation is due to \cite{hytonen2012non}. 

\begin{lem}\label{cubes}
	Let $X$ be a doubling metric space and $\{X_k\}_{k\in \Z}$ be a sequence of maximal $\varrho^k$-separated nets such that $X_{k+1} \subseteq X_k$, where $\varrho = 1/1000$ and let $\innercube = 1/500.$ Then, for each $k \in \Z$, there is a collection $\mathcal{D}_k$ of ``cubes'' such that the following conditions hold.
	\begin{enumerate}
		\item For each $k \in \Z, \ X = \bigcup_{Q \in \calD_k}Q.$
		\item If $Q_1,Q_2 \in \calD \coloneqq \bigcup_{k}\calD_k$ and $Q_1 \cap Q_2 \not= \emptyset,$ then $Q_1 \subseteq Q_2$ or $Q_2 \subseteq Q_1.$ 
		\item For $Q \in \calD_k$ set $\ell(Q) = 5\varrho^k.$ Then there is $x_Q \in X_k$ such that
		\begin{align*}
			B_X(x_Q,c_0\ell(Q)) \subseteq Q \subseteq B_X(x_Q , \ell(Q)). 
		\end{align*}
	\end{enumerate}
\end{lem}
For each $Q\in \calD$ we set $B_Q\coloneqq B_X(x_Q , \ell(Q)).$

Given a collection of cubes $\calD$ and $Q  \in \calD,$ define 
\begin{align}\label{d:subcubes}
	\calD(Q) = \{R \in \calD : R \subseteq Q\}.
\end{align}
Let $\text{Parent}(Q)$ be the unique cube $R \in \calD$ such that $R \supseteq Q$ and $\ell(Q) = \varrho \ell(Q).$ Let $\text{Child}(Q)$ denote the collection of cubes $R \subseteq Q$ such that $Q = \text{Parent}(R).$ Let $\text{Sibling}(Q)$ be the collection of cubes $R$ such that $\text{Parent}(R) = \text{Parent}(Q).$ 

\begin{rem}\label{r:contained-cubes}
	If $A > 2$ and $Q,R \in \calD$ are such that $Q \subseteq R$ it is easy to see that $AB_Q \subseteq AB_R.$ For a proof see \cite[Lemma 2.10]{azzam2018analyst}.
\end{rem}

Fix for the moment a doubling metric space and a system of Christ-David cubes $\cubes$ on $X$.

\begin{defn}
    A collection of cubes $S \subseteq \cubes$ is called a \textit{stopping-time region} if the following conditions holds. 
    \begin{enumerate}
        \item There is a \textit{top cube} $Q(S) \in S$ such that, if $Q \in S$, then $Q \subseteq Q(S)$.
        \item If $Q \in S$ and $Q \subseteq R \subseteq Q(S)$, then $R \in S$.
        \item If $Q \in S$, then $R \in S$ for all $R \in \text{Sibling}(Q)$. 
    \end{enumerate}
\end{defn}

\begin{defn}\label{d:packing}
	A collection of cubes $\calB \subseteq \calD$ satisfies a \textit{Carleson packing condition} if there exists $C > 0$ such that
	\begin{align}
		\sum_{\substack{Q \in \calB \\ Q \subseteq Q_0}} \ell(Q)^n \leq C\ell(Q_0)^n
	\end{align}
	for all $Q_0 \in \calD.$ 
\end{defn}

\begin{defn}
    For a collection of cubes $\calC \subseteq \cubes$, define the \textit{distance function} to $\calC$ by 
    \[ d_\mathcal{C}(x) = \inf\{\ell(R) + \text{dist}(x,R) : R \in \mathcal{C} \}. \] 
For $Q \in \cubes,$ set
\[ d_\mathcal{C}(Q) = \inf\{ d_\mathcal{C}(x) : x \in Q \} = \inf \{ \ell(R) + \dist(Q,R) : R \in \mathcal{C} \}. \]
\end{defn}

The following can be easily deduced from the definitions.
\begin{lem}\label{l:cubes-bounded overlap}
	Let $A,C \geq 1$ and $c \in (0,1)$. Suppose $\{B_i\}_{i \in I}$ is a collection of balls such that $cB_i \cap cB_j = \emptyset$ and $C^{-1}r_i \leq r_j \leq C r_i$ for all $i,j \in I$ such that $AB_i \cap AB_j \neq \emptyset$, then $\{AB_i\}_{i \in I}$ has bounded overlap depending on $d,c,A$ and $C$. Let $\{Q_i\}_{i \in I}$ be a collection of disjoint Christ-David cubes. If $C^{-1}\ell(Q_j) \leq \ell(Q_i) \leq C\ell(Q_j)$ for all $i,j \in I$ such that $AB_{Q_i} \cap AB_{Q_j} \neq \emptyset$, then $\{AB_{Q_i}\}_{i \in I}$ has bounded overlap with constant depending on $d,\innercube,A$ and $C$. 
\end{lem}

\subsection{Gluing Ahlfors regular sets} 

At many points in the paper we will need to glue various configurations of (upper, lower or upper and lower) regular sets. Here, we gather some useful result that will allow us to preserves these regularity properties after the gluing. We start with the following easy lemma regarding gluing finitely many spaces with non-empty intersection. 

\begin{lem}\label{l:ADR-gluing}
	Let $Z$ be a metric space and $\{X_i\}_{i=1}^N \subseteq Z$ be a finite collection of Ahlfors $d$-regular subsets with uniform regularity constant $C_0$. If $X_i \cap X_{i+1} \neq \emptyset$ for each $i \in {1,\dots,N-1}$, then $\bigcup_{i=1}^N X_i$ is Ahlfors $d$-regular with regularity constant depending only on $C_0$ and $N$. 
\end{lem}

\begin{proof}
    By induction, it is enough to show that $X \cup Y$ is Ahlfors $d$-regular whenever $X$ and $Y$ are Ahlfors $d$-regular with uniform constant $C_0$ and satisfy $X \cap Y \neq \emptyset$. Upper regularity with constant $2C_0$ follows directly from the upper regularity of $X$ and $Y$. For the lower regularity, fix a ball $B$ centred on $X \cup Y$ with $r_B \in (0,2\diam(X \cup Y))$ and assume (without loss of generality) that $B$ is centred at $x_B\in X$. If $r_B \le 4\diam(X)$, then the lower regularity follows from the lower regularity of $X$. Suppose instead that $r_B > 4\diam(X)$. Observe first that 
    \begin{align}\label{e:diameter-Y}
        \diam(Y) > r_B/4. 
    \end{align}
    Indeed, if this were not true, then we would have $\diam(X \cup Y) < 2r_B$, which is a contradiction. Now, take any point $z \in X \cap Y$ and consider the ball $B' = B(z,r_B/2)$. Since $d(x_B,z)<r_B/4$, we have that $B' \subseteq B$, so the lower regularity of $B$ follows from the lower regularity of the space $Y$ and the diameter estimate \eqref{e:diameter-Y}. 
\end{proof}

If the spaces satisfy certain properties, we are able to obtain a  gluing where the regularity constant is independent of the number of spaces in the collection. 

\begin{lem}\label{l:gluing-ADR}
	Let $\{C_i\}_{i=0}^3 \subseteq [1,\infty)$. Let $(Z,d)$ be a metric space and $X,\{Y_i\}_{i \in I} \subseteq Z$ be Ahlfors $d$-regular subsets with regularity constant $C_0$. Let $\{B_i\}_{i \in I}$ be a collection of balls such that the following conditions hold. 
	\begin{enumerate}
		\item \label{prop:bounded-overlap} The set $\{5C_1B_i\}_{i \in I}$ has overlap bounded by $C_2$. 
		\item \label{prop:cover} If $B_j = B(x_j,r_j)$, then 
		\begin{align}
			x_j \in Y_j, \quad  r_j \leq \max\{\diam(Y_j),C_1\diam(X)\} \quad  \mbox{and} \quad Y_i \subseteq C_1 B_i.
		\end{align}
		\item \label{prop:close} We have 
		\begin{align}
			\dist(x_j,X) \leq C_1 r_i. 
		\end{align} 
	\end{enumerate}
	Then, the set 
	\begin{align}
		Y \coloneqq \bigcup_{i \in I} Y_i 
	\end{align}
	is upper Ahlfors $d$-regular, with constant depending on $d,C_0,C_1$ and $C_2$. Furthermore, if
	\begin{enumerate}
		\setcounter{enumi}{3}
		\item \label{prop:cover-X} $X \subseteq \bigcup_{i \in I} C_1B_i$,
	\end{enumerate}
    and for all all $i,j \in I$ such that $3C_1B_i \cap 3C_1B_j \neq \emptyset$ we have
		\begin{align}
			C_3^{-1} \leq r_i/r_j \leq C_3,
		\end{align}
	then  $Y$ is also lower Ahlfors $d$-regular with $\diam(Y) \sim_{C_1} \diam(X)$ and constant depending on $d,C_0,C_1,C_2$ and $C_3$. 
\end{lem}

\begin{proof}
 
	Let $x \in Y$ and $r \in (0,2\diam(Y))$. We start with the lower regularity estimate.  Let $i(x) \in I$ be such that $x \in Y_i \subseteq C_1 B_{i(x)}$.
	
	If  $C_1 r_{i(x)} \geq r/{(100C_3)}$, then the lower regularity estimate in $Y$ for $B(x,r)$ (with constant depending on $C_1$) follows by considering a ball $B \subseteq B(x,r)$ centred on $Y_{i(x)}$, with radius comparable to $r$, and using the lower regularity of $Y_{i(x)}$ (which we have by assumption). Suppose instead that 
	\begin{align}\label{e:assump-LR}
		C_1r_{i(x)} < r/(100C_3).
	\end{align} 
    We also note that 
    $$\diam(X)\ge r/10.$$ Suppose indeed that $\diam(X)< r/10$. Then $C_1r_i<r/10$ for all $i$, because otherwise $X\subset 3C_1B_i$ by \eqref{prop:close} and so $3C_1B_i\cap 3C_1B_{i(x)}\neq \emptyset$ and so $C_1r_i\le C_1C_3 r_{i(x)}<C_1r/100$ which is a contradiction. Hence 
    \[
    \diam(Y)\le \diam(X)+2C_1r_i\le r/10+2/5r=r/2,
    \]
    which contradicts our choice of $r.$
	 Let $x^* \in X$ be any point such that $d(x,x^*) \le 2 \dist(x,X)$. By assumption \eqref{prop:close}, we have $|x - x^*| \leq 3C_1r_{i(x)} \leq r/50$.
	Let 
	\begin{align}
		\calI = \{i \in I \colon X \cap B(x^*,r/4) \cap C_1 B_i \neq \emptyset\}. 
	\end{align}
	By assumption \eqref{prop:cover-X}, $\{C_1B_i\}_{i \in \calI}$ form a cover of $X \cap B(x^*,r/4)$. 
    Thus, since $X$ is Ahlfors $d$-regular with $\diam(X)\ge r/10$ and each $Y_i$ is Ahlfors $d$-regular with $\diam(Y_i) \geq r_i$,  we have
	\begin{align}\label{e:LR}
		r^d \lesssim_{d,C_0} \calH^d(X \cap B(x^*,r/4)) \lesssim_{d,C_0,C_1} \sum_{i \in \calI} r_{B_i}^d \lesssim_{d,C_0} \sum_{i \in \calI} \calH^d(Y_i). 
	\end{align}
	We claim that $Y_i  \subseteq B(x,r)$ for all $i \in \calI$. This would be enough to conclude the lower regularity estimate. Indeed the balls $\{C_1B_i\}$ have bounded overlap and the last term  in the above inequality is at most $C_2\calH^d(Y \cap B(x,r))$.
    
    Suppose by contradiction that there exists $i\in \calI$ and  $y\in Y_i \setminus B(x,r)$. Then, since $B(x^*,r/4) \cap C_1 B_i\neq \emptyset$,
    \[
    r\le |x-y|\le |x-x^*|+|x^*-y|< r/50+r/4+2C_1r_i\le r/2+2C_1r_i.
    \]
    Hence $r< 4C_1r_i$ and thus also $|x-x_i|< r/2+C_1r_i\le 3C_1r_i$, which gives  $x\in  3C_1B_i.$ Therefore, $3C_1B_{i(x)}\cap3C_1B_{i}\neq \emptyset $ and by assumption $r_i\le C_3r_{i(x)}\le \frac{r}{100C_1}$, which is in contradiction  with $r_i>r/(4C_1)$ obtained above. \\
	
	We now prove the upper regularity. Let
	\begin{align}
		\mathcal{J} \coloneqq \{i \in I \colon Y_i \cap B(x,r)\neq \emptyset \}.
	\end{align} 
	By definition, the set $\{Y_i\}_{i \in \mathcal{J}}$ form a cover of $Y \cap B(x,r)$. Thus, 
	\begin{align}
		\calH^d(Y \cap B(x,r)) \leq \sum_{i \in \mathcal{J}} \calH^d(Y_i \cap B(x,r)). 
	\end{align} 
	To estimate the right-hand side, we split $\mathcal{J}$ in two sub-families, those whose corresponding ball $C_1B_i$ is large relative to $B(x,r)$ and those whose corresponding ball is relatively small. Indeed, let 
	\begin{align}
		\mathcal{J}_{\rm large} \coloneqq \{i \in \mathcal{J} \colon C_1r_i \geq r/2 \} \quad \mbox{and} \quad \mathcal{J}_{\rm small} \coloneqq \{i \in \mathcal{J} \colon C_1r_i < r/2\}. 
	\end{align}
	
	We start by showing that 
	\begin{align}
		\sum_{i \in \mathcal{J}_{\rm large}} \calH^d(Y_i \cap B(x,r)) \lesssim_{C_0,C_2} r^d. 
	\end{align}
	Since each $Y_i$ is upper regular with constant $C_0$, it suffices to show that 
	\begin{align}\label{e:-large-I}
		|\mathcal{J}_{\rm large}| \leq C_2. 
	\end{align}
	This follows from Lemma \ref{l:gluing-ADR} \eqref{prop:bounded-overlap}. Indeed, let $i_* \in \mathcal{J}_{\rm large}$ be the index whose corresponding ball has the largest radius. Since $C_1B_{i} \cap B(x,r) \neq \emptyset$ and $C_1r_i > r/2$ for all $i \in \mathcal{J}_{\rm large}$, it follows that $C_1B_i \subseteq 5C_1B_{i_*}$ for all $i \in \mathcal{J}_{\rm large}$. The estimate \eqref{e:-large-I} now follows directly from the bounded overlap of $\{5C_1B_i\}_{i \in I}$.
	
	We finish by showing that 
	\begin{align}
		\sum_{i \in \mathcal{J}_{\rm small}} \calH^d(Y_i \cap B(x,r)) \lesssim r^d. 
	\end{align}

    We further divide $\mathcal{J}_{\rm small}^1=\{j \in \mathcal{J} \ : \  \diam(X)>r_i\}$ and $\mathcal{J}_{\rm small}^2=\mathcal{J}\setminus \mathcal{J}_{\rm small}^1$.
	For each $i \in \mathcal{J}_{\rm small}^1$, let $z_i \in X \cap 2C_1B_i$ (which exists by Lemma \ref{l:gluing-ADR} \eqref{prop:close}) and let $B_i' = B(z_i,r_i)$. Note that $B_i' \subseteq 2C_1B_i$ which, by Lemma \ref{l:gluing-ADR} \eqref{prop:bounded-overlap}, implies
	\begin{align}
		\{B_i'\}_{i \in \mathcal{J}_{\rm small}} \mbox{ has overlap bounded by } C_2. 
	\end{align}
	Note also, since $C_1B_i \cap B(x,r) \neq \emptyset$ and $C_1r_i < r/2$, we have $B_i' \subseteq 2C_1B_i \subseteq B(x,3r)$. 
    Moreover, since $B_i'\cap X\neq \emptyset$ we have $B_i' \subseteq B(x^*,6r)$ for some $x^*\in X$. Putting this together, and using the Ahlfors regularity of $X$ and $\{Y_i\}$, $\diam(X)>r_i,$ we have 
	\begin{align}
		\sum_{i \in \mathcal{J}_{\rm small}} \calH^d(Y_i \cap B(x,r)) &\lesssim_{d,C_0,C_1} \sum_{i \in \mathcal{J}_{\rm small}} \calH^d(X \cap B_i') \\
		&\lesssim_{d,C_0,C_1} C_2\calH^d(X \cap B(x^*,4r)) \lesssim_{d,C_0,C_1,C_2}  r^d.
	\end{align}
	To conclude the proof it suffices to show $|\mathcal{J}_{\rm small}^2|\le C_2$ and then argue as for $\mathcal{J}_{\rm large}$.  Since $\diam(X)\le r_i$ and $d(X,Y_i)\le C_1r_i$,
    it holds $X\subset 2C_1B_i$ for all $i\in \mathcal{J}_{\rm small}^2.$ From this $|\mathcal{J}_{\rm small}^2|\le C_2$ follows from the bounded overlap assumption.\\

    Finally, we prove the diameter estimate $\diam(Y) \sim_{C_1} \diam(X)$ under the assumption \eqref{prop:cover-X}. The estimate $\diam(Y) \lesssim_{C_1} \diam(X)$ follows from the triangle inequality, using the facts that $\diam(Y_i) \leq 2C_1r_i$, $\dist(x_j,X) \leq C_1 r_i$ and $r_i \leq C_1\diam(X)$ for all $i \in I$. For the estimate $\diam(X) \lesssim_{C_1} \diam(Y)$, consider two point $x,x' \in X$ such that $d(x,x') \geq \diam(X)/2$. By assumption, there are indices $i,j \in I$ such that $x \in C_1B_i$ and $x' \in C_1B_j$. If one of $r_i$ or $r_j$ is larger than $\diam(X)/(8C_1)$, then the diameter estimate follows directly from the fact that $\diam(Y_k) \geq r_k$. In the other case that $r_i,r_j \leq \diam(X)/(8C_1)$, we have 
    \begin{align}
        d(x_i,x_j) &\geq d(x,x') - d(x_i,x) - d(x_j,x')\\
        &\geq \diam(X)/2 - C_1r_i - C_1 r_j \geq \diam(X)/4,
    \end{align}
    which finishes the estimate. 
\end{proof}

\section{Poincar\'e inequalities}\label{s:PI} 

There are a number of key results we will use to check a Poincar\'e inequality in the set $\envelope$ that we construct for Theorem \ref{t:main}. The first is a well-known characterisation of Poincar\'e inequalities via maximal function estimates. To state a more precise result, we will introduce some notation.

\begin{notation}
	Let $C,\lambda \geq 1$, $(X,d,\mu)$ be a metric measure space, $B$ a ball in $X$ and $u \colon X \to \R$ Lipschitz. 
	\begin{enumerate}
		\item We say $(u,B) \in \PI(C,\lambda)$ if the estimate 
		\begin{align}
			\dashint_{B'} |u - u_B| \, d\mu \leq C \dashint_{\lambda B'} \Lip  u \, d\mu 
		\end{align}
		holds for all $B' \subseteq B$. \\
		\item We say $(u,B) \in \MPI(C,\lambda)$ if the estimate
		\begin{align}
			\hspace{2em} |u(x) - u(y)| \leq C d(x,y)(M_{\lambda d(x,y)}[\Lip u]^p(x) + M_{ \lambda d(x,y)}[\Lip u]^p(y))^\frac{1}{p}
		\end{align}
		holds for all $x,y \in B$.
	\end{enumerate}  
\end{notation}
Above we used the maximal function notation $M_R g (x)\coloneqq \sup_{r<R}\fint_{B(x,r)}|g| d \mu.$

The following is proved in \cite[Lemma 5.15]{heinonen1998quasiconformal}. Although not stated as such in \cite{heinonen1998quasiconformal}, it can be deduced by examining the proof there.

\begin{lem}\label{l:PI-max}
	Let $C,\lambda > 1$, $(X,d,\mu)$ be a locally compact doubling metric space with doubling constant $\Cdoub$, $B$ a ball in $X$ and $u \colon X \to \R$.  Then,
	\begin{align}
		(u,B) \in \PI(C,\lambda) \implies (u,\tfrac{1}{10}B) \in \MPI(C',10\lambda).
	\end{align} 
	and 
	\begin{align}
		(u,B) \in \MPI(C,\lambda) \implies (u,B) \in \PI(C',2\lambda)
	\end{align}
	In both cases $C'$ depends quantitatively on $C$ and $\Cdoub$. In particular, $X$ admits a weak $(1,p)$-Poincar\'e inequality if and only if there is a constant $C,\lambda \geq 1$ such that 
	\begin{align}\label{e:PI-max}
		|u(x) - u(y)| \leq C d(x,y)(M_{ \lambda d(x,y)}[\Lip u]^p(x) + M_{ \lambda d(x,y)}[\Lip u]^p(y))^\frac{1}{p}
	\end{align}
	for all $x,y \in X$ and $u \colon X \to \R$ Lipschitz. 
\end{lem}

As a consequence of Lemma \ref{l:PI-max} one can show that the pointwise estimate from the $\MPI(C,\lambda)$ implies control of averages of Lipschitz functions on nearby balls. The proof (with slightly different constants) can be found in \cite{heinonen2012lectures}.  
Below we denote  $M_{R',R} g (x)\coloneqq \sup_{R'<r<R}\fint_{B(x,r)}|g| d \mu.$

\begin{lem}\label{l:variation-in-balls}
	Let $C,\lambda > 1$, $(X,d,\mu)$ be a locally compact doubling metric space with doubling constant $\Cdoub$, $B_0$ a ball in $X$ with radius $r_0$ and $u \colon X \to \R$ Lipschitz. If $(u,B_0) \in \MPI(C,\lambda)$ then
	\begin{align}
		| u_B - u_{B'}| \lesssim_{C',\Cdoub} r_{B_0} \left( M_{r_B,4\lambda r_{B_0}} \Lip (x_B) + M_{r_{B'},4\lambda r_{B_0}} \Lip(x_{B'}) \right) 
	\end{align}
	for all balls $B,B' \subseteq \tfrac{1}{4}B_0$. 
\end{lem}

The construction of the set $Z$ from Theorem \ref{t:main} also relies on gluing PI spaces together. The following tells us when a gluing preserves a Poincar\'e inequality and is proved in \cite{heinonen1998quasiconformal}.

\begin{lem}\label{l:gluing}
	Suppose $E,F \subseteq \R^n$ are Ahlfors $d$-regular, each admitting a weak $(1,1)$-Poincar\'e inequality. Suppose the set $A = E \cap F$ satisfies 
	\begin{align}
		\calH^d(A \cap B) \geq C^{-1} r_B^d
	\end{align}
	for all balls $B$ centred on $A$ with radius at most $\min\{\diam(E),\diam(F)\}$ and that, for any $x \in E$ and $y \in F$, we have 
	\begin{align}\label{e:gluing-metric-statement}
		\inf_{a \in A} |x-a| + |a-y| \lesssim |x-y|. 
	\end{align} 
	Then, $E \cup F$ admits a weak $(1,1)$-Poincar\'e inequality. 
\end{lem}

Lemmas \ref{l:diameter-PI} and \ref{l:checking-the-gluing-condition} below provide us with some convenient ways to check the condition  \eqref{e:gluing-metric-statement}.

\begin{lem}\label{l:diameter-PI}
	Let $(Z,d)$ be a metric space and $X,Y \subseteq Z$ such that $A \coloneqq X \cap Y \neq \emptyset$. If $x \in X$ and $y \in Y$ are such that 
	\begin{align}
		d(x,y) \gtrsim \min\{\diam(X),\diam(Y)\}
	\end{align}
	then \eqref{e:gluing-metric-statement} holds. 
\end{lem} 

\begin{proof}
	Suppose, without loss of generality, that $\diam(X) \leq \diam(Y)$. Choose an arbitrary point $a \in A$. Since $A \subseteq X$, we have $d(x,a) \leq \diam(X)$. Then, by the triangle inequality,
	\begin{align}
		d(x,a) + d(a,y) \leq 2d(x,a) + d(x,y) \lesssim d(x,y),
	\end{align}
	as required. 
\end{proof}

\begin{lem}\label{l:checking-the-gluing-condition}
	For each $c_1 \in (0,1)$ there is $\tau_1 \in (0,1)$ such that the following holds. Let $\tau \in (0,\tau_1)$, $g \colon \R^n \to \R$ be $\tau$-Lipschitz and $\Gamma$ be its graph. Let $x,y \in \R^{n+1}$ such that $\pi_1(x) \geq g(\pi_n(x))$ and $\pi_1(y) \leq g(\pi_n(y))$. If $a \in \Gamma$ is such that either 
	\begin{align}\label{e:unnatural}
		c_1|a-x| + g(\pi_n(a)) - g(\pi_n(x)) \leq |\pi_1(a-x)|,
	\end{align}
	or 
	\begin{align}
		c_1|a-y| + g(\pi_n(a)) - g(\pi_n(y)) \leq |\pi_1(a-y)|,
	\end{align}
	then \eqref{e:gluing-metric-statement} holds with this choice of $a$. 
\end{lem}

\begin{proof}
	We will prove \eqref{e:gluing-metric-statement} in the case that \eqref{e:unnatural} holds. The other case is similar, and we will not include the details. We begin by showing how \eqref{e:unnatural} implies the more natural condition
	\begin{align}\label{e:natural}
		c_1|a-x| \leq 2 | \pi_1(a-x) |.
	\end{align}
	Indeed, taking $g(\pi_n(a)) - g(\pi_n(x))$ to the right-hand side of \eqref{e:unnatural} and using that $g$ is $\tau$-Lipschitz and $\pi_n$ is 1-Lipschitz, we get
	\begin{align}
		c_1|a-x| \leq |\pi_1(a-x)| + \tau |a-x|. 
	\end{align}
	If $\tau_1 = \tau_1(c_1)$ is sufficiently small, this implies \eqref{e:natural}. We now use \eqref{e:natural} to prove \eqref{e:gluing-metric-statement}. By the triangle inequality, we have $|x-a| + |a-y| \leq |x-y| + 2|a-x|$. Thus, it suffices to show that 
	\begin{align}\label{e:a-y}
		|a-x| \lesssim |x-y|. 
	\end{align} 
	Assume to begin with that 
	\begin{align}
		\pi_n(y) \notin B \coloneqq B_n(\pi_n(x), |\pi_1(a - x)| ) \subseteq \R^n. 
	\end{align}
	In this case, we have by \eqref{e:natural} that  
	\begin{align}
		|x-y| \geq |\pi_n(x -y)| \geq |\pi_1(a-x)| \gtrsim |a-x|,
	\end{align}
	which proves \eqref{e:a-y}. Suppose instead that 
	\begin{align}
		\pi_n(y) \in B .
	\end{align}
	Recalling \eqref{e:natural}, we also have $\pi_n(a) \in 2C_1B$, where $C_1 \coloneqq c_1^{-1}$. In particular, $|\pi_n(y) - \pi_n(a)| \leq 4C_1|\pi_1(a-x)|$. Since  $g$ is $\tau$-Lipschitz, setting $y'=(\pi_n(y),g_n(\pi_n(y))$ we have $|x-y'|\le \sqrt{1+\tau^2}|x-y|$. Using these facts, and the fact that $g$ is $\tau$-Lipschitz,  if $\tau_1 = \tau_1(c_1)$ is small enough, we have 
	\begin{align}
		\sqrt{1+\tau^2}|x-y| &\geq |x-y'|\ge |\pi_1(x - y')| \geq |\pi_1(a-x)| - |\pi_1(y'-a)| \\
		&= |\pi_1(a-x)|  - |g(\pi_n(y)) - g(\pi_n(a))|\\
		&\geq |\pi_1(a-x)| - 4C_1\tau | \pi_1(a-x)| \geq (1/2)|\pi_1(a-x)| \gtrsim |a-x|,
	\end{align}
	which gives the required estimate provided $\tau_1$ is small enough. 
\end{proof}

\subsection{Poincar\'e inequalities via Harnack chains}

In this and the next subsection we collected two results which give us ways of transferring local Poincar\'e inequalities to global inequalities. The result below states, roughly speaking, that if every two points is connected by a Harnack chain of balls, where consecutive balls lie in a common PI space, then the whole space has PI.  

\begin{lem}\label{l:local-to-global}
	Let $A > 1$, $(X,d,\mu)$ be a metric measure space where $\mu$ is Ahlfors $d$-regular with constant $\Cdoub$, $u \colon X \to \R$ be Lipschitz, and $x,y \in X$. Suppose there exists a sequence of balls $\{B_i\}_{i=1}^N$ in $X$ and subsets $\{X_i\}_{i=1}^{N}$ of $X$ such that the following conditions holds. 
	\begin{enumerate}
		
		\item \label{l:local-to-global-PI} Each $X_i$ is Ahlfors $d$-regular with constant $\Cdoub$ and supports a $(1,p)$-Poincar\'e inequality with constants $(\CPI,\lambdaPI)$.
		\item \label{l:local-to-global-balls} For each $i \in \{1,\dots,N-1\}$, the space $X_i$ completely contains the balls $B_i$ and $B_{i+1}$ i.e., 
		\begin{align}
			B_i \subseteq X_i  \quad \mbox{ and } \quad B_{i+1} \subseteq X_i. 
		\end{align}
		\item \label{l:local-to-global-X} The end spaces contain the end points and balls i.e., 
		\[ \{x\} \cup B_{1} \subseteq X_1 \quad \mbox{ and } \quad \{y\} \cup B_{{N}} \subseteq X_N.\]
		\item \label{l:local-to-global-ends} If $z_i$ and $r_i$ denote the centre and radius, respectively, of the ball $B_i$, then there exists $K \in \{1,\dots,N\}$ such that 
		\begin{align}
			x &\in  B(z_i, A r_i) \quad \mbox{for all } i \in \{1,\dots,K\}; \label{e:x-in-ball} \\
			y & \in  B(z_i,A r_i) \quad \mbox{for all } i \in \{K,\dots,N\}.
		\end{align}
		\item \label{l:local-to-global-r} We have 
		\begin{align}
			\sum_{i=1}^N r_i \leq A d(x,y). 
		\end{align}
	\end{enumerate}
	Then, there is $\CPI' = \CPI'(A,\Cdoub,\CPI,d,p)$ and $\lambdaPI' = \lambdaPI'(A,\lambdaPI)$ such that 
	\begin{align}\label{e:maximal-estimate}
		\hspace{2em} |u(x) - u(y)| \leq \CPI' d(x,y) \left( M_{\lambdaPI' d(x,y)}[\emph{Lip}u]^p(x) + M_{\lambdaPI' d(x,y)}[\emph{Lip}u]^p(y) \right)^\frac{1}{p} .
	\end{align}
\end{lem}

\begin{proof}
	We begin by showing that
	\begin{align}
		| u_{B_i} - u_{B_{i+1}} | &\lesssim_{A,\Cdoub,\CPI} \max\{r_i, r_{i+1}\} \left(M_{\lambdaPI'  d(x,y)} [\text{Lip}u]^p(x)\right)^\frac{1}{p};  \label{e:av-1} \\
		| u_{B_j} - u_{B_{j+1}} | &\lesssim_{A,\Cdoub,\CPI} \max\{r_i, r_{i+1}\} \left(M_{\lambdaPI'  d(x,y)} [\text{Lip}u]^p(y)\right)^\frac{1}{p}, \label{e:av-2}
	\end{align}
	for each $i \in \{1,\dots,K-1\}$ and $j \in \{K,\dots,N-1\}$, where
	\begin{align}
		\lambdaPI' \coloneqq960\lambdaPI A^3.
	\end{align}
	We shall prove \eqref{e:av-1} only (the proof for \eqref{e:av-2} is the same, up to notational changes). Fix some $i \in \{1,\dots,K-1\}$ and assume, without loss of generality, that $r_{i+1} \geq r_i$. Under this assumption, we have by \eqref{l:local-to-global-ends} that $B_i \subseteq 3A B_{i+1}$. Hence $B_i,B_{i+1}\subset \frac{1}{4} 12AB_{i+1}$. Since $X_i$ contains $B_i \cup B_{i+1}$ and supports a $(1,p)$-Poincar\'e inequality, it follows from Lemma \ref{l:variation-in-balls}  and Lemma \ref{l:PI-max} that 
	\begin{align}\label{e:chain-1}
		| u_{B_i} - u_{B_{i+1}} | &\lesssim_{\Cdoub,\CPI} Ar_{i+1} \left( M^{X_i}_{r_i , 480\lambdaPI A r_{i+1}}[\text{Lip}u]^p(z_i) \right)^\frac{1}{p} \\
		&\hspace{3em} + Ar_{i+1}\left( M^{X_i}_{r_{i+1} , 480\lambdaPI A r_{i+1}}[\text{Lip}u]^p(z_{i+1}) \right)^\frac{1}{p}    
	\end{align}
	(recall that $z_i$ is the centre of $B_i$), where $M^{X_i}$ denotes the maximal function with distance and measure restricted to $X_i.$ For $r \in [r_i, 480\lambdaPI A r_{i+1}]$, it follows from \eqref{l:local-to-global-ends} and \eqref{l:local-to-global-r} that 
	\begin{align}
		B(z_i,r) \subseteq B(x,(A+1)r) \subseteq B(x,480(A+1) \lambdaPI Ar_{i+1}) \subseteq  B(x,\lambdaPI' d(x,y)).
	\end{align}
	Thus, since $B(z_i,r_i)\subset X_i$ and by the Ahlfors $d$-regularity of $X_i$ and $X$,
	\begin{align}\label{e:chain-2}
		\quad\fint_{B(z_i,r)\cap X_i} [\text{Lip}u]^p \, d\mu \lesssim_{A,\Cdoub} \dashint_{B(x,(A+1)r)} [\text{Lip}u]^p \, d\mu \leq M_{\lambdaPI'  d(x,y)} [\text{Lip}u]^p(x). 
	\end{align}
	Similarly, for $r \in [r_{i+1}, 480\lambdaPI A r_{i+1}]$,
	\begin{align}\label{e:chain-3}
		\dashint_{B(z_{i+1},r)\cap X_i} [\text{Lip}u]^p \, d\mu \lesssim_{A,\Cdoub} M_{\lambdaPI'  d(x,y)} [\text{Lip}u]^p(x). 
	\end{align}
	Inequality \eqref{e:av-1} now follows from \eqref{e:chain-1}, \eqref{e:chain-2} and \eqref{e:chain-3}. \\ 
	
	Since $x \in B(x_1,Ar_1) \subseteq B(x_1,A^2d(x,y))$ and $y \in B(x_N,Ar_N) \subseteq B(x_N,A^2d(x,y))$, it follows from Lemma \ref{l:variation-in-balls} that  
	\begin{align}
		|u(x) - u_{B_1}| &\lesssim_{A,\Cdoub,\CPI}  d(x,y)\left( M_{4\lambdaPI A^2d(x,y)} [\text{Lip}u]^p(x) \right)^\frac{1}{p}; \label{e:av-3} \\
		|u(y) - u_{B_N}| &\lesssim_{A,\Cdoub,\CPI}  d(x,y)\left( M_{4\lambdaPI A^2d(x,y)} [\text{Lip}u]^p(y) \right)^\frac{1}{p}. \label{e:av-4}
	\end{align}
	
	Now, putting \eqref{e:av-1}, \eqref{e:av-2} together with \eqref{e:av-3}, \eqref{e:av-4}, and applying \eqref{l:local-to-global-r}, we get 
	\begin{align}
		| u(x) - u(y) | &\leq |u(x) - u_{B_1} | + \sum_{i=1}^{N-1} | u_{B_i} - u_{B_{i+1}} | + | u_{B_N} - u(y) | \\
		&\hspace{-2em}\lesssim_{A,\Cdoub,\CPI} \big(d(x,y)+\sum_{i=1}^N r_i\big) \left[ \left( M_{\lambdaPI'  d(x,y)} [\text{Lip}u]^p(x) \right)^\frac{1}{p} + \left(M_{\lambdaPI'  d(x,y)} [\text{Lip}u]^p(y) \right)^\frac{1}{p} \right] \\
		&\hspace{-2em}\lesssim_A d(x,y) \left[ M_{\lambdaPI'  d(x,y)}[\text{Lip}u]^p(x) + M_{\lambdaPI'  d(x,y)}[\text{Lip}u]^p(y) \right]^\frac{1}{p}.
	\end{align}
	This completes the proof of \eqref{e:maximal-estimate}. 
\end{proof}

\subsection{Poincar\'e inequalities via Whitney balls} 
Let $Z$ be an ambient metric space. The next result says that if a space $X \subseteq Z$ is covered by a collection of balls resembling a Whitney decomposition (i.e., well-separated with nearby balls having comparable size), which have PI estimates locally around each ball, and cover nearby PI space $Y \subseteq Z$, then $X$ also has PI.

	\begin{prop}\label{p:local-to-global-prime}
	Let $\{C_i\}_{i=0}^{5} \subseteq [1,\infty)$ such that 
	\begin{align}\label{e:C_4}
		C_3 \geq 5C_1 \quad \mbox{ and } \quad  C_5 \geq 50C_1.
	\end{align} 
	Let $(Z,d)$, $X,\{Y\}_{i \in I} \subseteq Z$ and $\{B_i\}_{i \in I}$ satisfy the hypotheses Lemma \ref{l:gluing-ADR} \eqref{prop:bounded-overlap}--\eqref{prop:cover-X} with constants $C_0,C_1$ and $C_2$. Suppose further that $X$ support a $(1,1)$-Poincar\'e inequality with constant $(\CPI,\lambdaPI)$ and that the following condition holds
	\begin{enumerate}
		\item \label{prop:separated} The exists $\csep > 0$ such that, whenever $i \neq j \in I$, we have 
		\begin{align}
			\csep B_i \cap \csep B_j = \emptyset .
		\end{align}
		\item \label{prop:comp-size} For all $i,j \in I$ such that $C_3B_i \cap C_3B_j \neq \emptyset$ we have
		\begin{align}
			C_4^{-1} \leq r_i/r_j \leq C_4
		\end{align}
		\item \label{prop:small-balls} For all $i \in I$ and $u \colon X \to \R$ we have $(u,C_5B_i) \in \MPI(C_6,C_6)$. 
	\end{enumerate}
	Then, the set 
	\begin{align}
		Y = \bigcup_{i \in I} Y_i 
	\end{align}
	is Ahlfors $d$-regular, with regularity constant depending on $d,C_0,C_1$ and $C_2$, and supports a $(1,1)$-Poincar\'e inequality, with constant $\CPI'$ depending on all the above constants (with the exception of $C_3$) and $d$, and $\lambdaPI'$ depending on $\lambdaPI,C_1,C_5$ and $C_6$.
\end{prop}

The Ahlfors regularity follows directly from Lemma \ref{l:gluing-ADR}. We will let $C_0(Y)$ denote the regularity constant of $Y$ for remainder of the section. Our goal now it to prove a Poincar\'e inequality, which will be done over a number of lemmas. By Lemma \ref{l:PI-max}, it suffices to find $C_7 \geq 1$ such that 
\begin{align}\label{e:conclusion-all-balls}
	\hspace{1.5em} |u(x_0) - u(y_0)| \leq C_7 |x_0 - y_0| \left( M_{C_7|x_0-y_0|} \Lip  u(x_0) + M_{C_7|x_0-y_0|} \Lip  u(y_0)  \right) 
\end{align}
for any $u \colon Y \to \R$ Lipschitz and $x_0,y_0 \in Y$, which we now fix. If $x_0,y_0 \in C_5B_i$ for some $i \in I$, then this follows by assumption, so let us suppose for the moment that there is no such $i \in I$. Note that, by Lemma \ref{l:variation-in-balls} and the hypothesis \eqref{prop:small-balls} in Proposition \ref{p:local-to-global-prime}, 
\begin{align}\label{l:controlled-balls}
	|u_B - u_{B'} | \lesssim_{C_0(Y) , C_5, C_6} r_i \left( M_{r_B , \lambda r_i} \Lip  u(x_B) + M_{r_{B'}, \lambda r_i} \Lip  u(x_{B'}) \right)
\end{align}
for all indices $i \in I$ and balls $B,B' \subseteq (C_5/4)B_i$, where 
\begin{align}
	\lambda \coloneqq 4 C_5 C_6. 
\end{align}

We define functions $\tilde{u} , \tilde{\rho} \colon X \to \R$ as follows. First, let $\{\phi_i\}_{i\in I}$ be a collection of Lipschitz functions such that 
\begin{align}\label{e:phi_i}
	\supp(\phi_i) \subseteq 2C_1B_i , \ \Lip (\phi_i) \lesssim r_i^{-1} \mbox{ and } \sum_{i \in I} \phi_i = 1 \mbox{ on } \bigcup_{i \in I} C_1B_i \supseteq X ,
\end{align}
where the final inclusion follows from hypothesis \eqref{prop:cover} in Proposition \ref{p:local-to-global-prime}. The existence of such a collection follows from the fact that $Y$ is Ahlfors regular. Similarly, let $\{\varphi_i\}_{i \in I}$ be a collection of Lipschitz functions such that
\begin{align}
	\supp(\varphi_i) \subseteq 3C_1B_i , \ \Lip (\varphi_i) \lesssim r_i^{-1} \mbox{ and }  \varphi_i \gtrsim 1 \mbox{ on } \bigcup_{i \in I} 2C_1B_i. 
\end{align} 
Then, set 
\begin{align}
	\tilde{u} \coloneqq \sum_{i \in I} \langle u \rangle_{B_i} \phi_i \quad \mbox{ and } \quad \tilde{\rho} \coloneqq \sum_{i \in I} \langle \Lip  u \rangle_{2 \lambda B_i} \varphi_i.
\end{align}
We will show via several lemmas the following. 

\begin{lem}\label{l:Lipu-less-rho}
	The function $\tilde{u}$ is Lipschitz and satisfies 
	\begin{align}
		\Lip \tilde{u} \lesssim_{C_0(Y),C_2,C_4,C_5,C_6} \tilde{\rho}.
	\end{align} 
\end{lem}

To help with the proof of Lemma \ref{l:Lipu-less-rho}, we will need some notation. For each $y \in Y$, choose a point $p_X(y) \in X$ and an index $i(y) \in I$ such that 
\begin{align}
	d(y,p_X(y)) = \dist(y,X) \quad \mbox{ and } \quad y \in C_1 B_{i(y)}. 
\end{align}
For each $z \in Z$, let 
\begin{align}
	I(z) = \{i \in I \colon z \in 4C_1B_i\}
\end{align}
and 
\begin{align}\label{e:j(z)-prime}
	j(z) \coloneqq \arg \max_{i \in I(z)} r_i.  
\end{align}

\begin{lem}\label{l:j(z)}
	For each $z \in Z$ and $i \in I(z)$, we have 
	\begin{align}\label{e:j(z)}
		r_i \sim_{C_4} r_{j(z)}, \quad \bigcup_{i \in I(z)} 2B_i \subseteq (C_5/4) B_{j(z)} \quad \mbox{ and } \quad |I(z)| \leq C_2,
	\end{align}
	where $|\cdot|$ denotes the cardinality of a set
\end{lem}

\begin{proof}
	The right-most property follows from Lemma \ref{l:gluing-ADR} \eqref{prop:bounded-overlap} and the definition of $I(z)$. For the remaining properties, fix some $z \in Z$ and $i \in I(z)$. Since $i,j(z) \in I(z)$, the left-most property follows directly from hypothesis \eqref{prop:comp-size} in Proposition \ref{p:local-to-global-prime} \eqref{prop:comp-size} since $C_3 \geq 5C_1$. Since $4C_1B_i \cap 4C_{1}B_{j(z)} \neq\emptyset$ and $r_{j(z)} \geq r_{i}$, we have $4C_1B_i \subseteq 12C_1B_{j(z)}$. Thus, the middle claim in \eqref{e:j(z)} follows from the fact that $C_5 \geq 50C_1$.
\end{proof}

\begin{lem}\label{l:I(w)}
	For each $y \in Y$, $w \in B(p_X(y),r_{i(y)})$ and $i \in I(w)$, we have 
	\begin{align}
		B(y,r_{i(y)}) \cup B_i \subseteq (C_5/4) B_{j(w)}.  
	\end{align}
\end{lem}

\begin{proof}
	Fix $y,w$ and $i \in I(w)$ as above. By Lemma \ref{l:j(z)}, we immediately have $B_i \subseteq (C_5/4) B_{j(w)}$. Since $y \in C_1 B_{i(y)}$ and $\dist(x_{i(y)},X) \leq C_1r_{i(y)}$, it follows that
	\begin{align}\label{e:tilde-y}
		|y-p_X(y)| \leq 2C_1 r_{i(y)}.
	\end{align}
	Since $y \in C_1 B_{i(y)}$ and $w \in B(p_X(y),r_{i(y)})$, this gives $|w - x_{i(y)}| \leq (1+3C_1)r_{i(y)} \leq 4C_1 r_{i(y)}$ i.e., $i(y) \in I(w)$.  Since $B(y,r_{i(y)}) \subseteq 2B_{i(y)}$, the inclusion $B(y,r_{i(y)}) \subseteq (C_5/4)  B_{j(y)}$ follows from Lemma \ref{l:j(z)}.
\end{proof}

\begin{proof}[Proof of Lemma \ref{l:Lipu-less-rho}]
	That $\tilde{u}$ is Lipschitz is easy to see since it is the sum of Lipschitz functions whose supports have bounded overlap by $C_2$. To prove the estimate on $\Lip  \tilde{u}$, fix some $x \in X$ and let $j \coloneqq j(x)$ as in \eqref{e:j(z)-prime}. Recalling the properties of $\phi_i$ in \eqref{e:phi_i}, and applying Lemma \ref{l:variation-in-balls}, Lemma \ref{l:j(z)} and \eqref{l:controlled-balls}, we have  
	\begin{align}
		\Lip  \tilde{u}(x) &= \Lip  \left(\tilde{u} - \langle u \rangle_{B_j}  \right)(x) =\Lip \left( \sum_{i \in I} \left(\langle u \rangle_{B_i} - \langle u \rangle_{B_{j}}\right)\phi_i \right)(x) \\
		&\lesssim \sum_{i \in I(x)} r_{i}^{-1}  | \langle u \rangle_{B_i} - \langle u \rangle_{B_{j}} | \\
		& \lesssim_{C_2,C_4} \sup_{i \in I(x)} r_j^{-1} | \langle u \rangle_{B_i} - \langle u \rangle_{B_{j}} | \\
		&\lesssim_{C_0(Y),C_5,C_6} \sup_{i\in I(x)} \left( M_{r_i,\lambda  r_j} \Lip  u (x_i) + M_{r_j, \lambda r_j} \Lip  u (x_j) \right)
	\end{align}
    (recall that $x_i$ is the centre of $B_i$).
	Fix an arbitrary $i \in I(x)$. Since $2B_i \subseteq (C_5/4)B_{j}$, any ball $B(x_i,r)$ with $r \in [r_i,\lambda r_j]$ must satisfy $B(x_i,r) \subseteq 2\lambda B_j$ and $r \sim_{C_4,\lambda} r_j$. On the other hand, any ball $B(x_j,r)$ with $r \in [r_j , \lambda r_j]$ must satisfy $B(x_j,r) \subseteq B(x_j,2\lambda r_j)$ and $r \sim_\lambda r_j$. By comparing the averages of $\Lip  u$ in these balls, it follows from this that 
	\begin{align}
		M_{r_i,\lambda  r_j} \Lip  u (x_i) + M_{r_j, \lambda r_j} \Lip  u (x_j) \lesssim_{C_0(Y),C_4,\lambda} \langle \Lip  u \rangle_{2\lambda B_j}. 
	\end{align}
	Recalling that $x \in 2C_1B_j$ and $\varphi_j(x) \gtrsim 1$ on $2C_1B_j$, we now have 	
	\begin{align}
		\Lip  \tilde{u}(x) \lesssim_{C_0(Y),C_2,C_4,C_5,C_6}  \langle\Lip  u \rangle_{2\lambda B_j} \varphi_j(x) \leq \tilde{\rho}(x),
	\end{align}
	as required. 
\end{proof}

For the remainder of this section, we will denote $\nu \coloneqq \calH^d|_X$.  To avoid confusion over which measure various average integrals are being computed with respect to, for a function $f \colon Z \to \R$ and a ball $B \subseteq Z$, we will write 
\begin{align}
    \langle f \rangle_{B,\nu} \coloneqq \dashint_B f \, d\nu.   
\end{align}

Lemma \ref{l:Lipu-less-rho} and the assumptions on $X$ tell us that we can estimate the oscillations of $\tilde{u}$ by the integral of $\tilde{\rho}$ with respect to $\nu$. To prove the Poincar\'e inequality on $Y$, we must now prove estimates which related $u$ to $\tilde{u}$ and the integral of $\Lip u$ with repsect to $\mu$ to the integral of $\tilde{\rho}$ with respect to $\nu$. 

\begin{lem}\label{l:u-from-average}
	We have 
	\begin{align}\label{e:u-from-average}
		|u(y) - \langle \tilde{u} \rangle_{B(p_X(y),r_{i(y)}) , \nu} | \lesssim_{C_0(Y),C_4,C_5,C_6} r_{i(y)} M_{C_8 r_{i(y)}} \Lip  u(y)
	\end{align}
	for all $y \in Y$, where $C_8 \coloneqq 2C_4\lambda$. 
\end{lem}

\begin{proof}
	Fix some $w \in B(p_X(y),r_{i(y)})$ and $i \in I(w)$ for the moment. Recalling from Lemma \ref{l:I(w)} that $B(y,r_{i(y)}) \cup B_i \subseteq (C_5/4) B_{j(w)}$, by applying \eqref{l:controlled-balls}, we get 
	\begin{align}
		|u(y) - \langle u \rangle_{B_i}| &\leq |u(y) - \langle u \rangle_{B(y,r_{i(y)})}| + |\langle u \rangle_{B(y,r_{i(y)})} - \langle u \rangle_{B_i}|  \\
		&\lesssim_{C_0(Y),C_5,C_6} r_{i(y)} M_{\lambda r_{i(y)}} \Lip u(y) + r_{j(w)} M_{r_{i(y)},\lambda r_{j(w)}}  \Lip u(y) \\
		&\hspace{2em} +  r_{j(w)} M_{r_{i},\lambda r_{j(w)}}  \Lip u(x_i).
	\end{align}
	Arguing as in the proof of Lemma \ref{l:Lipu-less-rho}, any ball $B(x_{j(w)},r)$ which satisfies $r \in [r_{i(y)} , \lambda r_{j(w)}]$ is contained in $B(y,2\lambda r_{j(w)}) \subseteq B(y,2\lambda C_4 r_{i(y)})$ and satisfies $r \sim_{C_4,\lambda}  r_{i(y)}$. Using also the fact that $r_{j(w)} \leq C_4 r_{i(y)}$, we can now estimate
	\begin{align}
		|u(y) - \langle u \rangle_{B_i}| \lesssim_{C_0(Y),C_4,C_5,C_6} r_{i(y)} M_{2 C_4 \lambda r_{i(y)}} \Lip  u(y). 
	\end{align}
	Letting $C_8 \coloneqq 2C_4\lambda$, the above estimate implies 	
	\begin{align}
		|u(y) - \langle \tilde{u} \rangle_{B(p_X(y),r_{i(y)}) , \nu}| &\leq \dashint_{B(p_X(y),r_{i(y)})} | u(y) - \tilde{u}(w) | \, d\nu(w) \\
		&=\dashint_{B(p_X(y),r_{i(y)})} |u(y) - \sum_{i \in I(w)} \langle u \rangle_{B_i} \phi_i(w) | \, d\nu(w)  \\
		&\leq \dashint_{B(p_X(y),r_{i(y)})} \sum_{i \in I(w)} |u(y) - \langle u \rangle_{B_i} | \phi_i(w)  \, d\nu(w)\\
		&\lesssim \dashint_{B(p_X(y),r_{i(y)})}  \sum_{i \in I(w)} r_{i(y)} M_{{C_8}r_{i(y)}} \Lip  u(y) \phi_i(w)      \, d\nu(w) \\
		&= r_{i(y)} M_{{C_8}r_{i(y)}} \Lip  u(y),
	\end{align}
	where the implicit constant depends on $C_0(Y),C_4,C_5$ and $C_6$. This finishes the proof.  
\end{proof}

\begin{lem}\label{l:r-is-small}
	Let $y \in Y$ and $r \geq r_{i(y)}$. If $i \in I$ is such that $B(p_X(y),r) \cap B_i \neq\emptyset$, then we have $r_i  \leq \max\{3+2C_1,C_4\} r$. 
\end{lem}

\begin{proof}
	Suppose towards a contradiction that there exists $i \in I$ such that $B(p_X(y),r) \cap B_i \neq \emptyset$ and $r_i > Ar$, where $A = \max\{3+2C_1,C_4\}$. Recalling that $y \in B_{i(y)}$ and from \eqref{e:tilde-y} that $|p_X(y) - y| \leq 2C_1 r_{i(y)} \leq 2C_1 r$, it follows from that triangle inequality that 
	\begin{align}
		B_{i(y)} \subseteq B(y,2r_{i(y)}) \subseteq B(p_X(y),(2+2C_1)r_{i(y)}) \subseteq B(p_X(y),(2+2C_1)r). 
	\end{align}
	By our choice of $A$ and the triangle inequality, we have $B_{i(y)} \subseteq 2B_i$. Thus, it follows from Proposition \ref{p:local-to-global-prime} \eqref{prop:comp-size} that 
	\begin{align}
		C_4r \leq Ar < r_i \leq C_4 r_{i(y)} \leq C_4 r,
	\end{align}
	which is a contradiction.	
\end{proof}

\begin{lem}\label{l:rho-less-rho}
	For each $y \in Y$ and $r \geq r_{i(y)}$ we have 
	\begin{align}
		\dashint_{B(p_X(y),r)} \tilde{\rho} \, d\nu \lesssim_{C_0,C_0(Y),C_1,C_4,\csep,\lambda} \dashint_{B(y,C_9r)} \Lip u \, d\mu,
	\end{align}
	where $C_9$ depends only on $\lambda, C_1$ and $C_5$. 
\end{lem}

\begin{proof}
	By a standard volume argument, using Ahlfors regularity of $Y$, Proposition \ref{p:local-to-global-prime} \eqref{prop:separated} and \eqref{prop:comp-size}, we know that the balls $\{2 \lambda  B_i\}_{i \in I}$ have bounded overlap with constant depending only on $C_0(Y),C_4,\csep$ and $\lambda$. By Lemma \ref{l:r-is-small} and the estimate $|y - p_X(y)| \leq 2C_1 r_{i(y)} \leq 2C_1 r$ (see \eqref{e:tilde-y} for the first estimates), if $i \in I$ is such that $B_i \cap B(p_X(y),r) \neq\emptyset$ then it is easy to check from the triangle inequality that 
	\begin{align}
		2\lambda B_i \subseteq C_9B(y,r)
	\end{align}
	for some $C_9$ depending only on $\lambda, C_1$ and $C_4$. Combining these two facts we can estimate 
	\begin{align}
		\int_{B(p_X(y),r)} \tilde{\rho} \, d\nu &= \int_{B(p_X(y),r)} \sum_{i \in I} \langle \Lip u \rangle_{2\lambda B_i} \varphi_i \, d\nu \\ 
		&\leq \sum_{\substack{i \in I \\ B_i \cap B(p_X(y),r) \neq \emptyset}} \int_{B_i} \langle \Lip u \rangle_{2 \lambda B_i} \, d\nu  \\
		&\lesssim_{C_0,C_0(Y),\lambda } \sum_{\substack{i \in I \\ 2\lambda B_i \subseteq B(y,C_9r) \neq \emptyset}} \int_{2\lambda B_i} \Lip u \, d\mu \\
		&\lesssim_{C_0(Y),C_4,\csep,\lambda} \int_{B(y,C_9r)} \Lip u \, d\mu .
	\end{align}
	The desired estimate follows this since, by Ahlfors regularity of both $X$ and $Y$, $\mu(B(y,C_9 r)) \sim_{C_0,C_0(Y),C_9} \nu(B(p_X(y),r))$. 
\end{proof}

\begin{lem}
	The estimate \eqref{e:conclusion-all-balls} holds for the points $x_0,y_0$ with $C_7$ depending on all constants in the statement of Proposition \ref{p:local-to-global-prime}, except for $C_3$, and $C_8$ depending on $\lambdaPI, C_1, C_5$ and $C_6$. 
\end{lem}

\begin{proof}
	Since we have assumed that there is no $i \in I$ such that $x_0,y_0 \in C_5B_i$ and $C_5 \geq 2$, it follows that $(C_5/2)\max\{  r_{i(x_0)}, r_{i(y_0)} \} \leq |x_0-y_0|$. Furthermore, recalling from \eqref{e:tilde-y} that $|y_0 - p_X(y_0)| \leq 2C_1r_{i(y_0)}$ and  $|x_0 - p_X(x_0)| \leq 2C_1r_{i(x_0)}$, and from the statement that $C_5 \geq 50C_1$, we have by the triangle inequality that 
	\begin{align}\label{e:x-y-tilde}
		\max\{r_{i(x_0)}, r_{i(y_0)}\} \leq 	|p_X(x_0) - p_X(y_0)|  \lesssim |x_0-y_0|. 
	\end{align}
	We use this, the triangle inequality and Lemma \ref{l:u-from-average} to get 
	\begin{align}
		|u(x_0) - u(y_0)| &\leq |u(x_0) -  \langle \tilde{u} \rangle_{B(\tilde{x_0},r_{i(x_0)}) , \nu}| + | \langle \tilde{u} \rangle_{B(\tilde{x_0},r_{i(x_0)}) , \nu} -  \langle \tilde{u} \rangle_{B(p_X(y_0),r_{i(y_0)}) , \nu}| \\
		&\hspace{2em} + | \langle \tilde{u} \rangle_{B(p_X(y_0),r_{i(y_0)}) , \nu} - u(y_0)| \\
		&\lesssim r_{i(x_0)} M_{C_8r_{i(x_0)}} \Lip  u(x_0) + | \langle \tilde{u} \rangle_{B(p_X(x_0),r_{i(x_0)}) , \nu} -  \langle \tilde{u} \rangle_{B(p_X(y_0),r_{i(y_0)}) , \nu}|  \\
		&\hspace{2em} + r_{i(y_0)} M_{C_8r_{i(y_0)}} \Lip  u(y_0) \\
		&\leq |x_0-y_0|\left( M_{C_8|x_0-y_0|} \Lip  u(x_0) + M_{C_8|x_0-y_0|} \Lip  u(y_0) \right) \\
		&\hspace{2em} +  | \langle \tilde{u} \rangle_{B(p_X(x_0),r_{i(x_0)}) , \nu} -  \langle \tilde{u} \rangle_{B(p_X(y_0),r_{i(y_0)}) , \nu}|,
	\end{align}
	where the implicit constant depends on $C_0(Y),C_4,C_5$ and $C_6$.
	By applying the Poincar\'e inequality and Ahlfors regularity on $X$, the final term in the above inequality is bounded above by 
	\begin{align}
		&| \langle \tilde{u} \rangle_{B(p_X(x_0),r_{i(x_0)}) , \nu}  - \langle \tilde{u} \rangle_{B(p_X(x_0),|p_X(x_0) - p_X(y_0)|) , \nu} | \\
		&\hspace{2em}+ |\langle \tilde{u} \rangle_{B(p_X(x_0),| p_X(x_0) - p_X(y_0)|) , \nu} - \langle \tilde{u} \rangle_{B(p_X(y_0),| p_X(x_0) - p_X(y_0)|) , \nu} | \\
		&\hspace{4em} +|\langle \tilde{u} \rangle_{B(p_X(y_0),| p_X(x_0) - p_X(y_0)|) , \nu} - \langle \tilde{u} \rangle_{B(p_X(y_0),r_{i(y_0)}) , \nu} | \\
		&\lesssim_{C_0,\CPI} |p_X(x_0) - p_X(y_0)| M_{r_{i(x_0)} , 3\lambdaPI |p_X(x_0)-p_X(y_0)|} \Lip  \tilde{u}(p_X(x_0)) \\
		&\hspace{2em} + |p_X(x_0) - p_X(y_0)| M_{r_{i(y_0)} , 3\lambdaPI |p_X(x_0)-p_X(y_0)|} \Lip  \tilde{u}(p_X(y_0)).
	\end{align}
	By Lemma \ref{l:Lipu-less-rho}, Lemma \ref{l:rho-less-rho} and \eqref{e:x-y-tilde}, these final terms are  
	\begin{align}
		&\lesssim_{C_0(Y),C_2,C_4,C_5,C_6} |x_0-y_0|  M_{r_{i(y_0)} , 3\lambdaPI |p_X(x_0)-p_X(y_0)|} \Lip  \tilde{\rho}(p_X(y_0)) \\
		&\hspace{10em} + |p_X(x_0) - p_X(y_0)| M_{r_{i(y_0)} , 3\lambdaPI |p_X(x_0)-p_X(y_0)|} \Lip  \tilde{\rho}(p_X(y_0)) \\
		&\lesssim_{C_0,C_0(Y),C_1,C_4,\csep,\lambda} |x_0-y_0| \left(M_{24\lambdaPI C_9 |x_0-y_0| } \Lip u(x_0)+M_{24\lambdaPI C_9 |x_0-y_0| } \Lip u(y_0)\right).		
	\end{align}
	Recalling that $C_9 = C_9(\lambda,C_5,C_6)$, this finishes the proof. 
\end{proof}

\section{The stopping-time regions}\label{s:surfaces}
With the preliminaries out of the way, we can start making some inroads towards Theorem \ref{t:main}. Thus, for the rest of the paper, $E \subseteq \R^n$ will be a fixed set satisfying the hypotheses of that theorem, i.e., a uniformly $d$-rectifiable subset of $B(0,1) \subseteq \R^n$. We will also fix system of Christ-David cubes $\cubes$ on $E$. By scaling and translating, we can assume there is $Q_0 \in \cubes$ such that 
\begin{align}\label{e:largest0-ball}
	Q_0 = E \quad \mbox{ and } \quad B_{Q_0} = B(0,5). 
\end{align}

Our goal in this section is to define a relevant collection of stopping-time regions which will be used in the construction of the approximating surface.

\subsection{The cubes $\mathbf{Layer(N), Stop(N)}$ and $\mathbf{Top(N)}$} Before defining the stopping-time regions we need to introduce several collection of cubes in $\cubes(Q_0)$ which we call $\Layer(N)$, $\Stop(N)$ and $\text{Top}(N)$. We also gather various useful properties that will be used throughout the remainder of the construction. 

\begin{constants}
	Let $\alphastop, \vestop, \thetastop ,\taustop \in (0,1)$, $\Astop >  1$ and $\Mstop = 20\Astop.$ Each of $\alphastop, \vestop, \thetastop ,\taustop$ will be chosen sufficiently small with $\vestop = \vestop(\Astop,\alphastop, \thetastop, \taustop)$. We will choose $\Astop$ sufficiently large. 
\end{constants}

\begin{notation}
	Set 
\begin{align}
	\mathcal{G} =  \{Q \in \cubes : b\beta_E(\Mstop B_Q) \leq \vestop\}.  \quad \mbox{ and } \quad \calB = \cubes \setminus \mathcal{G}.
\end{align}
For $Q \in \calG$, let $P_Q$ be a $d$-plane such that 
\begin{align}\label{e:P_Q}
	b\beta_E(\Mstop B_Q,P_Q) \leq 2b\beta_E(\Mstop B_Q) \leq 2\vestop. 
\end{align}
\end{notation}

The starting point is the following result due to David and Semmes, which is proven in Sections 7 - 11 of \cite{david1991singular}.

\begin{lem}\label{l:stop-angle}
	For $\Astop$ sufficiently large, and $\alphastop$ and $\vestop = \vestop(\Astop,\alphastop)$ sufficiently small, there exists a collection of disjoint stopping-time regions $\calF$ such that the following holds.
	\begin{enumerate}
		\item The family $\calF$ partitions $\cubes(Q_0).$
		\item The collection of cubes $\{Q(S')\}_{S' \in \calF}$ satisfy a Carleson packing condition. 
		\item If $Q(S') \in \calB$ then $S' = \{Q(S')\}.$
		\item We have 
		\[ \bigcup_{\substack{S' \in \calF \\ Q(S') \in \calG}} S' \subseteq \mathcal{G}. \]
		\item For each $S' \in \calF$ such that $Q(S') \in \calG$ and each $Q \in S'$ we have 
		\begin{align}
			\angle(P_{Q(S')},P_Q) \leq \alphastop.
		\end{align}
	\end{enumerate}
\end{lem}

\begin{rem}
	The above formulation differs slightly from the one in \cite{david1991singular} since we are treating \textit{bad cubes} as stopping-time regions with a single cube. For clarity, we have chosen to denote stopping-time regions in $\calF$ with a prime and we will denote stopping-time regions in $\calT$ (to be defined) without a prime. 
\end{rem}

\begin{constants}
	For the remainder of the paper we will suppose the $\Astop,\alphastop$ and $\vestop$ are chosen to satisfy the conditions of Lemma \ref{l:stop-angle}. We will fix $\calF$ as in the conclusion of Lemma \ref{l:stop-angle}. 
\end{constants}

\begin{notation}\label{n:original-stopping-times}
	For $Q \in \cubes(Q_0),$ let 
	\[S_Q' = \{ R \in S' : R \subseteq Q\},\] 
	where $S' \in \calF$ is such that $Q \in S'.$
\end{notation}

\begin{notation}
     For each integer $N \geq 0$, define collections of cubes Layer($N$), Stop($N$) and Up($N$) as follows. First, set 
\[\Top(0) = \{Q_0\}.\] 
Then, assuming $\Top(N)$ has been defined for some $N \geq 0$, let 
\begin{align}
	\Layer(N) &= \bigcup \{S_Q' : Q \in \Top(N) \}; \\
	\Stop(N) &= \{Q \in \cubes : Q \ \text{max. such that there is} \ Q^\prime \in \text{Sibling}(Q)\\
	&\hspace{10em}\text{with} \ \ell(Q^\prime) < \taustop  d_{\Layer(N)}(Q^\prime)\}; \\
	\Top(N+1) &= \bigcup_{Q \in \Stop(N)} \text{Child}(Q). 
\end{align}
For brevity, in the remainder of the paper we will denote 
\begin{align}\label{e:d_N}
	d_N \coloneqq d_{\Layer(N)}.  
\end{align}
\end{notation}

The following lemma says that the size of cubes in $\Stop(N)$ is determined by $d_N$, similar to a Whitney decomposition. 

\begin{lem}[{\cite{azzam2018analyst}}]\label{l:stop-size}
	For $\taustop $ sufficiently small, if $R \in \Stop(N)$ then 
	\begin{align}\label{e:stop-size-1}
		\frac{1}{2}\ell(R) \leq \taustop  d_N(R) \leq \frac{2}{\rho} \ell(R).
	\end{align}
	Moreover, if $C > 1$ and $\taustop  =\taustop (C)$ is sufficiently small then 
	\begin{align}\label{e:sttop-size-2}
		\frac{1}{4} \ell(R) \leq \taustop  d_N(x) \leq\frac{4}{\rho} \ell(R)
	\end{align}
	for all $R \in \Stop(N)$ and $x \in CB_R.$ 
\end{lem}

The following result says that nearby cubes in $\Stop(N)$ have comparable sizes.

\begin{lem}[{\cite{azzam2018analyst}}]\label{l:comp-size}
	Let $C > 1$. If $\taustop  = \taustop (C)$ is sufficiently small, then $\Ccomp^{-1}\ell(R) \leq   \ell(R') \leq \Ccomp \ell(R)$ for all $R,R' \in \Stop(N)$ satisfying $CB_R \cap CB_{R'} \neq\emptyset$. Here $\Ccomp > 1$ is an absolute constant independent of both $C$ and $\taustop.$ 
\end{lem}

\begin{rem}\label{r:Top(N)}
	The above statement remains true with $\Top(N)$ in place of $\Stop(N).$ Indeed, if $R,R' \in \Top(N)$ satisfy $CB_R \cap CB_{R'} \neq \emptyset,$ then $\text{Parent}(R),\text{Parent}(R') \in \Stop(N-1)$ and $CB_{\text{Parent}(R)} \cap CB_{\text{Parent}(R')} \neq \emptyset.$ One then applies Lemma \ref{l:stop-size} for cubes in $\Stop(N-1)$ and uses the fact that $\ell(Q) = \rho \ell(\text{Parent}(Q))$ for any $Q \in \cubes$.
\end{rem}

\subsection{The extended stopping-time regions} We now introduce the stopping-time regions $S_Q$ for each $Q \in \Top(N)$, which is obtained from $S_Q'$ by adding all those cubes between $\min(S_Q')$ and $\Stop(N)$. More precisely,  
\begin{align}
	S_Q = \{R \subseteq Q : R \supseteq T \mbox{ for some } T \in \Stop(N)\} \supseteq S_Q'.
\end{align}
For $N \geq 0$, let  
\[\calT^N = \{S_Q : Q \in \Top(N)\} \quad \mbox{ and } \quad \calT = \bigcup_{N \geq 0 } \calT^N.  \]

Lemma \ref{l:Carleson} and Lemma \ref{l:lip-graphs} below tell us that the collection of stopping-time regions $\calT$ satisfy a Carleson packing condition and that there is a graph which well approximates the stopping-time regions near a \textit{good} stopping-time regions (i.e., a region in which the top cube is flat). Moreover, the approximation between the graph and the set is \textit{bi-lateral}. For the reader familiar with corona decompositions, this formulation implies that usual corona decomposition. The point is that (in this formulation) nearby top cubes have comparable sizes and a graph approximations nearby stopping-time regions as well. Since the proofs are essentially standard (and implicit in \cite{azzam2018analyst}), we have decided not to include them.

\begin{lem}\label{l:Carleson}
	The family $\calT$ partitions $\cubes(Q_0)$ and the collection of cubes $\{Q(S)\}_{S \in \calT}$ satisfies a Carleson packing condition. 
\end{lem}

\begin{lem}\label{l:lip-graphs}
	Let $M \geq 1$ and $\{Q_k\}_{k=0}^K \subseteq \Top(N)$ such that $Q_0 \in \calG$ and $MB_{Q_i} \cap MB_{Q_0} \neq \emptyset$. If $\Astop = \Astop(M)$ is large enough and $\taustop = \taustop(M)$, $\vestop = \vestop(\Astop,\taustop,\thetastop)$ are small enough, there exists a $C\alphastop$-Lipschitz map $g_0 \colon P_{Q_0} \to P_{Q_0}^\perp$ with $\|g\|_\infty \leq \thetastop \ell(Q_0)$ such that the following conditions hold. 
	\begin{enumerate}
		\item \label{prop:E-is-close} If $Q \in S_i$ and $x \in E \cap \Astop B_Q$, then 
		\begin{align}
			|x - \left(\pi_0(x) + g_0(\pi_0(x))\right) |  \leq \thetastop \ell(Q) ,
		\end{align}
		where $\pi_0$ is the orthogonal projection onto $P_{Q_0}.$ In particular, if $\Gamma_0$ is the graph of $g_0$ over $P_{Q_0},$ then 
		\begin{align}
			\dist(x, \Gamma_0)  \leq \thetastop \ell(Q). 
		\end{align}
		\item \label{prop:graph-is-close} If $Q \in S_i$ and $x \in \Gamma_0 \cap \Astop B_Q,$ then
		\begin{align}
			\dist(x,E) \leq \thetastop \ell(Q). 
		\end{align}
	\end{enumerate}
\end{lem}

\subsection{The truncated (extended) stopping-time regions}

In later sections it will be convenient to truncate the stopping-time regions from the previous section at some specified scale. Associated to these will be a truncated distance function. 

\begin{notation}
	For $\stopk \in \N$, let 
	\begin{align}
		\calT^{N}(\stopk) \coloneqq \{S \in \calT^N \colon Q(T) \in \cubes_{\leq \stopk} \}. 
	\end{align}
	For $S \in \calT^N(\stopk)$, let 
	\begin{align}
		S(\stopk) \coloneqq S \cap \cubes_{\leq \stopk}
	\end{align}
	and
	\begin{align}
	\Stop(N,\stopk) \coloneqq \bigcup_{S \in \calT^N(\stopk)} \min(S(\stopk)).
	\end{align}
	Also, let 
	\begin{align}
		\Top(N,\stopk) \coloneqq \{Q(T) \colon T \in \calT^N(\stopk)\}. 
	\end{align}
	Finally, define 
	\begin{align}
		d_{N,\stopk} \coloneqq \max\{d_N , 5\rho^\stopk/\taustop\}. 
	\end{align}
\end{notation}

We shall now prove variants of Lemma \ref{l:stop-size} and Lemma \ref{l:comp-size} for these truncated regions and the truncated distance function.

\begin{lem}\label{l:stop-size-cutoff}
	Let $S \in \calT^N(\stopk)$ and $Q \in \Stop(N,\stopk)$ such that $Q \subseteq Q(S)$. Then, 
	\begin{align}\label{e:stop-size-cutoff-1}
		\ell(Q) \sim \taustop  d_{N,\stopk}(Q) \coloneqq \max\{\taustop  d_N(Q) , 5\rho^{\stopk}\}.
	\end{align}
	More generally, if $C > 1$ and $\taustop = \taustop(C)$ is sufficiently small, then 
	\begin{align}\label{e:stop-size-cutoff-2}
		\taustop d_{N,\stopk}(x) \sim \ell(Q)
	\end{align}
	for all $x \in CB_Q$. 
\end{lem}

\begin{proof}
	We start with \eqref{e:stop-size-cutoff-1}. Since $Q$ contains some element of $\Stop(N)$, it follows from maximality that $\ell(Q) \geq \taustop  d_N(Q)$. Since $Q \in \cubes_{\leq \stopk}$, we also have $\ell(Q) \geq 5\rho^{\stopk}$. These two inequalities imply $\ell(Q) \geq d_{N,\stopk}(Q)$. For the converse, note that if $\ell(Q) = 5\rho^{\stopk}$ then the inequality $\ell(Q) \leq \taustop  d_{N,\stopk}(Q)$ is immediate. Now, let us suppose that $\ell(Q) > 5\rho^{\stopk}$. In this case, $Q \in \Stop(N)$. Thus, $\ell(Q) \lesssim \taustop  d_N(Q) \leq \taustop  d_{N,\stopk}(Q)$ by Lemma \ref{l:stop-size}.
	
	We now prove \eqref{e:stop-size-cutoff-2}. Let $C > 1$, $x \in CB_Q$ and suppose initially that $\taustop$ small enough for \eqref{e:stop-size-cutoff-1}. Let $z \in Q$ such that $d_{N,\stopk}(z) = d_{N,\stopk}(Q)$, so that 
	\begin{align}\label{e:z}
		c\ell(Q) \leq d_{N,\stopk}(z) \leq C' \ell(Q). 
	\end{align}
	Since $\taustop d_{N,\stopk}$ is $\taustop$-Lipschitz, we have 
	\begin{align}\label{e:d(x-z)}
		|\taustop d_{N,\stopk}(x) - \taustop d_{N,\stopk}(z)| \leq 2C\taustop \ell(Q) < \min\{c,C'\}\ell(Q)/2
	\end{align}
	(assuming $\taustop < \min\{c,C'\}/(4C)$). The estimate \eqref{e:stop-size-cutoff-2} then follows directly from \eqref{e:z} and \eqref{e:d(x-z)} after applying the triangle inequality. 
\end{proof}

\begin{lem}\label{l:comp-size-cutoff}
	Given $C > 1$, if $\taustop  = \taustop (C)$ is sufficiently small, then $\Ccomp^{-1}\ell(R) \leq   \ell(R') \leq \Ccomp \ell(R)$ for all $R,R' \in \Stop(N,\stopk)$ satisfying $CB_R \cap CB_{R'} \neq\emptyset$. Here $\Ccomp > 1$ is an absolute constant independent of both $C$ and $\taustop.$ In particular, the set $\{CB_R\}_{R \in \Stop(N,\stopk)}$ has bounded overlap with constant depending on $d,\innercube,\Ccomp$ and $C$.
\end{lem}

\begin{proof}
	The second statement (regarding bounded overlap) is an immediate consequence of the first and Lemma \ref{l:cubes-bounded overlap}. For the first, let $R,R' \in \Stop(N,\stopk)$ satisfy the above hypotheses. If $R,R' \in \cubes_{< \stopk}$ (so that $R,R' \in \Stop(N)$) then the result follows directly from Lemma \ref{l:comp-size}. If $R,R' \in \cubes_{\stopk}$ then $\ell(R) = \ell(R')$ and the result is obvious. Suppose instead that $R \in \cubes_{< \stopk}$ (so that $R \in \Stop(N)$) and $R' \in \cubes_{\stopk}$. By symmetry, this is the only remaining case that we need to consider. Since $\ell(R) > \ell(R')$, we only need to prove the reverse inequality. Since $CB_R \cap CB_{R'} \neq \emptyset$ (and $\ell(R) > \ell(R')$), it follows from the triangle inequality that $R' \subseteq 3CB_R$. Let $R'' \in \Stop(N)$ such that $R'' \subseteq R' \subseteq 3CB_R$. Choosing $\taustop  = \taustop (3C)$ as in Lemma \ref{l:comp-size}, it now follows that 
	\begin{align}
		\ell(R') \geq \ell(R'') \geq \Ccomp^{-1} \ell(R),
	\end{align}
	as required. 
\end{proof}

It turns out that graph constructed in Lemma \ref{l:lip-graphs} are (locally) covered by the corresponding minimal cubes. More precisely, we have the following.

\begin{lem}\label{l:lip-graphs-subset}
	Let $M \geq 1$ and $\{Q_k\}_{k=0}^K$ as in Lemma \ref{l:lip-graphs}, and choose $\Astop = \Astop(M)$ and $\taustop $ such that the conclusions there hold. If $B = B(x_{Q_0},r)$ is a ball such that $r \in ( \thetastop \ell(Q_0), \Astop \ell(Q_0))$ and 
	\begin{align}\label{e:Q_i-covers}
		E \cap 2B \subseteq \bigcup_{k=0}^K Q_k,
	\end{align}
	then 
	\begin{align}
		\Gamma_0 \cap B \subseteq \bigcup_{k=0}^K  \bigcup_{\substack{R \in \min(S_{Q_k}(\stopk)) \\ \Gamma_0 \cap B \cap 2B_R \neq \emptyset  }} 2B_R , 
	\end{align}
\end{lem}

\begin{proof}
	Fix a ball $B$ (as in the statement) and a point $z \in \Gamma_0 \cap B$. Let $m \geq 0$ be the largest integer such that 
	\begin{align}
		z \in \bigcup_{k=0}^K \bigcup_{Q \in S_{Q_k} \cap \cubes_m} \Astop B_Q, 
	\end{align}
	which exists since each $S_{Q_k}(\stopk)$ is a finite collection of cubes. It satisfies $5\rho^m \leq \ell(Q_0)$, since $z \in \Astop B_{Q_0}$. By Lemma \ref{l:lip-graphs} \eqref{prop:graph-is-close}, there exists $y \in E$ such that $|y-z| \leq 5\thetastop \rho^m \leq \thetastop \ell(Q_0)$. Since $z \in B$ and $r_B \geq \thetastop \ell(Q_0)$, this implies that $y \in E \cap 2B  \subseteq \bigcup_{k=0}^K Q_k.$ Let $T \in \calD_{m+1}$ such that $y \in T \subseteq Q_k$ for some $k$. For $\thetastop$ small enough, it follows from the estimate $|y-z| \leq 5\thetastop \rho^m$ that $z \in 2B_T \subseteq \Astop B_T$. So, by maximality of $m$, we must have $T \not\in S_{Q_k}$. In particular, there exists $R \in \min(S_{Q_k})$ such that $T \subseteq R$. This then satisfies $z \in \Gamma_0 \cap B \cap 2B_R,$ as required. 
\end{proof}

Below we introduce a slightly variant of the cubes $\Stop(N,k)$. Indeed, we need to identify those cubes in $\Stop(N,k)$ which have an element of $\Top(N,k)$ as a child. It is from nearby these cubes that we will build bridges between the approximating spaces.

\begin{notation}
	Let 
	\begin{align}
		\Stopprime(N-1,k) &\coloneqq \{Q \in \Stop(N-1,k) \colon {\rm Child}(Q) \cap \Top(N,k) \neq \emptyset\} \\
		&= \Stop(N-1,k) \cap \cubes_{<k}. 
	\end{align}
\end{notation}

Nearby cubes from $\Stopprime(N,k)$, it turns out that the distance functions $d_N$ and $d_{N,k}$ are comparable. 

\begin{lem}\label{l:N-and-k-com}
	Let $C > 1$, $Q \in \Stopprime(N,k)$ and $x \in CB_Q$. Then, for $\taustop = \taustop(C)$ sufficiently small, we have 
	\begin{align}\label{e:dN-dNk}
		d_N(x) = d_{N,k}(x).
	\end{align}
	In particular, if $R \in \Top(N,k)$ and $T \in \Stop'(N-1,k)$ are such that $Q \subseteq R$, then 
	\begin{align}\label{e:Q-very-small}
		\ell(Q) \lesssim \taustop \ell(R) \leq \taustop \ell(T).
	\end{align} 
\end{lem}

\begin{proof}
	Since $d_{N} \leq d_{N,k}$, the only inequality we need to prove is $d_N(x) \geq d_{N,k}(x)$. If $d_{N,k}(x) = d_N(x)$, then this is trivial. Suppose instead that $d_{N,k}(x) = 5\rho^{k}/\taustop$. Let $S \in \calT^N(k)$ such that $Q \in S(k)$. Since $Q \in \Stopprime(N,k)$ (so that $\ell(Q) \geq 5\rho^{k-1}$), we actually have $Q \in \min(S)$, which implies $Q \in \Stop(N)$. Thus, by Lemma \ref{l:stop-size},
	\begin{align}
		 d_N(x) \geq  \ell(Q)/(4\taustop) \geq 5\rho^k/\taustop = d_{N,k}(x),
	\end{align}
	as required for \eqref{e:dN-dNk}. Lemma \ref{l:stop-size-cutoff} and \eqref{e:dN-dNk} imply \eqref{e:Q-very-small} since 
	\begin{align}
		\ell(Q) \lesssim \taustop d_{N,k}(Q) \lesssim \taustop d_N(Q) \leq \taustop \ell(R) \leq \taustop \ell(T). 
	\end{align}
\end{proof}

Nearby cubes in $\Stopprime(N-1,k)$, we can also guarantee that the distance function is $d_{N,k}$ is smaller than $d_{N-1,k}$. 

\begin{lem}\label{l:monoton-in-N}
	Let $C > 1$, $N \in \N$ and $Q \in \Stopprime(N-1,k)$. Then, 
	\begin{align}
		d_{N,k}(x) < d_{N-1,k}(x) 
	\end{align}
	for all $x \in CB_Q$. 
\end{lem}

\begin{proof}
	By Lemma \ref{l:stop-size} and Lemma \ref{l:N-and-k-com}, we have 
	\begin{align}
		d_{N-1,k}(x) = d_N(x) \geq \ell(Q)/4. 
	\end{align}
	On the other hand, we first note that $5 \rho^k \leq \rho \ell(Q)$. To estimate $d_N(x)$ from above, consider an arbitrary cube $R \in \Top(N,k) \subseteq \Top(N)$ such that $R \subseteq Q$. Since $\{x\} \cup B_{R} \subseteq CB_Q$, we have 
	\begin{align}
		\taustop d_N(x) \leq \taustop \left( \ell(R) + \dist(x,R) \right) \leq \taustop (1+2C)\ell(Q) \leq \rho \ell(Q). 
	\end{align}
	Thus, we have $d_{N,k}(x) \leq \rho \ell(Q) < d_{N-1}(x)$, as required. 
\end{proof}

In fact, we can prove a more quantitative estimate regarding the distance between the distance functions at different generations. The following result ensure that the graphs we construct are well-separated after pushing them into a further dimension.

\begin{lem}\label{l:well-separated}
	Let $M,N \in \N$ such that $M > N$ and $C>1$. For each $Q \in \Stopprime(N,k)$ and $x,y \in CB_Q$ we have 
	\begin{align}\label{e:stop-prime}
		\taustop d_{N,k}(x) - 2\taustop d_{M,k}(y) \gtrsim \taustop d_{N,k}(x) \sim \ell(Q). 
	\end{align}
\end{lem}

\begin{proof}
	The right-most inequality is just Lemma \ref{l:stop-size-cutoff}, so the main point is to prove the left-most inequality. We start by giving a lower bound for $\taustop d_{N,k}(x)$. Indeed, since $\Stopprime(N,k) \subseteq \Stop(N)$ and $x \in CB_Q$, we have by Lemma \ref{l:stop-size} and Lemma \ref{l:stop-size-cutoff} that 
	\begin{align}\label{e:lower-x}
		\ell(Q)/4 \leq \taustop d_N(x) \leq \taustop d_{N,k}(x) .
	\end{align}
	We now show that 
	\begin{align}\label{e:upper-y}
		2\taustop d_{M,k}(y) \leq \ell(Q)/10. 
	\end{align}
	By Lemma \ref{l:monoton-in-N}, it suffices to check the above estimate with $N+1$ in place of $M$. Since $Q \in \Stopprime(N,k)$ we have $\ell(Q) \geq 5\rho^{\stopk-1}$ (otherwise then can be no child of $Q$ in $\Top(N+1,\stopk)$). Thus, 
	\begin{align}\label{e:rho}
		10\rho^\stopk \leq 2\rho \ell(Q) \leq \ell(Q)/10. 
	\end{align}
	On the other hand, if $R \in \Top(N+1,k) \subseteq \Top(N)$ is such that $R \subseteq Q$, we by the definition of $d_{N+1}$ that  
	\begin{align}\label{e:d}
		2\taustop d_{N+1}(y) \leq \taustop (\ell(R) + \dist(y,R)) \leq \taustop ( \rho \ell(Q) + C\rho^{-1}\ell(Q)) \leq \ell(Q)/10.
	\end{align}
	The inequality \eqref{e:rho} and \eqref{e:d} imply that $d_{N+1,k}(y) \leq \ell(Q)/10$, as required. Combing \eqref{e:lower-x} with \eqref{e:upper-y} and the right-most inequality in \eqref{e:stop-prime}, we get 
	\begin{align}
		\taustop d_{N,k}(x) - 2\taustop d_{M,k}(y) \geq \ell(Q)/10 \gtrsim \taustop d_{N,k}(x),
	\end{align}
	as required. 
\end{proof}

\subsection{Good and bad stopping-time regions}
In order to build the set $\envelope$ from Theorem \ref{t:main}, we will need to construct approximating sets for each $S \in \calT^N(\stopk)$. The construction will depend on whether $S$ lies in a certain \textit{good} class $\mathscr{G}^N_2(\stopk)$ of stopping-time regions (in which the top cube is near another \textit{flat} top cube which is not too small) or its complement $\mathscr{B}^N_2(\stopk)$ (where all flat top cubes are far away or very small). 
It will also be convenient to introduce here some slight variants of these classes $ \mathscr{G}_1^N(\stopk) \subseteq \mathscr{G}_2^N(\stopk)$ and $\mathscr{B}_1^N(\stopk) \supseteq \mathscr{B}_2^N(\stopk)$. 

\begin{constants}
	Let $\goodA \in \N$ to be chosen sufficiently large. We will suppose for the remainder of the section that $\taustop = \taustop(\goodA)$ is small enough such that 
	\begin{align}
		\Ccomp^{-1}\ell(R) \leq \ell(Q) \leq \Ccomp \ell(R)
	\end{align}
	for any $Q,R \in \Stop(N,\stopk)$ (resp. $\Top(N,\stopk)$) such that $2\goodA B_Q \cap 2\goodA B_R \neq \emptyset$, see Lemma \ref{l:comp-size}. 
\end{constants} 

\begin{constants}\label{constants-badk}
	Let $\badk = \badk(\Ccomp) \in \N$ be the small integer such that the following property holds: if $Q,R \in \Stop(N,\stopk)$ (resp. $\Top(N,\stopk)$) are such that $Q \in \cubes_{< k_2}$ and $2\goodA B_Q \cap 2\goodA B_R \neq \emptyset$, then $R \in \cubes_{\leq \stopk}$. 
\end{constants}

\begin{notation}
	Let 
	\begin{align}\label{e:good-bad}
		\calG(\stopk) \coloneqq \calG \cap \cubes_{<\badk} \quad \mbox{ and } \quad  \calB(\stopk) \coloneqq \calG(\stopk)^c = \calB \cup \cubes_{\geq \badk}. 
	\end{align}
\end{notation}

\begin{notation}
	For $i \in \{1,2\}$, let 
\begin{align}
	\hspace{3em}\mathscr{G}^N_i(\stopk) &\coloneqq \left\{
	S \in \calT^N(\stopk) \;\middle|\;
	\begin{aligned}
		& \mbox{there exists } T \in \calT^N(\stopk) \mbox{ satisfying }  \\
		& Q(T) \in \calG(\stopk) \mbox{ and } i\goodA B_{Q(T)} \cap i\goodA B_{Q(S)} \neq \emptyset 
	\end{aligned}
	\right\}; \label{e:good-stop} \\
	\hspace{3em}\mathscr{B}^N_i(\stopk) &\coloneqq \calT^N \setminus \mathscr{G}^N_i(\stopk). \label{e:bad-stop}
\end{align}
\end{notation}

We show below that any cube from $\Stop(N,k)$ which is sufficiently close to a good stopping-time region must be relatively small in comparision to the top cube. On the other hand, any cube in $\Stop(N,k)$ which is relatively close to some stopping-time region from $\mathscr{B}_1^N(k)$ must be relatively large compared to the top cube.

\begin{lem}\label{l:Stop-k}
	If $S \in \mathscr{G}^N_i(\stopk)$ and $Q \in \Stop(N,k_0)$ satisfy $3B_Q \cap 2\rho^{-1}B_{Q(S)} \neq \emptyset$, then $\ell(Q) \leq \rho \ell(Q(S))$. 
\end{lem}

\begin{proof}
	If $\ell(Q) = 5\rho^{\stopk}$ then the estimate follows from the fact that $\ell(Q(S)) > 5\rho^{\stopk}$. Suppose instead that $\ell(Q) > 5\rho^{\stopk}$, which implies $Q \in \Stopprime(N,\stopk)$. In this case, Lemma \ref{l:N-and-k-com} implies $\ell(Q) \lesssim \taustop \ell(Q(S))$, which gives the desired estimate for $\taustop$ small enough. 
\end{proof}

\begin{proof}
	If $Q \in \cubes_{\stopk}$ then it is clear that $\ell(Q) \leq \ell(Q(S))$ since $Q(S) \in \cubes_{\leq \stopk}$. Suppose instead that $Q \in \cubes_{\leq \stopk-1}$. In this case we know that $Q \in \Stop(N)$. By Lemma \ref{l:stop-size}, and the fact that $Q(S) \in \Layer(N)$ and $3B_Q \cap 2\rho^{-1}B_{Q(S)} \neq \emptyset$, we have 
	\begin{align}
		\ell(Q) \lesssim \taustop  d_N(Q) \leq \taustop  \left( \ell(Q(S)) + 3\ell(Q) + 2\rho^{-1}\ell(Q(S)) \right) .
	\end{align}
	Thus, for $\taustop $ sufficiently small (after taking $3\taustop  \ell(Q)$ to the left-hand side), the above inequality implies $\ell(Q) \leq \ell(Q(S))$. 
\end{proof}

\begin{lem}\label{l:stopped is pretty big}
	There is $\cLB \in (0,1)$ such that the following holds. Suppose $S \in \mathscr{B}_1^N(\stopk)$ and $Q \in \Stop(N,\stopk)$ are such that $Q \subseteq (\Mclose/2)B_{Q(S)}$. Then, $\ell(Q) \geq \cLB \taustop  \ell(Q(S))$. 
\end{lem}

\begin{proof}
	If $\ell(Q(S)) \leq 5\Ccomp \rho^{\badk}$ then the estimate is immediate since $\ell(Q) \geq 5\rho^{\stopk} \gtrsim  5\Ccomp \rho^{\badk}$. Assume instead that 
	\begin{align}\label{e:Q(S)-big}
		\ell(Q(S)) \geq 5\Ccomp \rho^{\badk-1}. 
	\end{align}
	Let $R \in \Stop(N)$ such that $R \subseteq Q \subseteq (\Mclose/2)B_{Q(S)}$. By Lemma \ref{l:stop-size}, it is sufficient to show  
	\begin{align}\label{e:comp-d_N}
		d_{N}(R) = \inf_{x \in R} d_N(x) \gtrsim  \ell(Q(S)) .
	\end{align}
	To get \eqref{e:comp-d_N}, fix an arbitrary $x \in R$, an arbitrary cube $P \in \Layer(N)$. It suffices to show 
	\begin{align}\label{e:suffices}
		\ell(P) + \dist(x,P) \gtrsim \ell(Q(S)).
	\end{align}
	Since $x \in R \subseteq (\Mclose/2)B_{Q(S)}$ and $P \subseteq B_{P}$, if $B_{P} \cap \Mclose B_{Q(S)} = \emptyset$, then \eqref{e:suffices} is obvious. Suppose instead that $B_{P} \cap (\Mclose/2) B_{Q(S)} \neq \emptyset$. Let $T \in \calT^N(\stopk)$ such that $P \in S_T'$ (recall Notation \ref{n:original-stopping-times}). This region satisfies
	\begin{align}
		B_{Q(T)} \cap (\Mclose/2)B_{Q(S)} \neq \emptyset. 
	\end{align}
	Thus, by Lemma \ref{l:comp-size-cutoff} and \eqref{e:Q(S)-big}, we have
	\begin{align}\label{e:not-in-B}
		\ell(Q(T)) \geq \Ccomp^{-1}\ell(Q(S)) \geq 5\rho^{\badk-1} \implies Q(T) \notin \cubes_{\geq \badk}. 
	\end{align}
	On the other case, by the definition of $\mathscr{B}_1^N(\stopk)$, 
	\begin{align}\label{e:not-in-G}
		Q(T) \notin \calG(\stopk) \implies Q(T) \in \calB(\stopk). 
	\end{align}
	Recalling \eqref{e:good-bad}, \eqref{e:not-in-B} and \eqref{e:not-in-G} imply $Q(T) \in \calB$. In particular, $S_T' = \{Q(T)\} = \{P\}$. Thus, by \eqref{e:not-in-B}, 
	\begin{align}
		\ell(P) = \ell(Q(T)) \gtrsim \ell(Q(S)),
	\end{align}
	as required. 
\end{proof}

\section{Approximations for individual cubes}\label{s:skeleta-per-cube}

Our goal for Sections \ref{s:skeleta-per-cube}--\ref{s:conse-layers} it to gather the necessary ingredients for defining a certain set $\envelope(\stopk)$. The set $\envelope$ from Theorem \ref{t:main} will then be obtained from $Z(\stopk)$ after sending $\stopk \to \infty$. For now, let us fix some $\stopk \in \N$. The basic building blocks of $\envelope_{\stopk}$ are element of $\{\grid{Q}\}_{Q \in \cubes}$, a collection $d$-dimensional boundaries of dyadic meshes associated to each Christ-David cube, which we describe in this section.

\begin{constants}\label{c:gridk}
	Let $\gridk \in \N$. Eventually, we will choose $\gridk = \gridk(\taustop)$ sufficiently small.
\end{constants}

\subsection{Definitions and key properties} For a cube $Q \in \cubes$, we construct a set $\grid{Q}$ as follows. First, let $\dyadic{Q}$ be the dyadic cube with $x_Q \in \dyadic{Q}$ and $\ell(\dyadic{Q})/2 < \ell(Q) \leq \ell(\dyadic{Q})$. Note, in particular, that 
\begin{align}\label{e:J-contains-ball}
	4\rho^{-1}B_{Q} \subseteq \lambdadyadic \dyadic{Q} \subseteq 2\sqrt{n}\lambda_1 B_Q, \mbox{ where } \lambdadyadic \coloneqq 4\rho^{-1}+1. 
\end{align} 
Note also (since $\lambda_1$ is odd) that, for any $k \in \N$ with $2^{-k} \leq \ell(\dyadic{Q})$, the cubes in $\Delta_k$ which intersect $\lambda_1 \dyadic{Q}$ also form a partition of it. Let 
\begin{align}
	\mathscr{I}_Q' &\coloneqq  \{I \in \Delta \colon I \subseteq \lambdadyadic \dyadic{Q} \mbox{ and } \ell(I) = 2^{-\gridk} \ell(\dyadic{Q}) \}, \\
	\mathscr{I}_Q & \coloneqq \mathscr{I}_Q' \cap \{I \in \Delta \colon \partial I \cap \partial (\lambdadyadic \dyadic{Q}) = \emptyset \},
\end{align}
and 
\begin{align}\label{e:grid-Q}
	\grid{Q} \coloneqq \bigcup_{I \in \mathscr{I}_Q} \partial_d I. 
\end{align}

It will be useful to be able to find a $d$-dimensional cubes in $G_Q$ be are close to the centres of sub-cubes of $Q$ and relatively large compared to $Q$. If we allow the cube to be relatively large with some dependency on $\taustop$, then this is just immediate from the definition. However, it will be useful to have the relative size independent of $\tau$.

\begin{lem}\label{l:subset-of-G_Q}
	Let $\kappa \in (0,1)$, $R \subseteq Q \in \cubes$ and $I \in \Delta$, and suppose that $\kappa \ell(Q) \leq \ell(I) \leq \ell(Q)$ and $x_R \in I$. Then, provided $\gridk = \gridk(\kappa)$ is sufficiently large, we have 
	\begin{align}
		\partial_d I \subseteq G_Q.
	\end{align}
\end{lem}

\begin{proof}
	Since $|x_R - x_Q| \leq \ell(Q) \leq \ell(\dyadic{Q})$ and $x_Q \in \dyadic{Q}$ there is a cube $J \in \Delta$ such that 
	\begin{align}
		x_R \in J \subseteq 3\dyadic{Q} \quad \mbox{ and } \quad \ell(J) = \ell(\dyadic{Q}). 
	\end{align}
	Since $x_R \in I$ and $\ell(I) \leq \ell(J)$, it follows from the properties of dyadic cubes that $I \subseteq J \subseteq 3\dyadic{Q}$. Since the cube $\mathscr{I}_Q$ partition $3\dyadic{Q}$, if $\gridk(\kappa)$ is chosen sufficiently large such that $2^{-\gridk}\ell(\dyadic{Q}) \leq \ell(I)$, then we can write 
	\begin{align}
		I = \bigcup_{\substack{K \in \mathscr{I}_Q \\ K \subseteq I}} K. 
	\end{align}
	This implies 
	\begin{align}
		\partial_d I = \bigcup_{\substack{K \in \mathscr{I}_Q \\ K \subseteq I}} \partial_d K \subseteq G_Q,
	\end{align}
	as required. 
\end{proof}

In a moment we will prove that the $G_Q$, and various configurations of them, support a Poincar\'e inequality. The following lemma will be useful for this. It shows that the grids corresponding to cubes with comparable sizes are either very disjoint or they have large intersection. Moreover, the gluing metric between grids with large intersection is bi-Lipschitz equivalent (with respect to the identity) to the Euclidean metric. 

\begin{lem}\label{l:face-glue}
	Let $C \geq 1$ and suppose $Q,R \in \cubes$ are such that $C^{-1}  \leq \ell(R)/\ell(Q) \leq C$. If $\gridk$ large enough (depending only on $C$) then we have either $\lambdadyadic\dyadic{Q} \cap \lambdadyadic\dyadic{R} = \emptyset$, in which case
	\begin{align}
		\dist(\grid{Q},\grid{R}) \geq 2^{-\gridk} \max\{\ell(Q),\ell(R)\},
	\end{align}
	or $\lambdadyadic\dyadic{Q} \cap \lambdadyadic\dyadic{R} \neq \emptyset$, in which $\grid{Q} \cap \grid{R}$ is non-empty, lower Ahlfors $d$-regular and 
	\begin{align}\label{e:gluing-metric}
		|x-y| \sim \inf_{z \in \grid{Q} \cap \grid{R}} |x-z| + |y-z| 
	\end{align}
	for all $x \in \grid{Q}$ and $y \in \grid{R}$. 
\end{lem}

\begin{proof}
	Since $\grid{Q} \subseteq \lambdadyadic\dyadic{Q}$ and $\dist(\grid{Q},\partial(\lambdadyadic\dyadic{Q})) \geq 2^{-\gridk} \ell(Q)$ (with the same properties holding for the cube $R$), if $\lambdadyadic\dyadic{Q} \cap \lambdadyadic\dyadic{R} = \emptyset$ then the second conclusion holds. So, let us suppose that 
	\begin{align}\label{e:non-empty-intersection}
		\lambdadyadic\dyadic{Q} \cap \lambdadyadic \dyadic{R} \neq \emptyset,
	\end{align}
	and let us assume (without loss of generality) that $\ell(Q) \leq \ell(R) \leq C \ell(Q)$. \\
	
	To see $G_Q \cap G_R$ is non-empty, note to begin that there exists a dyadic cube $I \subseteq \lambdadyadic\dyadic{Q} \cap \lambdadyadic\dyadic{R}$ with side length $\ell(I) = \ell(\dyadic{Q})$. Consider the cube $I'$ such that $x_I \in I'$ and $\ell(I') = 2^{-\gridk}\ell(\dyadic{R}) \leq C 2^{-\gridk}\ell(\dyadic{Q})$. As long as $\gridk = \gridk(C)$ is sufficiently small, it follows that 
	\begin{align}\label{e:far-from-boundary}
		I' \cap \left(\partial(\lambdadyadic\dyadic{Q}) \cup \partial(\lambdadyadic\dyadic{R})\right) = \emptyset,
	\end{align}
	i.e., $I' \in \mathscr{I}_R$. Consider also an arbitrary dyadic cube $I'' \subseteq I'$ with side length $\ell(I'') = 2^{-\gridk}\ell(\dyadic{Q})$ such that $\partial_d I' \cap \partial_d I'' \neq \emptyset$. Since $I'' \subseteq I'$, \eqref{e:far-from-boundary} also holds with $I''$ in place of $I'$. Thus, $I'' \in \mathscr{I}_Q$. Since the $d$-dimensional boundaries of $I'$ and $I''$ have non-empty intersection, it now follows that $\grid{Q} \cap \grid{R} \neq \emptyset$. \\

	Our next goal is to show that $\grid{Q} \cap \grid{R}$ is lower Ahlfors $d$-regular. We start by proving the following. \\
	
	\noindent \textbf{Claim 1.} For each $x \in \grid{Q} \cap \grid{R}$ there are cubes $I \in \mathscr{I}_Q$ and $J \in \mathscr{I}_R$ such that $I \subseteq J$ and $x \in \partial_d I \cap \partial_d J$.   \\
	
	Let $T_Q = \bigcup_{I \in \mathscr{I}_Q} I$ and $T_R = \bigcup_{I \in \mathscr{I}_R} I$. Since $T_Q$ and $T_R$ are cubes in $\R^n$ (so that intersections and closures commute), we have 
	\begin{align}
		\grid{Q} \cap \grid{R} \subseteq \overline{T_Q} \cap \overline{T_R} = \overline{T_Q \cap T_R}.
	\end{align}
	Since every cube in $\mathscr{I}_Q$ is smaller than every cube in $\mathscr{I}_R$, we can write 
	\begin{align}
		\overline{T_Q \cap T_R} =\overline{\bigcup_{\substack{I \in \mathscr{I}_Q \\ I \subseteq T_R}} I} = \bigcup_{\substack{I \in \mathscr{I}_Q \\ I \subseteq T_R}} \overline{I}.
	\end{align}
	In particular, for each $x \in \grid{Q} \cap \grid{R}$ there are cubes $I \in \mathscr{I}_Q$ and $J \in \mathscr{I}_R$ such that $x \in \overline{I} \subseteq \overline{J}$ (so that, $I \subseteq J$). Since $x \in \grid{Q}$ (i.e., it is in the $d$-dimensional boundary of some cube $I'$ with $\ell(I') = \ell(I)$), we conclude that $x \in \partial_d I$. By similar reasoning, we also conclude that $x \in \partial_d J$, as required. \\
	
	Fix some $x \in \grid{Q} \cap \grid{R}$ and $0 < r < \ell(R)$. By Claim 1, there are cubes $I \in \mathscr{I}_Q$ and $J \in \mathscr{I}_R$ such that $I \subseteq J$ and $x \in \partial_d I \cap \partial_d J$. By Lemma \ref{l:sub-face}, there are $d$-dimensional sub-faces $F_I \subseteq \partial_d I$ and $F_J \subseteq \partial_d J$ such that $x \in F_I \subseteq F_J$. Thus, by the definition of $\grid{Q}$ and $\grid{R}$, we have 
	\begin{align}
		F_{I} \subseteq \grid{Q} \cap \grid{R}. 
	\end{align} 
	Since $\ell(F_{I}) = \ell(I) \sim \ell(R)$ and $F_{I}$ is lower Ahlfors $d$-regular, it follows that 
	\begin{align}
		\calH^d((\grid{Q} \cap \grid{R}) \cap B(x,r)) \geq \calH^d(F_{I} \cap B(x,r)) \gtrsim r^d,
	\end{align}
	as required. \\
	
	Finally, let us check \eqref{e:gluing-metric}. Fix some $x \in \grid{Q}$ and $y \in \grid{R}$. Let $I \in \mathscr{I}_Q$ such that $x \in \partial_d I$ and let $J \in \mathscr{I}_R$ such that $y \in \partial_d J$. If $\partial_d I \cap \partial_d K' = \emptyset$ for all $K' \in \mathscr{I}_R$, then $|x-y| \geq \ell(I)$. Thus, by taking an arbitrary $z \in \grid{Q} \cap \grid{R}$ and using the fact that $\diam(\grid{Q} \cup \grid{R}) \lesssim \ell(R) \lesssim \ell(I)$, we get 
	\begin{align}
		|x-z| + |z-y| \leq 2|x-z| + |x-y| \lesssim \ell(I) + |x-y| \lesssim |x-y|. 
	\end{align}
	Thus, we can assume that there exists some $K' \in \mathscr{I}_R$ such that $\partial_d I \cap \partial_d K' \neq \emptyset$. Applying Claim 1, we find another $K \in \mathscr{I}_R$ such that $I \subseteq K$. The estimate \eqref{e:gluing-metric} now follows from Lemma \ref{l:reverse-triangle} (recalling that $\partial_d I \cap \partial_d K \subseteq \grid{Q} \cap \grid{R}$). 
\end{proof}

\subsection{Poincar\'e inequalities for the $G_Q$} We now collect various properties of the set we constructed above in terms of Poincar\'e inequalities. We start with the most basic, that each $G_Q$ supports a Poincar\'e inequality.

\begin{lem}\label{l:grid-PI}
	The set $\grid{Q}$ is Ahlfors $d$-regular and supports a $(1,1)$-Poincar\'e inequality with constants depending on $\gridk$ and $n$.
\end{lem}

\begin{proof}
 We first observe that for any cube $I\subset \mathbb{R}^n$ it holds that $\partial_d I$ is Ahlfors $d$-regular with regularity constant  depending only on $n$. Up to rescaling $I$, we can assume that $\ell(I)=1.$  The $d$-regularity follows e.g.\ from Lemma \ref{l:ADR-gluing}, indeed $\partial_d I$ is the union of less thatn $3^n$ closed $d$-dimensional cubes $Q_i$ of side-length one and such that $F_i\coloneqq Q_i\cap Q_{i+1}\neq \emptyset$.

   Recall that $G_Q$ is the union of the $d$-skeleton of all the dyadic cubes of side-length $2^{-\gridk}\ell(J_Q)$ contained in the dyadic cube $\lambda_1J_Q\in \Delta$ minus the outermost layer (recall that $\lambda_1$ is odd). Hence $G_Q=\cup_{i=1}^N\partial_d I_i$, with $I_i$ disjoint dyadic cubes such that $\ell(I_i)=2^{-\gridk}\ell(J_Q)$  with $N\le (2^{\gridk}\lambda_1)^n\lesssim_n (2^{\gridk})^n$ (in fact $N= L^n$ with $L=2^{\gridk}\lambda_1-2$) because $\lambda_1=4\rho^{-1}+1=4001$. Up to a rescaling, we can assume that $\ell(J_Q)=2^\gridk$ and so $G_Q$ coincides with the $d$-skeleton of a cubical one-by-one grid in $\mathbb{R}^n$ of side-length $L.$ In other words, up to translation, $G_Q=\{(x_1,\dots,x_n)\in [0,L]^n \ : \  \text{$x_i$ is a non-negative integer for at least $n-d$ choices of $i$}\}$.   We can also choose an enumeration so that  $\partial_d I_i\cap \partial_d I_{i+1}$ is non-empty. Applying Lemma \ref{l:ADR-gluing} $N-1$ times follows that $G_Q$ is  Ahlfors $d$-regular with regularity constant depending only on $\gridk$.

    To show the validity of the Poincar\'e inequality we aim to apply the results  in \cite{PivarskiSaloffeCoste}. We need to check that $G_Q$ satisfies (1)-(5) in \cite[Assumtions 1.10]{PivarskiSaloffeCoste}. Note that in the terminology of \cite{PivarskiSaloffeCoste} the set $G_Q$ equipped with the intrinsic length distance $d_l$ is a connected Euclidean homogeneous polyhedral complex of dimension $d.$ 

    Since $G_Q$ is homogeneous (the maximal polyhedra have have dimension $d$), to verify (1) we need only to check that if $F,E$ are two closed $d$-dimensional faces of $G_Q$ such that $F\cap E$ is a $k$-dimensional face $P$, $k\le d-1$, then there is a chain of closed $d$-dimensional faces $E=G_1,\dots,G_m=F$ such that $G_i\cap G_{i+1}$ is a $(d-1)$-dimensional face and $P\subset G_i$. We have that $E=\prod_{i=1}^n E_i$  and $F=\prod_{i=1}^n F_i$, with $E_i,F_i$ either closed unit intervals or singletons and $E_i\cap F_i\neq \emptyset$ for all $i.$ Set $G_1\coloneqq E$ and proceed inductively as follows. Take the first $i$ such that $E_i\neq F_i$. If both $E_i$ and $F_i$ are intervals  replace $E_i$ in $E$ by $F_i$. If  $E_i$ is an interval and $F_i$ is a point, since the  $E_i$'s and $F_i$'s are  intervals exactly $d$ times, there must be $j$ such that $E_j$ is a point while $F_j$ is not. Then  we  replace $E_i$ with $F_i$ and $E_j$ by $F_j$. If  $E_i$ is a point and $F_i$ is an interval we do the same but reversed. This produces a  $d$-dimensional face (still in $G_Q$) $G_2$ such that $E\cap G_2$ is a $(d-1)$-dimensional face and $P\subset G_2$, indeed $P\subset \prod_{i=1}^n E_i\cap F_i\subset G_2$ by construction because $G_2$ is obtained by swapping some of the $E_i$'s with some of the $F_i$'s. After finitely many step we end up with $G_m=F$ as desired.
    
    Assumption (2) is trivially true  with $\Lambda=1$ because each polyhedron is equipped with the Euclidean metric.
    
    For (3) we have that for each $j< k\le d$ and any $j$-dimensional face $F$ of $G_Q$, the maximal number of $k$-dimensional faces that contain $F$ is clearly bounded above by $M=3^n$. Indeed $F=\prod_{i=1}^n F_i$ where each $F_i$  is either a unit intervals with integer endpoints or an integer. If a $k$-dimensional face $E$  of $G_Q$ contains $F$ then $E=\prod_{i=1}^n E_i$, where $E_i$ are again either a unit intervals with integer endpoints or an integer and $F_i\subset E_i$. However, independently  of $F_i$ there are clearly at most three choices of  $E_i.$

    For (4) we can immediately check that the distance between two disjoint $k$-dimensional faces of $G_Q$ is at least $\ell=1$.

    Assumption (5) requires a \textit{uniform solid angle lower bound $\alpha$}. We recall the meaning of this.  Take any $0\le k\le d$ and any ball $B\subset G_Q$  with respect to the length distance with center $x_B$ lying on a $k$-dimensional  face $F$ and radius $r<d_l(x_B,G_Q^{k-1})$ (where  $G_Q^{k-1}$ is the $(k-1)$-skeleton of $G_Q$), assuming simply $r<1$ if $k=0$. A wedge  $W$ is the closure of a  connected component  of $B\setminus G_Q^{d-1}$. Having a uniform solid angle lower bound $\alpha$ means that $\mathcal H^d(W)/(\omega_d r^d) \ge \alpha$ for all wedges $W$ (note that this ratio is always $\le 1$). In our case, for any $B$ as above, we have that $B$ is isometric to the union of copies of $B^{\mathbb R^d}(0,r)\cap \mathbb R^k\times [0,\infty)^{d-k}$. Therefore each wedge $W$ is a copy of $B^{\mathbb R^d}(0,r)\cap \mathbb R^k\times [0,\infty)^{d-k} $ and so $\mathcal H^d(W)/(\omega_d r^d) \ge 2^{-(d-k)}\ge 2^{-d}\ge 2^{-n}$. Hence we can choose $\alpha=2^{-n}.$
   
    Therefore we can apply  \cite[Corollary 2.10]{PivarskiSaloffeCoste} choosing $R_0'\coloneqq 2\diam_{d_l}(G_Q)$ (see also \cite[Theorem 2.7]{PivarskiSaloffeCoste}) to deduce that $G_Q$, \textit{endowed with the length metric $d_l$}, supports a $(1,1)$-Poincar\'e inequality with constants depending only on $\gridk,d$ and $n$. Indeed note that their notion of gradient for a Lipschitz function $f$ in $G_Q$  satisfies $|\nabla f|\le \Lip(f)$ $\mathcal H^d$-a.e. To pass to $G_Q$ with the Euclidean distance $d_e$ we only need to observe that $d_e\le d_l\lesssim_{n} d_e.$ The first inequality is obvious. For the second fix $x=(x_1,\dots,x_n),y=(y_1,\dots,y_n)\in G_Q$.  Recall that at least $n-d$ entries for both $x$ and $y$ must be integers. Suppose first that $|x-y|< 1$. Then if $x_i$ and $y_i$ are both integers they must coincide. We join $x$ and $y$ by   
     the curve which is the union of the segments formed by  subsequent linear interpolation from $x_i$ to $y_i$, starting first from all the non-integer entries $x_i$. This  is contained in $G_Q$, since at all times there are at least $n-d$-integer coordinates. Moreover this curve has length $\sum_{i=1}^n|x_i-y_i|\le \sqrt n |x-y|.$ Suppose now that $|x-y|\ge  1$. Since both $x$ and $y$ lie in  $d$-dimensional faces of $G_Q$ we can take vertices of such faces $\bar x,\bar y$ having all integer coordinates such that $|x-\bar x|,|y-\bar y|\le \sqrt n$. Moreover the segments $x,\bar x$ and $y,\bar y$ are contained in $G_Q$ (because $d$-dimensional faces are convex). As above we can build a curve from $\bar x$ to $\bar y$ contained in $G_Q$ of  length at most $\sqrt{n}|\bar x-\bar y|$ which,  joined with the two segments $x,\bar x$ and $y,\bar y$, produces a curve from $x$ to $y$ of length at most $\sqrt{n}|\bar x-\bar y|+2\sqrt n\le (3\sqrt n+2n)|x-y|$, as desired.
\end{proof}

More generally, we would like to show that the union of sets of this form support a Poincar\'e inequality (under certain hypotheses). This will be our goal for the remainder of the section. In Lemma \ref{l:local-PI} below we show this to be true provided the union is connected (however, we allow the constant to depend on the number of elements in the union).

\begin{lem}\label{l:local-PI}
	Let $C \geq 1$, $\calQ$ be a finite collection of cubes such that $C \leq \ell(Q)/\ell(R) \leq V$ for all $Q,R \in \calQ$, and suppose that $G \coloneqq \bigcup_{Q \in \calQ} \grid{Q}$ is connected. If $\gridk$ large enough (depending only on $C$) then $G$ is Ahlfors $d$-regular and supports a $(1,1)$-Poincar\'e inequality with constant depending only on $\gridk,n$ and $|\calQ|$.
\end{lem}

\begin{proof}
	Ahlfors regularity follows directly from Lemma \ref{l:ADR-gluing}. Thus, our main goal is to prove prove a Poincar\'e inequality. Let $K \coloneqq |\calQ|$ and let $\gridk = \gridk(C)$ be as in Lemma \ref{l:face-glue}. Inductively, we will define a sequence of cubes $\calQ = \{Q_i\}_{i=1}^K$ such that 
	\begin{align}\label{e:PI}
		G_k \coloneqq \bigcup_{i=1}^k \grid{Q_k} \mbox{ supports a } (1,1)\mbox{-Poincar\'e inequality }
	\end{align} 
	for each $k \in \{1,\dots,K\}$ with constant depending on $k$. Indeed, to begin with we let $Q_1$ be an arbitrary cube in $\calQ$. Then, \eqref{e:PI} follows immediately for $G_1$ by Lemma \ref{l:grid-PI}. Inductively, suppose that we have constructed cubes $\{Q_i\}_{i=1}^k$ with $k < K$ such that $G_k$ satisfies \eqref{e:PI}. Since $G$ is connected, there exists some cube in $\calQ$, which we label $Q_{k+1}$, such that $G_k \cap \grid{Q_{k+1}} \neq \emptyset$. Let 
	\begin{align}
		\calI_1 \coloneqq \{i \in \{1,\dots,k\} \colon \grid{Q_i} \cap \grid{Q_k} \neq \emptyset\} \quad \mbox{and} \quad \calI_2 \coloneqq \{1,\dots,k\} \setminus \calI_1. 
	\end{align}
	By Lemma \ref{l:face-glue}, we know that 
	\begin{align}
		A \coloneqq G_k \cap \grid{Q_{k+1}} = \bigcup_{i \in \calI_1} \grid{Q_i} \cap \grid{Q_{k+1}}
	\end{align}
	is lower Ahlfors $d$-regular. By Lemma \ref{l:gluing}, the only thing left to show is that 
	\begin{align}\label{e:T'}
		|x-y| \sim \inf_{z \in A} |x-z| + |y-z| 
	\end{align} 
	for any pair of points $x \in G_k$ and $y \in \grid{Q_{k+1}}$, which we now fix.  \\
	
	Let $i \in \{1,\dots,k\}$ such that $x \in \grid{Q_i}$. If $i \in \calI_1$, then \eqref{e:T'} follows immediately from Lemma \ref{l:face-glue} since $A \supseteq \grid{Q_i} \cap \grid{Q_{k+1}}$. Suppose now that $i \in \calI_2$. In this case, it follows from Lemma \ref{l:face-glue} that $|x-y| \gtrsim \ell(Q_{k+1})$. Furthermore, since $A \cap \grid{Q_{k+1}} \neq \emptyset$ and $\diam(\grid{Q_{k+1}}) \lesssim \ell(Q_{k+1})$, there exists $z \in A$ such that $|y-z| \lesssim \ell(Q_{k+1})$. Combing these two facts, we get 
	\begin{align}
		|x-z| + |y-z| \leq |x-y| + 2|y-z| \lesssim |x-y| + \ell(Q_{k+1}) \lesssim |x-y|,
	\end{align} 
	as required. 
\end{proof}

We now show how to remove this dependency on the number of elements in the union provided they lie close to a graph generated by the class $\mathscr{G}_2^N(\stopk)$ of good stopping-time regions from \eqref{e:good-stop}. In order to fix constants, we introduce the following.

\begin{constants}
	Let $\Afirst \coloneqq 3000\sqrt{n}\lambdadyadic \Ccomp$ and $\Cprop = (\cLB\taustop)^{-1}$. 
\end{constants}

\begin{prop}\label{l:sigma-good}
	Let 
	\begin{align}
		k \in \N, \quad A \in [\Afirst, 2\Afirst^2] \quad \mbox{ and } \quad C \in [1,\Cprop]. 
	\end{align}
	Let $S \in \mathscr{G}^N_2(k)$ and $\{R_j\}_{j \in J} \subseteq \bigcup_{T \in \calT^N(\stopk)} T(\stopk)$ be a disjoint collection of cubes such that 
	\begin{enumerate}
		\item $E \cap AB_{Q(S)} \subseteq \bigcup_{j \in J} R_j$; \label{p:R-covers-E}
		\item $E \cap 3B_{R_j} \cap A B_{Q(S)} \neq \emptyset$ for all $j \in J$; \label{p:R-intersects-E}
		\item $C^{-1} \leq \ell(R_i)/\ell(R_j) \leq C$ for all $i,j \in J$ such that $AB_{R_i} \cap AB_{R_j} \neq \emptyset$. \label{p:R-comp}
	\end{enumerate}
	Then, provided 
	\begin{align}
		\taustop = \taustop(\Ccomp,\Mclose),  \ \thetastop = \thetastop(\taustop) \mbox{ and } \alphastop
	\end{align}
	are sufficiently small and
	\begin{align}
		\Mclose = \Mclose(\Afirst ), \  \Astop = \Astop(\Ccomp,\Mclose) \mbox{ and } \gridk = \gridk(\taustop)
	\end{align} 
	are sufficiently large, we have that $\bigcup_{j \in J} G_{R_j}$ is Ahlfors $d$-regular with constants depending on $d,n,C_0,C$,  and supports a $(1,1)$-Poincar\'e inequality with constant $\CPI = \CPI(d,n,A,C,C_{\rm grid})$ and $\lambdaPI = \lambdaPI(n , C, \lambda_{\rm grid} )$.
\end{prop}

We will prove Proposition \ref{l:sigma-good} by applying Proposition \ref{p:local-to-global-prime} with the collection of balls $\{B_j\}_{j \in J} \coloneqq \{B_{R_j}\}_{j \in J}$, the space
\begin{align}
	Y = \bigcup_{j \in J} G_{R_j},
\end{align}
and $X$ a (truncated) Lipschitz graph associated to the nearby stopping-time with flat root. To be more precise regarding $X$, let $T \in \calT^N(\stopk)$ such that $Q(T) \in \calG(\stopk)$ and $2\Mclose B_{Q(T)} \cap 2\Mclose B_{Q(S)} \neq \emptyset$ (which exists by the definition of $\mathscr{G}^N_2(\stopk)$, recall \eqref{e:good-stop}). Let 
\begin{align}
	\{Q_i\}_{i \in I} \coloneqq \{ Q \in \Top(N,\stopk) \colon Q \cap 5\Ccomp \Mclose B_{Q(T)} \neq \emptyset \},
\end{align} 
where $\Ccomp$ is as in Lemma \ref{l:comp-size-cutoff}. Assume that $\Astop = \Astop(\Ccomp,\Mclose)$ is large enough and $\taustop = \taustop(\Ccomp, \Mclose)$ is small enough so that we can apply Lemma \ref{l:lip-graphs} with the collection $\{Q_i\}_{i \in I}$ and $Q(T)$. Let $\Gamma$ be the resulting graph, then set 
\begin{align}
	X \coloneqq \Gamma \cap \pi_T^{-1}\left(P_{Q(T)} \cap AB_{Q(S)} \right) .
\end{align}

It only remains to check the hypotheses of Proposition \ref{p:local-to-global-prime} with this data. We will check it with  
\begin{align}\label{e:C_1}
		C_1 \coloneqq 6\max\{\sqrt{n}\lambdadyadic,\Ccomp\}
\end{align} 
(the remaining constants shall be found shortly). This is done in Lemmas \ref{l:1&4}--\ref{l:union-of-the-Gs} below, recall that some of these hypotheses are stated in Lemma \ref{l:gluing-ADR}. 

\begin{lem}\label{l:1&4}
	Proposition \ref{p:local-to-global-prime} \eqref{prop:separated} and \eqref{prop:comp-size} with constants $\csep = \innercube$, $C_3 = A$ and $C_4  = C$, and Lemma \ref{l:gluing-ADR} \eqref{prop:bounded-overlap} and \eqref{prop:cover} hold with constants $C_1$ (as above) and $C_2$ depending on $d,C,C_0,C_1$ and $c_0$.    
\end{lem}

\begin{proof}
	Proposition \ref{p:local-to-global-prime} \eqref{prop:separated} holds with constant $\csep = \innercube$ since $\{R_j\}_{j \in J}$ is a disjoint collection of cubes. Proposition \ref{p:local-to-global-prime} \eqref{prop:comp-size} holds with $C_3 = A$ and $C_4 = C$ by assumption Lemma \ref{l:sigma-good} \eqref{p:R-comp}. Lemma \ref{l:gluing-ADR} \eqref{prop:cover} follows from Lemma \ref{l:contained-in-ball}. Lemma \ref{l:gluing-ADR} \eqref{prop:bounded-overlap} follows from Proposition \ref{p:local-to-global-prime} \eqref{prop:separated} and \eqref{prop:comp-size} by a volume argument. 
\end{proof}

\begin{lem}
	Lemma \ref{l:gluing-ADR} \eqref{prop:close} holds with $C_1$.
\end{lem}

\begin{proof}
	Fix some $j \in J$ and $x \in E \cap 3B_{R_j} \cap AB_{Q(S)}$, which exists by assumption (Lemma \ref{l:sigma-good} \eqref{p:R-intersects-E}). By Lemma \ref{l:lip-graphs}, there exists some $y \in \Gamma \cap 4B_{R_j}$ such that $|\pi_T(y) - \pi_T(x_{Q(S)})| \leq \ell(R_j)$. By sliding along $\Gamma$, we can find a further point $z \in \Gamma \cap 5B_{R_j}$ such that $\pi_T(z) \in P_{Q(T)} \cap AB_{Q(S)}$ i.e., $z \in X$.
\end{proof}

\begin{lem}
	Lemma \ref{l:gluing-ADR} \eqref{prop:cover-X} holds with constant $C_1$.  
\end{lem}

\begin{proof}
	Observe to begin with that
	\begin{align}\label{e:inclusion-prime}
		E \cap 2AB_{Q(S)} \subseteq E \cap 5\Ccomp \Mclose B_{Q(T)} \subseteq \bigcup_{i \in I} Q_i. 
	\end{align}
	Indeed, the first inclusion uses Lemma \ref{l:comp-size-cutoff} along with the fact that $2\Mclose B_{Q(S)} \cap 2\Mclose B_{Q(T)} \neq \emptyset$. The second inclusion follows since $E \cap 5\Ccomp \Mclose B_{Q(T)}$ is covered by the set of $Q \in \Top(N)$ such that $Q \cap 5\Ccomp \Mclose B_{Q(T)} \neq \emptyset$. Since $Q(T) \in \mathcal{G}(\stopk)$, our choice of $\gridk$ (see Constants \ref{c:gridk}) implies that any such cube lies in $\cubes_{\leq \stopk}$ i.e., $Q \in \{Q_i\}_{i \in I}$. 
	
	By \eqref{e:inclusion-prime} and Lemma \ref{l:lip-graphs-subset}, we have 
	\begin{align}
		X \subseteq  \bigcup_{i \in I} \bigcup_{\substack{ R \in \min(S_{Q_i}(\stopk)) \\ \Gamma \cap AB_{Q(S)} \cap 2B_R \neq \emptyset }} 2B_R.
	\end{align}
	Fix some $R \in \min(S_{Q_i}(\stopk))$ such that $\Gamma \cap AB_{Q(S)} \cap 2B_R \neq \emptyset$. We claim that there exists $j \in J$ such that 
	\begin{align}\label{e:R-R_j}
		2B_R \subseteq 6\Ccomp B_{R_j}, 
	\end{align}
	which is sufficient for Lemma \ref{l:gluing-ADR} \eqref{prop:cover-X}. Indeed, by travelling along $\Gamma$ slightly, we can find a point $y \in \Gamma \cap 3B_R$ such that $B(y,5 \thetastop \ell(R)) \subseteq AB_{Q(S)}$. By Lemma \ref{l:lip-graphs} there exists a point $z \in E \cap AB_{Q(S)} \cap 4B_R$ such that $|y-z| \leq \thetastop \ell(R)$. By assumption (Lemma \ref{l:sigma-good} \eqref{p:R-covers-E}), there exists $j \in J$ such that $z \in R_j \cap 4B_R$. It now follows from Lemma \ref{l:comp-size-cutoff} that $\Ccomp \ell(R_j) \geq \ell(R)$. Then, \eqref{e:R-R_j} follows from the triangle inequality. 
\end{proof}

It only remains to check that Proposition \ref{p:local-to-global-prime} \eqref{prop:small-balls} holds with constant $C_5 = 50C_1$. Note that $C_5 \leq \Afirst \leq A$.\\

Fix an index $j \in J$ and a pair of points $x,y \in C_5B_{j}$. Let $Q,R \in \{R_i\}_{i \in J}$ such that $x \in \grid{Q}$ and $y \in \grid{R}$. Before proving Proposition \ref{p:local-to-global-prime} \eqref{prop:small-balls}, we will prove the following useful lemma.

\begin{lem}\label{l:chain}
	If $\Astop = \Astop(C_5)$ is sufficiently large, and $\thetastop = \thetastop(\taustop)$ is sufficiently small, there is a finite sequence of cubes $\{T_k\}_{k=0}^K$, with $K \lesssim_{C,\Cmax} 1$, such that $T_0 = Q$, $T_K = R$ and $6B_{T_{k}} \cap 6B_{T_{k+1} }\neq \emptyset$ for each $k \in \{0,\dots,K-1\}$.
\end{lem}

\begin{proof}
	Since $\grid{T} \subseteq 2\sqrt{n} \lambdadyadic B_T$ for each $T \in \cubes$, we have 
	\begin{align}\label{e:x-y}
		x \in 2\sqrt{n}\lambdadyadic B_Q \cap C_5 B_{R_j} \quad \mbox{ and } \quad y \in 2\sqrt{n} \lambdadyadic B_R \cap C_5 B_{R_j},
	\end{align}
	Since $A \geq C_5$, it follows by assumption (Lemma \ref{l:sigma-good} \eqref{p:R-comp}) that 
	\begin{align}\label{e:Q,R,R_j}
		\ell(Q)  \sim_{C} \ell(R_j) \sim_{C} \ell(R) .
	\end{align}
	Let 
	\begin{align}
		Q^* \in \{Q,R,R_j\} \mbox{ have the maximal side-length.} 
	\end{align}
	One can check (from \eqref{e:x-y} and the triangle inequality) that 
	\begin{align}\label{e:Q^*}
		5B_Q \cup 5B_R \subseteq \lambda B_{Q^*},
	\end{align}
	where 
	\begin{align}
		\lambda = 5 + 4\sqrt{n}\lambdadyadic + 600\max\{\sqrt{n}\lambdadyadic,\Ccomp\} \leq 100\sqrt{n}\lambdadyadic \Ccomp \leq \Afirst/3. 
	\end{align}

	By arguing as in the proof of Lemma \ref{l:gluing-ADR} \eqref{prop:close}, there are points $z_Q \in X \cap 5B_Q$ and $z_R \in X \cap 5B_R$. Let $L$ be the line segment from $z_Q$ to $z_R$ such that $\pi_T(\Gamma_L) = \pi_T(L)$. Let $\beta > 0$ (to be chosen small enough, depending on $C$) and $\{w_k\}_{k=0}^K$ be a collection of points in $L$ such that $w_0 = \pi_T(z_Q)$, $w_K = \pi_T(z_R)$ and
	\begin{align}\label{e:w_i}
		\beta \ell(Q^*) / 2 \leq |w_{i+1} - w_i| \leq \beta \ell(Q^*).
	\end{align}
	By \eqref{e:Q^*} and \eqref{e:w_i}, we have 
	\begin{align}
		K \lesssim_{\Cmax,\beta} 1.
	\end{align}
	
	Let $\Gamma_L \subseteq \Gamma$ such that $\pi_T(\Gamma_L) = \pi_T(L)$, and let $\{y_k\}_{k=0}^K$ be the corresponding points in $\Gamma_L$ i.e., $\pi_T(y_k) = w_k$. By convexity, it follows from \eqref{e:Q^*} that $L \subseteq \lambda B_{Q^*}$. Since $\Gamma$ is a $C\alphastop$-Lipschitz graph, we have 
	\begin{align}\label{e:dist-L}
		\dist(y,L) \leq (1+C\alphastop)|z_Q - z_R| \mbox{ for all } y \in \Gamma_L.
	\end{align}
	These two facts imply 
	\begin{align}
		\{y_k\}_{k=0}^K \subseteq 2\lambda B_{Q^*}. 
	\end{align}
	
	If $\Astop = \Astop(\lambda)$ is large enough, it follows from Lemma \ref{l:lip-graphs} that for each $y_i$ there is $z_i \in E \cap AB_{Q(S)} \cap 3\lambda B_{Q^*}$ such that $|y_i - z_i| \leq \thetastop \ell(Q^*)$. Let $T_0 =Q$ and $T_K = R$. For the remaining indices, let $T_i \in \{R_j\}_{j \in J}$ such that $z_i \in T_i$. Since $z_i \in 3\lambda B_{Q^*}$, it follows by assumption that $\ell(T_i) \sim_{C} \ell(Q^*)$ (recall that $A \geq \Afirst \geq 3\lambda$). Thus, using \eqref{e:w_i} and the fact that $X$ is the (truncated) graph of an $C\alphastop$-Lipschitz function, we have 
	\begin{align}
		|z_{i+1} - z_i| &\leq 2\thetastop \ell(Q^*) + |y_{i+1}-y_i| \leq 2\thetastop \ell(Q^*) + 2|w_{i+1} - w_i| \\
		&\lesssim_{C} (\thetastop + \beta)\ell(T_i). 
	\end{align}
	Since $z_i \in 5B_{T_i}$, if $\thetastop$ and $\beta$ are small enough (depending on $C$), this implies $z_{i+1} \in 6B_{T_i}$. Since $z_{i+1} \in 5B_{T_{i+1}}$, this finishes the proof. 
\end{proof}

\begin{lem}\label{l:union-of-the-Gs}
	Proposition \ref{p:local-to-global-prime} \eqref{prop:small-balls} holds with constant $C_6$ depending only on $C,C_5$ and $n$. 
\end{lem}

\begin{proof}
	Let $\{T_k\}_{k=0}^K$ be the sequence obtain by applying Lemma \ref{l:chain} and consider the set 
	\begin{align}
		G \coloneqq \bigcup_{k=0}^{K} \grid{T_k}. 
	\end{align}
	By \eqref{e:J-contains-ball} and Lemma \ref{l:face-glue}, we see that $\grid{T_{k+1}} \cap \grid{k_i} \neq \emptyset$ for each $k \in \{0,\dots,K-1\}$. In particular, $\grid{T_{k+1}} \cup \grid{T_k}$ is connected. Since this hold for arbitrary $k$, we have that $G$ is connected. Thus, by Lemma \ref{l:local-PI}, $G$ supports a $(1,1)$-Poincar\'e inequality with constant depending on $K = K(C,C_5)$ and $n$. Since $x,y \in G \subseteq Y$, this completes the proof. 
\end{proof}

\begin{constants}
	We will assume for the remainder of the paper that the constants $(\taustop,\thetastop,\Mclose,\Astop,\gridk)$ are chosen such that the Proposition \ref{l:sigma-good} holds. 
\end{constants}

\section{Approximations for layers of cubes}\label{s:skeleta-per-tree}

In this short section we define the sets $\surface{S,k}$ and $\surface{N,k}$. The first set will be a skeletal approximation for a full stopping-time region $S \in \calT^N(k)$ (as constructed in Section \ref{s:skeleta-per-cube}). The second set will simply be the union each of the skeletal approximation at layer $N$ i.e., 
\begin{align}
	\surface{N,k} \coloneqq \bigcup_{S \in \calT^N(k)} \surface{S,k}. 
\end{align}
These will be constructed from the collection $\{\grid{Q}\}_{Q \in \cubes}$. However, the definition will depend on whether $S$ lies $\mathscr{G}_1^N(k)$ or $\mathscr{B}_1^N(k)$ (recall \eqref{e:good-stop} and \eqref{e:bad-stop}). Indeed, for $S \in \mathscr{G}^N_1(k)$, let 
\begin{align}
	\surface{S,k} \coloneqq  \bigcup_{\substack{Q \in \Stop(N,k) \\ 3B_Q \cap 2\rho^{-1}B_{Q(S)} \neq \emptyset}} \grid{Q}. 
\end{align}
For $S \in \mathscr{B}^N_1(k)$, let 
\begin{align}
	\surface{S,k} \coloneqq  \grid{Q(S)} . 
\end{align}

Before proving more serious properties of these skeletal approximations in the next section, we end this section by showing two relatively easy properties. First, we show that each set $\surface{S,k}$ is upper regular and contained in a sufficiently large dilate of the ball $B_{Q(S)}$.

\begin{lem}\label{l:contained-in-ball}
	For each $S \in \calT^N(k)$, we have that $\surface{S,k}$ is upper Ahlfors $d$-regular and
	\begin{align}\label{e:contained-in-ball}
		\surface{S,k} \subseteq 2\sqrt{n}\lambdadyadic B_{Q(S)}.
	\end{align}
\end{lem}

\begin{proof}
	The upper regularity of $\surface{S,k}$ follows directly from Proposition \ref{l:sigma-good}. For \eqref{e:contained-in-ball}, if $S \in \mathscr{B}_1^N(k)$, this follows directly from \eqref{e:J-contains-ball}. If $S \in \mathscr{G}_1^N(k)$ and $Q \in \Stop(N,k)$ are such that $3B_Q \cap 2\rho^{-1}B_{Q(S)} \neq \emptyset$, it follows from Lemma \ref{l:Stop-k} that $\ell(Q) \leq \rho\ell(Q(S))$. Thus, for each $x \in \grid{Q} \subseteq 2\sqrt{n}\lambdadyadic B_Q$, we have 
	\begin{align}
		|x - x_{Q(S)}| &\leq (3+2\sqrt{n}\lambdadyadic)\ell(Q) + 2\rho^{-1}\ell(Q(S)) \\
		&\leq (5\sqrt{n}(1+4\rho^{-1}) \rho + 2\rho^{-1})\ell(Q(S)) \\
		&\leq (25\sqrt{n} + 2\rho^{-1})\ell(Q(S)) \leq 2\sqrt{n}\lambdadyadic \ell(Q(S)),
	\end{align}
	as required. 
\end{proof}

\begin{lem}\label{l:I_R-sigma}
	Let $S \in \calT^N(k)$, 
	\begin{align}
		Q \in 
		\begin{cases}
			\min(S(k)) & \mbox{ if } S \in \mathscr{G}_1^N(k); \\
			\min(S(k)) \cup \{Q(S)\} & \mbox{ if } S \in \mathscr{B}_1^N(k). 
		\end{cases}
	\end{align}
	and $R \in {\rm Child}(Q)$ such that $x_R = x_Q$. Then, supposing $\gridk = \gridk(\taustop)$ is sufficiently large, there is a $d$-dimensional $I_Q$ such that 
	\begin{align}
		I_Q \subseteq \surface{S,k} \cap (\innercube/8)B_R \quad \mbox{ and } \quad \ell(I_Q) \geq \innercube n^{-1/2} \ell(R)/10. 
	\end{align}
\end{lem}

\begin{proof}
	Let $k_3 \in \N$ such that 
	\begin{align}
	\rho \innercube n^{-1/2} \ell(Q)/20	< 2^{-k_3} \leq \rho \innercube n^{-1/2} \ell(Q)/10 .
	\end{align}
	Let $K$ be the dyadic cube such that $K \in \Delta_{k_3}$ and $x_Q \in K$, and let $I_Q$ be an arbitrary $d$-dimensional subface. By our choice of $k_3$, we have $\diam(K) \leq (\rho \innercube/10)\ell(Q)$. Since $x_Q \in K$ and $I_Q \subseteq K$, this implies 
	\begin{align}
		I_Q & \subseteq (\rho \innercube/8)B_Q = (\innercube/8)B_R; \\
		\ell(I_Q) &\geq \rho \innercube n^{-1/2} \ell(Q)/10 = \innercube n^{-1/2} \ell(R)/10.  \label{e:size-I_Q}
	\end{align}
	It only remains to show that 
	\begin{align}\label{e:I_Q-subset} 
		I_Q \subseteq \surface{S,k}. 
	\end{align}
	
	Suppose to begin with that $S \in \calG_1^N(k)$ and $Q \in \min(S(k))$. By Lemma \ref{l:subset-of-G_Q} and \eqref{e:size-I_Q} we have $I_Q \subseteq G_Q$ for $\gridk$ sufficiently large. Since, $G_Q \subseteq \surface{S,k}$ in this case, this implies \eqref{e:I_Q-subset}. Suppose now that $S \in \mathscr{B}_1^N(k)$. If $Q = Q(S)$ then the same argument as above give \eqref{e:I_Q-subset}, so let us suppose $Q \in \min(S(k))$. By Lemma \ref{l:stopped is pretty big} and \eqref{e:size-I_Q}, we have 
	\begin{align}
		\ell(I_Q) \gtrsim \ell(Q) \gtrsim \taustop \ell(Q(S)).
	\end{align}
	Thus, for $\gridk = \gridk(\taustop)$ sufficiently small, Lemma \ref{l:subset-of-G_Q} implies $I_Q \subseteq G_{Q(S)} = \surface{S,k},$ as required.  
\end{proof}

Next, we show every point in $Q \in S(k)$, for some $S \in \calT^N(k)$, lies close to $\surface{S,k}$. 

\begin{lem}\label{l:close-point}
	Let $S \in \calT^N(k)$, $Q \in S(k)$ and $x \in Q$. Then, 
	\begin{align}
		\dist(x,\surface{S,k}) \leq 2\ell(Q). 
	\end{align}
\end{lem}

\begin{proof}
	By definition, there exists $R \in \min(S(k))$ such that $R \subseteq Q$. By Lemma \ref{l:I_R-sigma}, there is a point $z \in \surface{S,k}  \cap  (\innercube/8) B_R \subseteq \surface{S,k} \cap B_Q$. Since $x \in Q \subseteq B_Q$, the lemma follows by the triangle inequality. 
\end{proof}

\section{Local Poincar\'e inequalities for a single layer}

In later sections we will want to glue the sets $\{\surface{N,k}\}_{N \in \N}$ that we defined in the previous section. We do this by first lifting these sets into $\R^{n+1}$, and then building bridges between consecutive layers. It will be important to know that any point in $\surface{S,k}$ (more accurately, a lifted version of this set) can be connected to one of these bridges within some larger subspace of $\surface{N,k}$ which supports a Poincar\'e inequality. We call this set $\Gamma_{S,k}$. There will be at least one neighbour of $S$ (which we call $U_S$) which is connected to a bridge, so we only need to know that $\Gamma_{S,k}$ connects $\surface{S,k}$ and $\surface{U_S,k}$ in a nice way. It will also be useful to know how these sets $\{\Gamma_{S,k}\}_{S \in \calT^N(k)}$ interact with each other. It turns our that we can build these sets in such a ways that any pair is either contained in an even larger space supporting a Poincar\'e inequality or is well-separated. Let us be more precise. 

\begin{notation}
	For $S \in \calT^N(k)$, let $U_S \in \calT^N(k)$ (possibly equal to $S$) be the stopping-time region such that $x_{Q(U_S)} = x_{Q}$, where $Q \in \Stop(N-1)$ is the parent of $Q(S)$. 
\end{notation}

\begin{constants}
	Let $\Asecond = \Afirst^5$. 
\end{constants}

Our main goal in this section is to prove the following. 

\begin{prop}\label{p:Gamma}
	Let $k \in \N$. There is a collection of sets $\{\Gamma_{S,k}\}_{S \in \calT^N(k)}$ in $\surface{N,k}$, each of which is Ahlfors $d$-regular and supports a $(1,1)$-Poincar\'e inequality with uniform constants such that 
	\begin{align}\label{e:Gamma1}
		\surface{S,k} \cup \surface{U_S,k} \subseteq \Gamma_{S,k}
	\end{align}
	for each $S \in \calT^N(k)$. If $S,T \in \calT^N(k)$ are such that $\Afirst B_{Q(S)} \cap \Afirst B_{Q(T)} \neq \emptyset$, either there exists a set 
	\begin{align}\label{e:Gamma-S-T} 
		\Gamma_{S,T,k} \subseteq \surface{N,k} \cap \Asecond B_{Q(S)} 
	\end{align}
	which contains $\Gamma_{S,k} \cup \Gamma_{T,k}$ and supports a $(1,1)$-Poincar\'e inequality, or 
	\begin{align}\label{e:Gamma2}
		\dist(\Gamma_{S,k} , \Gamma_{T,k}) \gtrsim \ell(Q(S)) + \ell(Q(T)). 
	\end{align}
	In the case that $2\rho^{-1} B_{Q(S)} \cap 2\rho^{-1} B_{Q(T)} \neq \emptyset$, the first alternative always holds. 
\end{prop}

Before defining the sets $\Gamma_{S,k}$ and $\Gamma_{S,T,k}$, we show that, for any given $S \in \mathscr{G}_2^N(k)$, there exists a family of sets satisfying something like the conclusion of Proposition \ref{p:Gamma} (and, indeed, $\Gamma_S$ will be chosen to be one of these sets for some particular choice of parameter).

\begin{lem}\label{l:good}
	For each $A \in [\Afirst,\Afirst^2]$ and $S \in \mathscr{G}^N_2(k)$ there is a set \begin{align}\label{e:Gamma-in-ball}
		\Gamma_{S,k}^A \subseteq \surface{N,k} \cap 10A^2 B_{Q(S)}
	\end{align}
	which is Ahlfors $d$-regular, supports a $(1,1)$-Poincar\'e inequality and satisfies 
	\begin{align}\label{e:sigma-T}
		\surface{T,k} \subseteq \Gamma_{S,k}^A
	\end{align}
	for all $T \in \calT^N(k)$ such that $\surface{T} \cap AB_{Q(S)} \neq \emptyset$. Moreover, if $A,A' \in [\Afirst,\Afirst^2]$ and $S,S' \in \mathscr{G}_2^N(k)$ are such that $AB_{Q(S)} \subseteq A'B_{Q(S')}$, then   \begin{align}\label{e:nested}
		\Gamma_{S,k}^{A} \subseteq \Gamma_{S',k}^{A'}.
	\end{align} 
\end{lem}

\begin{proof}
	Fix some $A \in [\Afirst,\Afirst^2]$ and $S \in \mathscr{G}^N_2(k)$. Let $A^* \coloneqq 4\sqrt{n}\lambdadyadic \Ccomp + A \in [\Afirst,2\Afirst^2]$ and set 
	\begin{align}
		\mathscr{S} &\coloneqq \{T \in \calT^N(k) \colon E \cap A^*B_{Q(S)} \cap Q(T) \neq \emptyset \}; \\
		\mathscr{T} &\coloneqq \{T \in \calT^N(k) \colon \surface{T,k} \cap AB_{Q(S)} \neq \emptyset \}. 
	\end{align}
	Recalling from Lemma \ref{l:contained-in-ball} that $\surface{T} \subseteq 2\sqrt{n}\lambdadyadic B_{Q(T)}$, by applying Lemma \ref{l:comp-size} and the triangle inequality, one has 
	\begin{align}\label{e:scrT}
		4\rho^{-1}B_{Q(T)} \subseteq 2\sqrt{n}\lambdadyadic B_{Q(T)} \subseteq A^*B_{Q(S)} \mbox{ for each } T \in \mathscr{T}. 
	\end{align}
	In particular,  
	\begin{align}
		\mathscr{T} \subseteq \mathscr{S}. 
	\end{align}
	The same argument which is used to establish \eqref{e:scrT} can also be used to show that  
	\begin{align}\label{e:scrS}
		2\sqrt{n}\lambdadyadic B_{Q(T)}\subseteq (4\sqrt{n}\lambdadyadic \Ccomp +A^*)B_{Q(S)} \subseteq 3A^2 B_{Q(S)} \mbox{ for each } T \in \mathscr{S}. 
	\end{align} 
	Let 
	\begin{align}
		\mathscr{R}_{T} &\coloneqq \{R \in \min(T(k)) \colon E \cap A^*B_{Q(S)} \cap R \neq \emptyset \} \mbox{ for } T \in \mathscr{G}^N_1(k), \\
		\mathscr{R}_{T} &\coloneqq \{Q(T)\} \mbox{ for } T \in \mathscr{B}^N_1(k),
	\end{align}
	and set 
	\begin{align}
		\{R_j\}_{j \in J} \coloneqq \bigcup_{T \in \mathscr{S}} \mathscr{R}_{T} \quad \mbox{ and } \quad \Gamma \coloneqq \Gamma_{S,k}^A \coloneqq \bigcup_{j \in J} \grid{R_j}. 
	\end{align}
	It follows directly from \eqref{e:scrS} that $\Gamma \subseteq  10A^2 B_{Q(S)}$. Property \eqref{e:nested} is also easy to see since the definitions of $\mathscr{S}$ and $\mathscr{R}_T$ are monotone with respect to set inclusions of balls.

	To prove that $\Gamma$ is Ahlfors regular and supports a Poincar\'e inequality we check that $\{R_j\}_{j \in J}$ satisfies the hypotheses of Lemma \ref{l:sigma-good} with constants $A^*$ and $C \leq \Cprop$. Lemma \ref{l:sigma-good} \eqref{p:R-covers-E} and \eqref{p:R-intersects-E} are immediate from the above the definitions. Lemma \ref{l:sigma-good} \eqref{p:R-comp} follows from Lemma \ref{l:comp-size-cutoff} whenever $R_i \in \mathscr{R}_{T}$ and $R_j \in \mathscr{R}_{T'}$ for some $T,T' \in \mathscr{G}_1^N(k)$ and from Lemma \ref{l:comp-size} (recall also Remark \ref{r:Top(N)}) whenever  or $T,T' \in \mathscr{B}^N_1(k)$. Thus, it only remains to check Lemma \ref{l:sigma-good} \eqref{p:R-comp} in the case that $R_i \in \min(T(k)) \subseteq \Stop(N,k)$ for some $T \in \mathscr{G}_1^N(k)$ and $R_j = \{Q(T')\}$ for some $T' \in \mathscr{B}_1^N(k)$. Since $A^*B_{Q(T)} \cap A^*B_{Q(T')} \neq \emptyset$, Lemma \ref{l:comp-size} implies
	\begin{align}\label{e:S'-T'}
		\ell(Q(T)) \sim \ell(Q(T')) \quad \mbox{ and } \quad B_{Q(T)} \subseteq (M/2)B_{Q(T')}
	\end{align}
	for $M$ sufficiently large. In particular, $R_i \subseteq (M/2)B_{Q(T')}$. Lemma \ref{l:stopped is pretty big} now implies that $\ell(R_i) \geq \cLB \taustop\ell(Q(T')) = \Cprop \ell(Q(T'))$, as required. \\ 
	
	We finish the proof of Lemma \ref{l:good} by checking \eqref{e:sigma-T}. Fix some $T \in \mathscr{T}$. If $T \in \mathscr{B}_1^N(k)$, then $\{Q(S)\} = \mathscr{R}_{T}  \subseteq \{R_j\}_{j \in J}$. Since $\surface{T,k} = \grid{Q(T)}$, the inclusion \eqref{e:sigma-T} now follows. Suppose instead that $T \in \mathscr{G}_1^N(k)$. Recalling the definition of $\surface{T,k}$ in this case, it suffices to show that  
	\begin{align}\label{e:R_S}
		\{Q \in \Stop(N,k) \colon 3B_Q \cap 2\rho^{-1}B_{Q(T)} \neq \emptyset \} \subseteq \bigcup_{T \in \mathscr{S}} \mathscr{R}_{T}. 
	\end{align}
	Fix some $Q \in \Stop(N,k)$ as above. It follows from Lemma \ref{l:Stop-k} that $Q \subseteq 4\rho^{-1}B_{Q(T)}$. Recalling \eqref{e:scrT}, this implies $Q \subseteq A^*B_{Q(S)}$. Thus, if $T' \in \calT^N(k)$ is such that $Q \in \min(T'(k))$, then $E \cap A^*B_{Q(S)} \cap Q(T') \neq \emptyset$ i.e., $T' \in \mathscr{S}$. The property $Q \subseteq A^*B_{Q(S)}$ also implies that $Q \in \mathscr{R}_{T'}$, which finishes the proof of \eqref{e:R_S}.
\end{proof}

We now define
\begin{align}
	\Gamma_{S,k} \coloneqq
	\begin{cases}
		\Gamma_{S,k}^{\Afirst} & \mbox{ if } S \in \mathscr{G}_2^N(k); \\
		\surface{S,k} \cup \surface{U_S,k} & \mbox{ if } S \in \mathscr{B}_2^N(k). 
	\end{cases}
\end{align}

\begin{proof}[Proof of Proposition \ref{p:Gamma} for $S \in \mathscr{G}_2^N(k)$]
	The inclusion \eqref{e:Gamma1} and the Poincar\'e inequality follow directly from Lemma \ref{l:good}. Fix now some some $T \in \calT^N(k)$ such that 
	\begin{align}
		\Afirst B_{Q(S)} \cap \Afirst B_{Q(T)} \neq \emptyset. 
	\end{align}
	It follows Lemma \ref{l:comp-size-cutoff} and the triangle inequality that 
	\begin{align}\label{e:inclusion-good}
		2\rho^{-1}B_{Q(T)} \subseteq \Afirst B_{Q(T)} \subseteq (1+2\Ccomp)\Afirst B_{Q(S)} \subseteq A_1^2 B_{Q(S)}. 
	\end{align}
	We set $\Gamma_{S,T,k} \coloneqq \Gamma_{S,k}^{A_1^2}$ (as in Lemma \ref{l:good}). It is clear that $\Gamma_{S,T,k}$ is Ahlfors $d$-regular, supports a $(1,1)$-Poincar\'e inequality, contains $\Gamma_S$ and that 
	\begin{align}\label{e:better-estimate}
		\Gamma_{S,T,k} \subseteq 10A_1^4 B_{Q(S)} \subseteq \Asecond B_{Q(S)}
	\end{align} 
	(we note the slightly better first inclusion for later on). If $T \in \mathscr{G}_2^N(k)$, then $\Gamma_{T,k} \subseteq \Gamma_{S,T,k}$ by \eqref{e:nested} and \eqref{e:inclusion-good}. Since $\surface{T,k}$ and $\surface{U_T,k}$ have non-empty intersection with $2\rho^{-1}B_{Q(T)}$, if $T \in \mathscr{B}_2^N(k)$, then the inclusion $\Gamma_{T,k} \subseteq \Gamma_{S,T,k}$ follows by \eqref{e:sigma-T} and \eqref{e:inclusion-good}. 
\end{proof}

\begin{proof}[Proof of Proposition \ref{p:Gamma} for $S \in \mathscr{G}_2^N(k)$]
	The inclusion \eqref{e:Gamma1} follows from the definition. Since $2\sqrt{n}\lambdadyadic \dyadic{Q(S)} \cap 2\sqrt{n}\lambdadyadic \dyadic{Q(U_S)} \neq \emptyset$, the Poincar\'e inequality and Ahlfors regularity for $\Gamma_S$ follows from Lemma \ref{l:face-glue}, Lemma \ref{l:local-PI} and Lemma \ref{l:ADR-gluing} once we show that 
	\begin{align}\label{e:U_S}
		U_S \in \mathscr{B}_1^N(k).
	\end{align} 
	To see \eqref{e:U_S}, note that $\Mclose B_{Q(U_S)} \subseteq 2\Mclose B_{Q(S)}$ for $\Mclose = \Mclose(\Ccomp)$ sufficiently large, by Lemma \ref{l:comp-size}. The inclusion \eqref{e:U_S} now follows since $S \in \mathscr{B}_2^N(k)$. 
	
	Now, fix some $T \in \calT^N(k)$ such that 
	\begin{align}
		\Afirst B_{Q(S)} \cap \Afirst B_{Q(T)} \neq \emptyset. 
	\end{align}
	If $T \in \mathscr{G}_2^N(k)$ then we set $\Gamma_{S,T,k} \coloneqq \Gamma_{T,S,k}$, where $\Gamma_{T,S,k}$ is the set obtained by applying Proposition \ref{p:Gamma} for $T$. It follows from Lemma \ref{l:comp-size} and the stricter inclusion in \eqref{e:better-estimate} that 
	\begin{align}
		\Gamma_{S,T,k} \subseteq 10\Afirst^4 B_{Q(T)} \subseteq \Afirst^5 B_{Q(S)}. 
	\end{align} Suppose instead that $T \in \mathscr{B}_2^N(k)$. In this case, we have 
	\begin{align}
		\Gamma_{S,k} \cup \Gamma_{T,k} &= \surface{S,k} \cup \surface{U_S,k} \cup \surface{T,k} \cup \surface{U_T,k} \\
		&= \grid{S} \cup \grid{U_S} \cup \grid{T} \cup \grid{U_T}
	\end{align}
	(note that \eqref{e:U_S} also holds for $T$). Thus, if $\Gamma_{S,k} \cup \Gamma_{T,k}$ is connected, we have by Lemma \ref{l:local-PI} that it is Ahlfors regular and supports a Poincar\'e inequality. Note that $\Gamma_{S,k} \cup \Gamma_{T,k}$ is always connected in the case that $2\rho^{-1}B_{Q(S)} \cap 2\rho^{-1}B_{Q(T)} \neq \emptyset$. If it is not connected, then Lemma \ref{l:face-glue} implies \eqref{e:Gamma2}.  
\end{proof}

\section{Bridges between consecutive layers}\label{s:bridges}

For $N \in \N$ and $S \in \calT^N(\stopk)$, consider the sets 
\begin{align}\label{e:surface-lift}
	\surfacelift{S,\stopk} = \text{Graph}_{\surface{S,\stopk}}(\taustop d_{N,\stopk}) \mbox{ and } \surfacelift{N,\stopk} = \text{Graph}_{\surface{N,\stopk}}(\taustop d_{N,\stopk})
\end{align}
(recall the notation $\text{Graph}_A(f)$ from \eqref{e:graph-notation}). In this section we describe a collection of bridges $\{b_Q\}_{Q \in \text{Stop}(N-1,\stopk)}$ between consecutive layer $\surfacelift{N-1,\stopk}$ and $\surfacelift{N,\stopk}$. In particular, we will prove the following.

\begin{prop}\label{p:b_Q}
	For each $k,N \in \N$ there is a collection of set $\{b_{Q,k}\}_{k \in \Stopprime(N-1,k)}$ in $\R^{n+1}$ satisfying the following conditions. 
	\begin{enumerate}
		\item \label{prop:b_Q-ADR-PI} Each $b_{Q,k}$ is Ahlfors $d$-regular with \begin{align}
			\diam(b_{Q,k}) \sim \ell(Q),
		\end{align}
		 and supports a $(1,1)$-Poincar\'e inequality with uniform constants, and satisfies 
		\item \label{prop:b_Q-subset} We have
		\begin{align}
			\pi_n(b_{Q,k}) \subseteq (\innercube/2)B_Q
		\end{align}
		and 
		\begin{align}
			b_{Q,k} \subseteq  \{z \in \R^{n+1} \colon \pi_1(z) \in [\taustop d_{N,k}(\pi_n(z)) ,\taustop d_{N-1}(\pi_n(z))]\}.
		\end{align}
		\item \label{prop:b_Q-upper} The set $\surfacelift{N-1,k} \cap b_{Q,k}$ is lower Ahlfors $d$-regular with diameter comparable to $\ell(Q)$. Moreover, if $S \in \calT^{N-1}(k)$ is such that $Q \in S$, then 
		\begin{align}
			\surfacelift{N-1,k} \cap b_{Q,k} = \surfacelift{S,k} \cap b_{Q,k}.
		\end{align}
		\item \label{prop:b_Q-lower} The set $\surfacelift{N,k} \cap b_{Q,k}$ is lower Ahlfors $d$-regular with diameter comparable to $\ell(Q)$. Moreover, if $T \in \calT^N(k)$ is such that $x_{Q(T)} = x_Q$, then 
		\begin{align}
			\surfacelift{N,k} \cap b_{Q,k} = \surfacelift{T,k} \cap b_{Q,k}.
		\end{align}
		\item \label{prop:b_Q-gluing} Denoting
		\begin{align}
			A_1 \coloneqq b_{Q,k} \cap \surfacelift{N-1,k}\quad \mbox{ and } \quad A_2 \coloneqq b_{Q,k} \cap \surfacelift{N,k},
		\end{align}
		for each $x \in b_{Q,k}$, 
		\begin{align}
			y_1 &\in \{z \in \R^{n+1}\colon \pi_1(z) \geq  \taustop d_{N-1,k}(\pi_n(z)) \} \supseteq \surfacelift{N-1,k}, \\
			y_2 &\in \{z \in \R^{n+1}\colon \pi_1(z) \leq  \taustop d_{N,k}(\pi_n(z)) \} \supseteq \surfacelift{N,k},
		\end{align}
		we have 
		\begin{align}\label{e:comp-metric} 
			\inf_{a \in A_i} | x - a | + | a - y_i | \lesssim | x-y_i|.
		\end{align}
	\end{enumerate}
\end{prop}

\subsection{Preliminary results} Before detailing on the construction of the $b_Q$, we state a useful technical preliminaries. Lemma \ref{l:useful} below give us a way to construct a bi-Lipschitz image of $d$-dimensional cube which connects the graph (over $\R^n$) of a $(d-1)$-segment in $\R^n$ to a different graph (over $\R^n$) of a different $(d-1)$-segment in $\R^n$. Assuming the graphs are sufficiently well-separated with small enough Lipschitz constant, we are able to do this in such a way that these bi-Lipschitz images interact well (keeping in mind Lemma \ref{l:gluing}) with each of the graphs. This will be important for proving a Poincar\'e inequality later. The proofs are technical, but rather easy. As such, we have decided to include them in Appendix \ref{a:useful}.

\begin{lem}\label{l:useful}
	Let $ 0 \leq \eta,\theta , \lambda <1,$ $0 < \gamma, \delta <1$ and $C_1 >1$. Let $Q,R \subseteq \R^n$ be $d$-dimensional cubes with $\ell(Q) \geq \ell(R)$, $I \subseteq Q$ and $J \subseteq R$ be $(d-1)$-dimensional faces, $P \subseteq \R^n$ be the $d$-dimensional plane containing $R$, $u \colon P \to P^\perp$ be $1$-Lipschitz and $K$ be the graph of $u$ over $R$. Let $g_1 \colon \R^n \to [0,\infty)$ be $\eta$-Lipschitz and $g_2 \colon \R^n \to [0,\infty)$ be $\theta$-Lipschitz such that
	\begin{align}\label{e:useful-1}
		\gamma \ell(I) \leq \inf_{x \in I} g_1(x) - \sup_{y \in u(J)} g_2(y) \quad \mbox{and} \quad  \sup_{x \in I} g_1(x) \leq C_1\ell(I).
	\end{align} 
	Suppose there exists a map $h \colon I \times [0,\ell(I)] \to \R^n$ which takes $I$ to $J$ such that 
	\begin{align}\label{e:useful-2}
		\delta | x-y| \leq | h(x,t) - h(y,t) | \leq |x-y|
	\end{align}
	for all $x,y \in I$ and $t \in [0,\ell(I)]$, 
	\begin{align}\label{e:useful-2.5}
		| h(x,\ell(I)) - h(y,\ell(I)) | \leq C_1  \delta |x-y|
	\end{align}
	and 
	\begin{align}\label{e:useful-3}
		|h(x,t) - h(x,s)| \leq C_1|t-s|
	\end{align}
	for all $x \in I$ and $t,s \in [0,\ell(I)]$. Finally, suppose that is a ball $B \subseteq \R^n$ such that 
	\begin{align}\label{e:all-in-ball}
		Q , u(R) , h(I \times [0,\ell(I)]) \subseteq B. 
	\end{align}
	Then, for $\eta = \eta(\gamma,\delta,C_1)$ and $\theta = \theta(\gamma,C_1)$ small enough, there exists a bi-Lipschitz map $H \colon I \times [-\ell(I),2\ell(I)] \to B \times \R$ such that  
	\begin{align}
		H(I,[-\ell(I),0]) &= \emph{Graph}_Q(g_1)   \subseteq \R^{n+1}, \label{e:useful-4} 
	\end{align} 
	and \begin{align}
		H(I,[\ell(I),2\ell(I)) &= \emph{Graph}_{K}(g_2)  \subseteq \R^{n+1}. \label{e:useful-4.5}
	\end{align}
	Furthermore, for each $x \in I$ and $t \in [0,\ell(I)]$,  
	\begin{align}\label{e:useful-5.6}
		g_2(\pi_{\R^n}(H(x,t))) + \tfrac{\gamma \cdot (\ell(I) - t)}{2}  \leq \pi_{\R}(H(x,t)) \leq g_1(\pi_{\R^n}(H(x,t))) - \tfrac{\gamma t}{2} .
	\end{align}
	Finally, denoting
	\begin{align}
		X \coloneqq  H(I,[-\ell(I),2\ell(I)]), \ X_{\rm up} \coloneqq H(I,[-\ell(I),0]), \  X_{\rm down} \coloneqq H(I,[\ell(I),2\ell(I)]), 
	\end{align}
	if $\eta = \eta(\gamma)$ and $\theta = \theta(\gamma)$ are sufficiently small, we have 
	\begin{align}\label{e:banana-4}
		\inf_{a \in X_{\rm up}} | x - a| + |a - y| \lesssim |x-y|
	\end{align}
	for all $x \in X$ and $y \in \R^{n+1}$ such that $\pi_1(y) \geq g_1(\pi_n(y))$, and 
	\begin{align}\label{e:banana-5}
		\inf_{a \in X_{\rm down}}| x - a| + |a - y| \lesssim |x-y|
	\end{align}
	for all $x \in X$ and $y \in \R^{n+1}$ such that $\pi_1(y) \leq \lambda g_2(\pi_n(y))$.
\end{lem}

\subsection{The part of the bridge in \texorpdfstring{$\surfacelift{N-1,\stopk}$}{TQ-star}} Fix some $N \in \N$ and $Q \in \Stopprime(N-1,\stopk)$ for the remainder of the section. Let 
\begin{align}
	S &\in \calT^{N-1}(\stopk) \mbox{ such that } Q \in S; \\
	T &\in \calT^N(\stopk) \mbox{ such that } x_{Q(T)} = x_Q. 
\end{align}
The construction of the bridge $b_Q$ will depend on whether $T \in \mathscr{G}_1^N(\stopk)$ or $T \in \mathscr{B}_1^N(\stopk)$. However, in either case, the intersection of $b_Q$ with $\surfacelift{N-1,\stopk}$ will be the same. In this subsection, we identify this set. \\

Throughout this section, $\upperface_R$ will denote the denote the $d$-dimensional cube from Lemma \ref{l:I_R-sigma}, and $\upperboundary_R$ will denote an arbitrary $(d-1)$-dimensional sub-face, which we fix. In this way, we have 
\begin{align}\label{e:I_R-F_R}
	\upperboundary_R \subseteq \upperface_R \subseteq (\rho\innercube/8)B_R \quad \mbox{ and } \quad  \ell(\upperboundary_R) = \ell(\upperface_R) \sim \ell(R). 
\end{align}
By scaling the metric, we will assume without loss of generality that 
\begin{align}
	\ell(\upperface_Q) = \ell(\upperboundary_Q) = 1. 
\end{align}
For the cube $Q$ (and only $Q$), we set 
\begin{align}\label{e:T_Q^*-prime}
	\upperface_Q^* \coloneqq \text{Graph}_{\upperface_Q}(\taustop d_{N-1,\stopk}) \quad \mbox{ and }  \quad
	\upperboundary_Q^* \coloneqq \text{Graph}_{\upperboundary_Q}(\taustop d_{N-1,\stopk}) ,
\end{align}
which are both subsets of $\surfacelift{S,\stopk}$. As we will see later, $\upperface_Q^*$ will be the intersection of $b_Q$ with $\surfacelift{N-1,\stopk}$. For now, let us collect some basic facts about these sets

\begin{lem}
	The set $\upperface_Q^*$ is Ahlfors $d$-regular with diameter comparable to $\ell(Q)$. For each $x,y \in (\innercube/2)B_Q$, we have 
	\begin{align}\label{e:d_N-T_Q-prime}
		\taustop d_{N-1,\stopk}(x) \lesssim \ell(\upperboundary_Q)
	\end{align}
	and 
	\begin{align}\label{e:well-sep-bad}
		\taustop d_{N-1,\stopk}(x) - 2\taustop d_{N,\stopk}(y) \gtrsim \ell(\upperboundary_Q)
	\end{align}
\end{lem}

\begin{proof}
	The Ahlfors regularity of $\upperface_Q^*$ follows since it is bi-Lipschitz equivalent to a $d$-dimensional cube with side-length comparable to $\ell(Q)$. Inequality \eqref{e:d_N-T_Q-prime} follows from Lemma \ref{l:stop-size-cutoff} and \eqref{e:I_R-F_R}. Inequality \eqref{e:well-sep-bad} follows from Lemma \ref{l:well-separated} and \eqref{e:I_R-F_R}.  
\end{proof}

We now split to construction into our two cases.

\subsection{Construction of \texorpdfstring{$b_Q$}{bQ} when \texorpdfstring{$T \in \mathscr{B}_1^N(\stopk)$}{Q(T) in B}.} For brevity, let us denote $\upperface_T \coloneqq I_{Q(T)}$ and $\upperboundary_T \coloneqq \upperface_{Q(T)}$. Then, set 
\begin{align}\label{e:F-lift-prime}
	\lowerfacebad_{T}^* = \text{Graph}_{\lowerfacebad_{T}}(\taustop d_{N,\stopk}) \mbox{ and } \lowerboundarybad_{T}^* = \text{Graph}_{\lowerboundarybad_{T}}(\taustop d_{N,\stopk}),
\end{align}
which are both subsets of $\surfacelift{T,\stopk} \subseteq \surfacelift{N,\stopk}$ (recall Lemma \ref{l:I_R-sigma}).

\begin{rem}
	We note that $\lowerfacebad_{T}^*$ is the graph of the function $\taustop d_{N,\stopk}$ as opposed to $\upperface_Q^*$ which is the graph of the function $\taustop d_{N-1,\stopk}$ (similar for $\lowerboundarybad_T$ versus $\upperboundary_Q$). We hope that this slight abuse of notation does not cause significant misunderstanding. 
\end{rem} 

Since $Q(T)$ is a child of $Q$, we have by \eqref{e:I_R-F_R} that 
\begin{align}\label{e:FsimT-prime}
	\ell(\lowerboundarybad_T) \sim \ell(\upperboundary_Q) =1 . 
\end{align}
Since $I_R \subseteq (\innercube/8)B_R$ for each $R \in \cubes$, and $(\innercube/8)B_{Q(T)} \subseteq (\innercube/8) B_Q$, there is a rigid motion $h \colon \upperboundary_Q \times [0,1] \to (\innercube/8)B_Q$ taking $\upperboundary_Q$ to $\lowerboundarybad_T$ such that 
\begin{align}\label{e:h_1-bi-lip}
	\frac{\ell(\lowerboundarybad_T)}{\ell(\upperboundary_Q)}|x-y| \leq |h(x,t) - h(y,t) | \leq |x-y|
\end{align}
for all $x,y \in \upperboundary_Q$ and $t \in [0,1]$, 
\begin{align}\label{e:isom-1}
	| h(x,1) - h(y,1) | = \frac{\ell(\lowerboundarybad_T)}{\ell(\upperboundary_Q)}|x-y|
\end{align}
for all $x,y \in \upperboundary_Q$, and
\begin{align}\label{e:h-in-time-2}
	|h(x,t) - h(x,s) | \lesssim |t-s|
\end{align}
for all $x,y \in \upperboundary_Q$ and $t,s \in [0,1].$ 

\begin{lem}\label{l:hypotheses}
	The hypotheses of Lemma \ref{l:useful} are satisfied with the data 
	\begin{align}
		(Q,R) &= (\upperface_Q,\lowerfacebad_T), \ (I,J) = (\upperboundary_Q,\lowerboundarybad_T), \ u \equiv 0,  \\
		(g_1, g_2) &= (\taustop d_{N-1,\stopk}, \taustop d_{N,k} ), \ h =h  \ \mbox{and} \ B = (\innercube/8)B_Q, 
	\end{align}
	and parameters 
	\begin{align}
		\eta = \theta = \taustop, \ \gamma \sim \delta \sim C_1 \sim 1. 
	\end{align}
\end{lem}

\begin{proof}
	The properties in the preamble of Lemma \ref{l:useful} are immediate from the definitions with the stated parameters. The inequalities in \eqref{e:useful-1} with parameters $\gamma \sim C_1 \sim 1$ follow from \eqref{e:d_N-T_Q-prime} and \eqref{e:well-sep-bad}  since $\upperboundary_Q,\lowerboundarybad_T \subseteq (\innercube/8)B_Q$. That $h$ satisfies \eqref{e:useful-2}--\eqref{e:useful-3} with parameters $\delta \sim C_1 \sim 1$, follows from the definition above and \eqref{e:FsimT-prime}. The inclusions in \eqref{e:all-in-ball} also follows directly from the construction since $(\innercube/8)B_{Q(T)} \subseteq (\innercube/8)B_Q$. 
\end{proof}

\begin{proof}[Proof of Proposition \ref{p:b_Q} when $T \in \mathscr{B}_1^N(\stopk)$]
	By Lemma \ref{l:hypotheses}, if we choose $\taustop = \taustop(1)$ small enough there is a bi-Lipschitz map $H \colon \upperboundary_Q \times [-1,2] \to (\innercube/8)$ satisfying the conclusion of Lemma \ref{l:useful} with that data. We set 
	\begin{align}
		b_Q \coloneqq H(\upperboundary_Q \times [-1,2]). 
	\end{align}
	Proposition \ref{p:b_Q} \eqref{prop:b_Q-ADR-PI} holds since $b_Q$ is the bi-Lipschitz image of a $d$-dimensional rectangle with each side-lengths comparable to $1 \sim \ell(Q)$. Proposition \ref{p:b_Q} \eqref{prop:b_Q-subset} holds by \eqref{e:useful-5.6} and the fact that $H$ maps into $(\innercube/8)B_Q$. Proposition \ref{p:b_Q} \eqref{prop:b_Q-upper} follows from \eqref{e:useful-4} and \eqref{e:useful-5.6} since these imply
	\begin{align}
		\surfacelift{N-1,\stopk} \cap b_Q = H(\upperboundary_Q \times [-1,0]) = \text{Graph}_{\upperface_Q}(\taustop d_{N-1,\stopk}) = \upperface_Q^* \subseteq \surfacelift{S,\stopk}.
	\end{align}
	Proposition \ref{p:b_Q} \eqref{prop:b_Q-lower}  follows from \eqref{e:useful-4.5} and \eqref{e:useful-5.6} since these imply 
	\begin{align}
		\surfacelift{N,\stopk} \cap b_Q = H(\upperboundary_Q \times [1,2]) = \text{Graph}_{\lowerfacebad_T}(\taustop d_{N,\stopk}) = \lowerfacebad_T^* \subseteq \surfacelift{T,\stopk}.
	\end{align}
	Proposition \ref{p:b_Q} \eqref{prop:b_Q-gluing} follows directly from \eqref{e:banana-4} and \eqref{e:banana-5}. 
\end{proof}

\subsection{Construction of \texorpdfstring{$b_Q$}{bQ} when \texorpdfstring{$T \in \mathscr{G}_1^N(\stopk)$.}{Q(T) in G}} We start by setting some notation. 

\begin{notation}
	Let $\Gamma$ denote the graph constructed from Lemma \ref{l:lip-graphs} with $M = 10$ and $\{Q_k\} = \{Q(T)\}$. Let $P \subseteq \R^n$ and $g \colon P \to P^\perp$ be the corresponding $d$-plane and map, respectively, such that $\Gamma$ is the graph of $g$ over $P$. Through this section $\pi$ will denote the orthogonal projection onto $P$. Since $E \cap \innercube B_{Q(T)} \subseteq Q(T)$ note, by Lemma \ref{l:lip-graphs-subset}, that 
	\begin{align}\label{e:gamma-is-covered}
		\Gamma \cap (\innercube/2)B_{Q(T)} \subseteq \bigcup_{\substack{R \in \min(T(\stopk))  \\ \Gamma \cap (\innercube/2)B_{Q(T)} \cap 2B_R \neq \emptyset  }} 2B_R.
	\end{align}
\end{notation}

\begin{notation}
Let $I_0 \subseteq P$ be the $d$-dimensional cube in $P$ such that 
\begin{align}\label{e:J'-prime}
\ell(I_0) = \ell(I_Q) \mbox{ and } \pi(x_Q)  = x_{I_0},
\end{align}   
and let $F_0$ be an arbitrary $(d-1)$-dimensional face.
\end{notation}

\begin{lem}\label{l:F_0}
	We have 
	\begin{align}
		F_0 \subseteq I_0 \subseteq (\innercube/4)B_{Q(T)}.
	\end{align}
	In particular, there exists a rigid motion $h \colon \upperboundary_Q \times [0,1] \to (\innercube/4)B_{Q(T)}$ which takes $\upperboundary_Q$ to $F_0$ and satisfies \eqref{e:useful-2}--\eqref{e:useful-3} with parameters $\delta \sim C_1 \sim 1$. 
\end{lem}

\begin{proof}
	Recall from Lemma \ref{l:I_R-sigma} that $\upperface_Q \subseteq (\innercube/8)B_{Q(T)}$. Since $I_0$ is isometric to $\upperface_Q$, this implies $I_0 \subseteq B(x_{I_0},\innercube \ell(Q(T))/8)$. Since $\|g\|_\infty \leq \thetastop \ell(Q(T))$ and $\dist(x_Q,\Gamma) = \dist(x_{Q(T)},\Gamma) \leq \thetastop \ell(Q(T))$ (see Lemma \ref{l:lip-graphs}), we have $|x_Q - x_{I_0}| \leq 2\thetastop \ell(Q(T))$. The triangle inequality now implies the desired inclusion, so long as $\thetastop$ is sufficiently small. 
\end{proof}

\begin{notation}
	Let  
	\begin{align}
		\calQ \coloneqq \left\{R \in \min(T(\stopk)) \colon \Gamma \cap \pi^{-1}(I_0) \cap (\innercube/8) B_R \neq\emptyset  \right\}.
	\end{align}
\end{notation}

\begin{lem}\label{l:Gamma-is-covered}
	We have 
	\begin{align}
		\Gamma \cap \pi^{-1}(I_0) \subseteq (\innercube/4)B_{Q(T)} \cap  \bigcup_{R \in \calQ} \Ccover B_R, 
	\end{align}
	where $\Ccover \coloneqq 17\sqrt{d}\Ccomp^2 \rho^{-1}$. 
\end{lem}

\begin{proof}
	Let $x \in \Gamma \cap \pi^{-1}(I_0)$. Since $\|g\|_\infty \leq \thetastop \ell(Q(T))$ and $I_0 \subseteq (\innercube/4)B_{Q(T)}$, we already have 
	\begin{align}
		x  \in (\innercube/4)B_{Q(T)}
	\end{align}
	for $\thetastop$ sufficiently small. Thus, by \eqref{e:gamma-is-covered}, there is $R \in \min(T(\stopk))$ such that $x \in 2B_R$. Let $m \in \N$ such that $R \in \cubes_m$ and let $l \in \N$ be the large integer less than $m$ such that $\rho^l \geq 2\Ccomp \rho^m$, with $\Ccomp$ as in Lemma \ref{l:comp-size-cutoff}. By maximality,
	\begin{align}
		2\Ccomp \rho^m \leq \rho^l < 2\Ccomp \rho^{m-1}.
	\end{align}
	
	Consider a point $x' \in \Gamma \cap \pi^{-1}(I_0)$ such that 
	\begin{align}
		|x-x'| \leq 10\sqrt{d}\Ccomp\ell(R) \quad \mbox{ and } \quad \dist(\pi(x'),\partial I_0) \geq 4\Ccomp \ell(R),
	\end{align}
	(imagine sliding $\pi(x)$ towards the centre of $I_0$ and letting $x'$ be the associated graph point). By Lemma \ref{l:lip-graphs}, there exists a point $y \in E$ such that $|x'-y| \leq \thetastop \ell(R)$. Let $R' \in \cubes_l$ such that $y \in R'$ and note that 
	\begin{align}\label{e:I_0-inclusion}
		\pi(\innercube B_{R'}) \subseteq I_0 \quad \mbox{ and } \quad B_{R'} \subseteq  15\sqrt{d}\Ccomp \rho^{-1} B_R. 
	\end{align}  
	Indeed, for each $z \in B_{R'}$ we have (by our choice of $l$) that   
	\begin{align}
		|z - x'| &\leq |z-x'| + |x'-y| \leq (2+\thetastop)\ell(R') \leq 3\Ccomp \rho^{-1} \ell(R);\\
		|z - x_R| &\leq |z - x'| + |x'-x| + |x-x_R| \leq (3\Ccomp \rho^{-1} + 10\sqrt{d}\Ccomp + 2)\ell(R).
	\end{align}
	The second inclusion in \eqref{e:I_0-inclusion} is immediate from the second estimate above, the first inclusion follows since $\dist(\pi(x'),\partial I_0) \geq 4\Ccomp \ell(R)$. We note also that there exists $R'' \in \min(T(\stopk))$ such that 
	\begin{align}\label{e:R''}
		\mbox{$\innercube B_{R''} \subseteq \innercube B_{R'}$ and $\Ccomp\ell(R'') \geq \ell(R)$.}
	\end{align}
	Indeed, by the first inclusion in \eqref{e:I_0-inclusion} and Lemma \ref{l:F_0}, we must have $R' \subseteq Q(T)$. However, $R'$ cannot be contained in any minimal cube from $\min(T(\stopk))$ since $\ell(R') \geq 2\Ccomp \ell(R)$ (by our choice of $l$) and $B_{R'} \subseteq  15\sqrt{d}\Ccomp \rho^{-1} B_R$ (otherwise we contradict Lemma \ref{l:comp-size-cutoff}). The cube $R''$ from \eqref{e:R''} is then just the cube in $\min(T(\stopk))$ with the same centre as $R'$. Finally, we note that
	\begin{align}
		R'' \in \calQ. 
	\end{align} 
	Indeed, Lemma \ref{l:lip-graphs} implies $\Gamma \cap (\innercube/8)B_{R''} \neq \emptyset$ which, when combined with \eqref{e:I_0-inclusion} and \eqref{e:R''}, give $\Gamma \cap \pi^{-1}(I_0) \cap (\innercube/8)B_{R''} \neq \emptyset$. 
	
	The proof of Lemma \ref{l:Gamma-is-covered} will be complete once we show 
	\begin{align}
		x \in 17\sqrt{d}\Ccomp \rho^{-1} C B_{R''}. 
	\end{align}
	This follows from the second inclusion in \eqref{e:I_0-inclusion} and \eqref{e:R''} since 
	\begin{align}
		|x - x_{R''}| \leq |x- x_{R}| + |x_R-x_{R''}| &\leq (2+15\sqrt{d}\Ccomp \rho^{-1})\ell(R) \\
		&\leq 17\sqrt{d}\Ccomp^2 \rho^{-1}\ell(R''),
	\end{align}
	as required. 
\end{proof}

\begin{lem}\label{l:calQ-subset}
	For each $R \in \calQ$, we have 
	\begin{align}
		(\innercube/4)B_R \subseteq (\innercube/2)B_Q.
	\end{align}
\end{lem}

\begin{proof}
	Let $x \in (\innercube/4)B_R$ and $y \in \Gamma \cap \pi^{-1}(I_0) \cap (\innercube/8)B_R$ (which exists by definition). Then, by applying Lemma \ref{l:Gamma-is-covered} and the fact that $\ell(R) \leq \ell(Q(T)) = \rho \ell(Q)$, we have 
	\begin{align}
		|x-x_Q| \leq |x-y| + |y-x_Q| \leq \innercube \ell(R)/2 + \innercube\ell(Q(T))/4 \leq \ell(Q)/2,
	\end{align} 
	which implies the state inclusion. 
\end{proof}

\begin{notation}
	For each $R \in \calQ$, let  
	\begin{align}
		z_R \in \Gamma \cap \pi^{-1}(I_0) \cap  (\innercube/8)B_R.
	\end{align}
	Let $\planecube_R$ be a $d$-dimensional cube in $P$ such that 
	\begin{align}\label{e:planecube}
		y_R \coloneqq \pi(z_R) \in \planecube_R \subseteq I_0 \quad \mbox{ and } \quad \diam(\planecube_R) = (\innercube/8)\ell(R).
	\end{align}
	Let 
	\begin{align}\label{e:sigma-I-prime}
		\graphcube_R = \{ x + z_R - y_R : x \in D_R\} 
	\end{align} 
	be the translate of $\planecube_R$ through $z_R$, and let $\graphboundary_R$ be an arbitrary $(d-1)$-dimensional sub-face. Finally, let 
	\begin{align}
		\planeset = \bigcup_{R \in \calQ} \planecube_R \subseteq I_0 .
	\end{align}
	\end{notation}

	\begin{lem}\label{l:graphcube}
		For each $R \in \calQ$, we have 
		\begin{align}\label{e:props-sigma-R-prime}
				\graphboundary_R \subseteq \sigma_R \subseteq (\innercube/4)B_R  \quad \mbox{ and } \quad  \ell_R \coloneqq \ell(\graphboundary_R) \sim \ell(R). 
			\end{align}
		In particular, there exists a rigid motion $h_R \colon \graphboundary_R \times [0,\ell_R] \to (\innercube/4)B_R$ which takes $\graphboundary_R$ to $\upperboundary_R$ and satisfies \eqref{e:useful-2}--\eqref{e:useful-3} with parameters $\delta \sim C_1 \sim 1$.
	\end{lem}
	\begin{proof}
		Since $\graphcube_R$ is isometric to $\planecube_R$, it follows that $\diam(\graphcube_R) = (\innercube/8)\ell(R)$ (in particular, this implies $\ell(\graphcube_R) \sim \ell(R)$). Since $z_R \in \graphcube_R$ and $z_R \in (\innercube/8)B_R$, we deduce the left-most inclusion in \eqref{e:props-sigma-R-prime} after applying the triangle inequality. 
	\end{proof}

\begin{lem}\label{l:tilde-u-prime}
The set $\sigma \coloneqq \bigcup_{R \in \calQ} \graphcube_R$ is the graph of a $C\alphastop$-Lipschitz map $\tilde{u} \colon \planeset \to P^\perp$ which satisfies
\begin{align}\label{e:tilde-u-size}
	|\tilde{u}(x) | \leq \thetastop \ell(Q)
\end{align}
for all $x \in \planeset$. 
\end{lem}

\begin{proof}
We define $\tilde{u} \colon \planeset \to P^\perp$ by the formula $\tilde{u}(x) = z_R - y_R$ whenever $x \in D_R$ (i.e., constant on each $\planecube_R$). In this way, we have 
\begin{align}\label{e:g-y_I}
	\tilde{u}(y_R) = g(y_R) \quad \mbox{for each } R \in \calQ. 
\end{align}
By construction, $\graphset$ is the graph of $\tilde{u}$ over $D$. Recall from Lemma \ref{l:lip-graphs} that $\|g\|_\infty \leq\thetastop \ell(Q(T)) \leq \thetastop\ell(Q)$. Since $z_R \in \Gamma$, it follows that 
\begin{align}
	|\tilde{u}(x)| = |\tilde{u}(y_R)| = |g(y_R)| \leq \thetastop \ell(Q)
\end{align}
for each $x \in \planecube_R$. Thus, it only remains to check that it is $C\alphastop$-Lipschitz. Fix $w,x\in \planeset$ and let $R , R' \in \calQ$ such that $w \in \planecube_R$ and $x \in \planecube_{R'}$. If $R = R'$ then $\tilde{u}(w) = \tilde{u}(x)$ and there is nothing to show. Suppose instead that $R \neq R'.$ First, by \eqref{e:g-y_I} and that fact that $g(y_R),g(y_{R'}) \in \Gamma$ (which is $C\alphastop$-Lipschitz), we have 
\begin{align}\label{e:u-1}
	|\tilde{u}(w) - \tilde{u}(x)| = |\tilde{u}(y_R) - \tilde{u}(y_{R'})| = |g(y_R) - g(y_{R'})| \leq C\alphastop |y_R - y_{R'}|. 
\end{align}
Thus, it suffices to show that 
\begin{align}
	|y_R - y_{R'}| \lesssim |w-x|. 
\end{align}
For this, note to begin with that, since $z_{R} \in \graphcube_R \subseteq (\innercube/4)B_R$, $z_{R'} \in \graphcube_{R'} \subseteq (\innercube/4)B_{R'}$ and $\innercube B_R \cap \innercube B_{R'} = \emptyset$, we have
\begin{align}
	|z_R - z_{R'}| \geq (3\innercube/4)\left( \ell(R) + \ell(R') \right). 
\end{align}
Since $y_R, y_{R'}$ are the points in $\Gamma$ above $z_R,z_{R'}$, respectively, and $\Gamma$ is $C\alphastop$-Lipschitz, this implies that 
\begin{align}
	|y_{R} - y_{R'}| \geq (\innercube/2)\left( \ell(R) + \ell(R')\right). 
\end{align}
Thus, after recalling the choice of diameter for $\planecube_R$ and $\planecube_{R'}$ from \eqref{e:planecube}, we have by the triangle inequality that 
\begin{align}\label{e:u-2}
	|w - x| &\geq |y_R - y_{R'}| - |w-y_R| - |x - y_{R'} | \\
	&\geq |y_{R} - y_{R'} | - (\innercube/8)\ell(R) - (\innercube/8)\ell(R') \\
	&\geq |y_{R} - y_{R'}| - \tfrac{1}{2}|y_R - y_{R'}| \geq \tfrac{1}{2} |y_R - y_{R'}|,
\end{align}
as required.
\end{proof}

\begin{notation}\label{n:G_Q}
	Let $u \colon P \to P^\perp$ be a $C\alphastop$-Lipschitz extension of $\tilde{u}$ to the whole of $P$. Let 
	\begin{align}\label{e:G_Q-prime}
		G_Q = \text{Graph}_{I_0}(u)\quad  \mbox{ and } \quad G^*_Q = \text{Graph}_{G_Q}(2\taustop d_{N,\stopk}).
	\end{align}
\end{notation}

We emphasise here the factor of 2 appearing on the right-hand side. This will be important for us as it gives a little wiggle room between the sets $G_Q^*$ and $\surfacelift{N,\stopk}$.

\begin{lem}
	We have 
	\begin{align}\label{e:u-size}
		\| u \|_{L^\infty(I_0)} \lesssim (\alphastop + \thetastop) \ell(Q). 
	\end{align}
	Furthermore, for each $x \in G_Q^*$ there exists $R \in \calQ$ such that 
	\begin{align}\label{e:G-close}
		\dist(x,\graphcube_R^*) \lesssim \ell(\graphboundary_R). 
	\end{align}
\end{lem} 

\begin{proof}
	The estimate \eqref{e:u-size} follows from \eqref{e:tilde-u-size} and the fact that $\diam(I_0) \leq \ell(Q)$ and $u$ is $C\alphastop$-Lipschitz. To see \eqref{e:G-close}, fix some $x \in \graphcube$, let $y = \pi_n(x) \in G_Q$ and $z = \pi(y) \in I_0$. By Lemma \ref{l:Gamma-is-covered}, there exists $R \in \calQ$ such that $z \in \pi(\Ccover B_R)$. Since $\graphcube_R \subseteq (\innercube/4)B_R$, this implies $\dist(z,\pi(\graphcube_R)) \lesssim \ell(R)$. Since $\{y\} \cup \graphcube_R \subseteq G_Q$ and $G_Q$ is a Lipschitz graph with small constant, this implies $\dist(y,\graphcube_R) \lesssim \ell(R)$. For similar reasons, we also have \eqref{e:G-close}. 
\end{proof}

\begin{lem}\label{l:C_Q-bilip}
The hypotheses of Lemma \ref{l:useful} are satisfied with the data 
\begin{align}
	(Q,R) &= (\upperface_Q,I_0), \ (I,J) = (\upperboundary_Q,F_0), \ u = u,  \\
	(g_1, g_2) &= (\taustop d_{N-1,\stopk}, 2\taustop d_{N,\stopk} ), \ h =h  \ \mbox{and} \ B = (\innercube/2)B_Q, 
\end{align}
and parameters 
\begin{align}
	\eta = \taustop, \ \theta = 2\taustop, \ \gamma \sim \delta \sim C_1 \sim 1. 
\end{align}
In particular, for $\taustop$ small enough, there exists a bi-Lipschitz map $H \colon \upperboundary_Q \times [-1,2] \to (\innercube/2)B_Q \times \R$ such that, if we denote 
\begin{align}
	C_Q \coloneqq H(\upperboundary_Q \times [-1,2]),
\end{align}
then the following conditions hold. First, 
\begin{align}
	\upperface_Q^* &= H(\upperboundary_Q \times [-1,0]) = C_Q \cap \surfacelift{N-1,\stopk}; \label{e:H-connect-2-1} \\
	G_Q^* &= H(\upperboundary_Q \times [1,2]) = C_Q \cap \emph{Graph}_{\R^n}(2\taustop d_{N,\stopk}). \label{e:H-connect-2-2}
\end{align}
Second, for each $x \in \upperboundary_Q$ and $t \in [0,1],$ 
\begin{align}\label{e:H-between-2*}
	2\taustop d_{N,\stopk}(\pi_{n}(H(x,t))) + c(1-t) &\leq \pi_1(H(x,t)) \\
	&\leq \taustop d_{N-1,\stopk}(\pi_{n}(H(x,t))) -  ct.
\end{align}
Finally, we have  
	\begin{align}\label{e:pears-1}
		\inf_{a \in \upperface_Q^*} | x- a | + | a -y | \lesssim |x-y|
	\end{align}
	for all $x \in C_Q$ and $y \in \R^{n+1}$ such that $\pi_1(y) \geq \taustop d_{N-1,\stopk}(\pi_n(y))$. 
\end{lem}

\begin{proof}
	Once we verify the hypotheses of Lemma \ref{l:useful} with the stated data and parameters, properties \eqref{e:H-connect-2-1}--\eqref{e:pears-1} follow directly from Lemma \ref{l:useful} after some definition chasing. The conditions in the preamble of Lemma \ref{l:useful} with the stated parameters are immediate from the definitions above. By \eqref{e:I_R-F_R}, Lemma \ref{l:F_0} and \eqref{e:u-size}, we have 
	\begin{align}\label{e:C_Q-inclusion}
		\upperboundary_Q,u(F_0),h(\upperboundary_Q\times [0,1]) \subseteq (\innercube/2)B_Q.
	\end{align}
	In particular, \eqref{e:all-in-ball} holds with $B = (\innercube/2)B_Q$. The inequalities in \eqref{e:useful-1} with parameters $\gamma \sim C_1 \sim 1$ follow from \eqref{e:d_N-T_Q-prime}, \eqref{e:well-sep-bad} and \eqref{e:C_Q-inclusion}. Finally, properties \eqref{e:useful-2}--\eqref{e:useful-3} follow from Lemma \ref{l:F_0} with parameters $\delta \sim C_1 \sim 1$. 
\end{proof}

\begin{notation}
For each $R \in \calQ$, let 
\begin{align}
\graphcube_R^* &\coloneqq \text{Graph}_{\graphcube_R}(2\taustop d_{N,\stopk}) \quad \mbox{ and } \quad \graphboundary_R^* \coloneqq \text{Graph}_{\graphboundary_R}(2\taustop d_{N,\stopk}),
\end{align}
and
\begin{align}
\llowerface_R^* &\coloneqq  \text{Graph}_{\llowerface_R}(\taustop d_{N,\stopk}) \quad \mbox{ and } \quad \llowerboundary_R^*\coloneqq \text{Graph}_{\llowerboundary_R}(\taustop d_{N,\stopk}). 
\end{align}
\end{notation}

\begin{lem}\label{l:L_I-bilip}
For each $R \in \calQ$, the hypotheses of Lemma \ref{l:useful} are satisfied with the data 
\begin{align}
	(Q,R) &= (\graphcube_R,\upperface_R), \ (I,J) = (\graphboundary_R,\upperboundary_R), \ u = 0,  \\
	(g_1, g_2) &= (2\taustop d_{N,\stopk}, \taustop d_{N,\stopk} ), \ h =h_R  \ \mbox{and} \ B = (\innercube/4)B_R, 
\end{align}
and parameters 
\begin{align}
	\eta = 2\taustop, \ \theta = \taustop, \ \gamma \sim \delta \sim C_1 \sim 1. 
\end{align}
In particular, for $\taustop$ small enough, there exists a bi-Lipschitz map $H_R \colon \graphboundary_R \times [-\ell_R,2\ell_R] \to (\innercube/4)B_R \times \R$ such that, if we denote 
\begin{align}
	L_R \coloneqq H_R(\graphboundary_R \times [0,\ell_R]),
\end{align} 
then the following conditions holds. First, 
\begin{align}
	\graphcube_R^* &= H_R(\graphboundary_R\times [-\ell_R ,0]) = L_R \cap \emph{Graph}_{\R^n}(2\taustop d_{N,\stopk}) \label{e:lemons-1}; \\ 
	\llowerface_R^* &= H_R(\graphboundary_R \times [\ell_R,2\ell_R]) = L_R \cap \surfacelift{N,\stopk}\label{e:lemons-2}.
\end{align}
Second, for each $x \in \graphboundary_R$ and $t \in [0,\ell_R]$, 
\begin{align}\label{e:H_I-between}
	\taustop d_{N,\stopk}(\pi_{n}(H_R(x,t))) + c(\ell_R-t) &\leq \pi_{1}(H_R(x,t))  \\
	&\leq 2\taustop d_{N,\stopk}(\pi_{n}(H_R(x,t))) -  c t .
\end{align}
Finally, we have
\begin{align}\label{e:lemons-3}
	\inf_{a \in \llowerface_R^*} | x- a | + | a -y | \lesssim |x-y|.
\end{align}
$x \in L_R$ and $y \in \R^{n+1}$ with $\pi_1(y) \leq \taustop d_{N,\stopk}(\pi_n(y))$.
\end{lem}

\begin{proof}
	Fix $R \in \calQ$. Once we verify the hypotheses of Lemma \ref{l:useful} with the stated data and parameters, properties \eqref{e:lemons-1}--\eqref{e:lemons-3} follow directly from Lemma \ref{l:useful} after some definition chasing. The conditions in the preamble of Lemma \ref{l:useful} with the stated parameters are immediate from the definitions above. To check \eqref{e:useful-1}, let us fix $x \in \graphboundary_R$ and $y \in \upperboundary_R$. By \eqref{e:I_R-F_R} and \eqref{e:props-sigma-R-prime}, we know
	\begin{align}
	x,y \in (\innercube/4)B_R \implies	|x-y| \lesssim \ell(R) \sim \ell(\graphboundary_R). 
	\end{align}
	Thus, by applying Lemma \ref{l:stop-size-cutoff}, we have 
	\begin{align}
		2\taustop d_{N,\stopk}(x) - \taustop d_{N,\stopk}(y) &= \taustop d_{N,\stopk}(x) - \taustop (d_{N,\stopk}(y) - d_{N,\stopk}(x) ) \\
		&\geq c\ell(\graphboundary_R) - C\taustop\ell(\graphboundary_R) \gtrsim \ell(\graphboundary_R)
	\end{align}
	and
		\begin{align}
			2\taustop d_{N,\stopk}(x) \lesssim  \ell(\graphboundary_R). 
		\end{align}
	The above two inequalities imply \eqref{e:useful-1} with parameters $\gamma \sim C_1 \sim 1$. Properties \eqref{e:useful-2}--\eqref{e:useful-3} follows from Lemma \ref{l:graphcube} with parameters $\delta \sim C_1 \sim 1$. The inclusions in \eqref{e:useful-4} with $B= (\innercube/4)B_R$ follow from \eqref{e:I_R-F_R} and \eqref{e:props-sigma-R-prime}.
\end{proof}

We define
\begin{align}
	b_Q \coloneqq C_Q \cup \bigcup_{R \in \calQ} L_R. 
\end{align}

\subsection{Proof of Proposition \ref{p:b_Q} when \texorpdfstring{$T \in \mathscr{G}_1^N(\stopk)$.}{Q(S) in G}} Proposition \ref{p:b_Q} \eqref{prop:b_Q-ADR-PI} is a consequence of Lemma \ref{l:b_Q-ADR} and Lemma \ref{l:b_Q-Poincare} below. 

\begin{lem}\label{l:b_Q-ADR}
	The bridge $b_Q$ is Ahlfors $d$-regular with $\diam(b_Q) \sim \ell(Q)$. 
\end{lem}

\begin{proof}
	By Lemma \ref{l:ADR-gluing}, it suffices to check that $\bigcup_{R \in \calQ} L_R$ is Ahlfors regular, which we will check this Lemma \ref{l:gluing-ADR}. For each $R \in \calQ$, let $w_R \in L_R$ and $r_R \coloneqq \ell_R$. By Lemma \ref{l:L_I-bilip}, $\graphboundary_R \subseteq L_R$, which implies $r_R \leq \diam(L_R) \lesssim r_R$. In particular, $\{L_R\}$ and $\{B(w_R,r_R)\}$ satisfy Lemma \ref{l:gluing-ADR} \eqref{prop:cover} with $C_1' \sim 1$. By Lemma \ref{l:Gamma-is-covered}, $\Gamma \cap \pi^{-1}(I_0)$ and $\{B(w_R,r_R)\}$ satisfy Lemma \ref{l:gluing-ADR} \eqref{prop:close} and \eqref{prop:cover-X} with $C_1'' \sim 1$. Then, for $C_1 \coloneqq \max\{C_1',C_1''\}$, it follows from Lemma \ref{l:cubes-bounded overlap} that Lemma \ref{l:gluing-ADR} \eqref{prop:bounded-overlap} is satisfied with constant $C_1$ and $C_2$ depending on $d$ and $C_1$. This finishes the proof.
\end{proof}

Our next goal is the check that $b_Q$ supports a Poincar\'e inequality. The following lemma is key. 

\begin{lem}\label{l:special-Poincare}
	For each distinct pair $R,R' \in \calQ$, each of the sets $C_Q, \ C_Q \cup L_R$ and $C_Q \cup L_R \cup L_{R'}$ support a $(1,1)$-Poincar\'e inequality with uniform constants.
\end{lem}

\begin{proof}
	By Lemma \ref{l:C_Q-bilip} and Lemma \ref{l:L_I-bilip}, each of $C_Q, L_R$ and $L_{R'}$ support a $(1,1)$-Poincar\'e inequality since they are bi-Lipschitz to $d$-dimensional cubes. To verify the other cases we  appeal to Lemma \ref{l:gluing}. \\
	
	We start with $C_Q \cup L_R.$ Observe to begin with that 
	\begin{align}\label{e:int}
		C_Q \cap L_R = \graphcube_R^*. 
	\end{align} 
	Indeed, by \eqref{e:H-between-2*} and \eqref{e:H_I-between} we have $C_Q \cap L_R \subseteq \text{Graph}_{\R^n}(2\taustop d_{N,\stopk})$. Recalling from \eqref{e:lemons-1} that $\graphcube_R^* = L_R \cap \text{Graph}_{\R^n}(2\taustop d_{N,\stopk})$, the forward inclusion now follows. For the reverse inclusion, note first that $\graphcube_R^* \subseteq G_Q^*$. The reverse inclusion then follows by \eqref{e:H-connect-2-2} and \eqref{e:lemons-2}. The identity \eqref{e:int} implies  
	\begin{align}\label{e:int-LR}
		C_Q \cap L_R \mbox{ is Ahlfors $d$-regular with } \diam(C_Q \cap L_R) \sim \diam(L_R). 
	\end{align}
	
	We now show that  
	\begin{align}\label{e:apple-2}
		\inf_{w \in C_Q \cap L_R} |y-w| + |w-z| \lesssim |y-z|
	\end{align}
	for points $y \in C_Q$ and $z \in L_R$. Let $y' \in \upperboundary_Q$ and $t \in [0,\ell]$ such that $y = H(y,t)$, and let $z' \in \sigma_R$ and  $s \in [0,\ell_R]$ such that $z = H(z',s).$ Furthermore, let $w \coloneqq H_R(z',0)$. Since $H_R$ is bi-Lipschitz, we have 
	\begin{align}\label{e:w-z}
		| w -z | \leq Cs. 
	\end{align}
	On the other hand, by \eqref{e:H_I-between}, 
	\begin{align}
		|\pi_1(w-z)| &
		\geq 2\taustop d_{N,\stopk}(\pi_n(w)) - (2\taustop d_{N,\stopk}(\pi_n(x)) - c s) \\
		&\geq 2\taustop d_{N,\stopk}(\pi_n(w)) - 2\taustop d_{N,\stopk}(\pi_n(x)) + c |w-z|. 
	\end{align}
	The estimate \eqref{e:apple-2} now follows directly from Lemma \ref{l:checking-the-gluing-condition} so long as $\taustop$ is sufficiently small. Properties \eqref{e:int-LR}, \eqref{e:apple-2} and Lemma \ref{l:gluing} imply that $C_Q \cup L_{R}$ supports a $(1,1)$-Poincar\'e inequality. \\

	We finish by checking that $C_Q \cup L_R \cup L_{R'}$ supports a $(1,1)$-Poincar\'e inequality. Let $A \coloneqq (C_Q \cup L_R)\cap L_{R'}$. By Lemma \ref{l:L_I-bilip}, 
	\begin{align}\label{e:legs-far}
		\pi_{\R^n}(L_R) \subseteq (\innercube/4)B_R \quad \mbox{ and } \quad \pi_{\R^n}(L_{R'}) \subseteq (\innercube/4)B_{R'}.
	\end{align}
	In particular, $L_R$ and $L_{R'}$ are disjoint, which implies $A = C_Q \cap L_{R'}$. Thus, by the same argument used to establish \eqref{e:int-LR}, we get 
	\begin{align}
		A \mbox{ is Ahlfors $d$-regular with } \diam(A) \sim \diam(L_{R'}).
	\end{align}
	We now fix $y \in C_Q \cup L_R$, $z \in L_{R'}$ and check that 
	\begin{align}\label{e:legs-metric}
		\inf_{w \in A} |y-w| + |w - z| \lesssim |y-z|. 
	\end{align}
	If $y \in C_Q$, this follows by the same argument used to establish \eqref{e:apple-2}. Suppose instead that $y \in L_{R'}$. In this case we have $|x-y| \gtrsim \ell(R') \sim \diam(L_{R'})$ by \eqref{e:legs-far}. Thus, Lemma \ref{l:diameter-PI} implies the estimate \eqref{e:legs-metric}. 
\end{proof}

\begin{lem}\label{l:b_Q-Poincare}
	The bridge $b_Q$ supports a $(1,1)$-Poincar\'e inequality. 
\end{lem}

\begin{proof}
	By Lemma \ref{l:special-Poincare}, any pair of points in $b_Q$ are contained in some Ahlfors regular subset of $b_Q$ which supports a $(1,1)$-Poincar\'e inequality with uniform constants. This property implies that $b_Q$ supports a $(1,1)$-Poincar\'e inequality. 
\end{proof}

With Proposition \ref{p:b_Q} \eqref{prop:b_Q-ADR-PI} out of the way, we continue with the remaining points.

\begin{lem}
Proposition \ref{p:b_Q} \eqref{prop:b_Q-subset} holds.
\end{lem}

\begin{proof}
	The bounds on $\pi_1(x)$ for each $x \in b_Q$ follow from \eqref{e:H-between-2*} and \eqref{e:H_I-between}. Since $H$ takes values in $(\innercube/2)B_Q \times \R$, we know that $\pi_n(C_Q) \subseteq (\innercube/2)B_Q$. Since $H_R$ takes values in $(\innercube/4)B_R \times \R$, we know that $\pi_n(L_R) \subseteq (\innercube/4)B_R \subseteq (\innercube/2)B_Q$ for each $R \in \calQ$ (recall Lemma \ref{l:calQ-subset} for this final inclusion).
\end{proof}

\begin{lem}
	Proposition \ref{p:b_Q} \eqref{prop:b_Q-upper} holds.
\end{lem}

\begin{proof}
	By Lemma \ref{l:well-separated} and \eqref{e:H_I-between}, we have $L_R \cap \surfacelift{N-1,\stopk} = \emptyset$ for each $R \in \calQ$. Thus, by \eqref{e:H-connect-2-1}, 
	\begin{align}
		\surfacelift{N-1,\stopk} \cap b_Q = \upperface_Q^* \subseteq \surfacelift{S,\stopk}. 
	\end{align} 
	Since $\upperface_Q$ is bi-Lipschitz to a $d$-dimensional cube with side legnth comparable to $\ell(Q)$, we conclude the lower regularity and diameter properties of $\surfacelift{N-1,\stopk} \cap b_Q$. 
\end{proof}

\begin{lem}
	Proposition \ref{p:b_Q} \eqref{prop:b_Q-lower} holds.
\end{lem}

\begin{proof}
	By Lemma \ref{l:well-separated} and \eqref{e:H-between-2*}, we have $C_Q \cap \surfacelift{N,\stopk} = \emptyset$. If $R \in \calQ$, it then follows from \eqref{e:lemons-2} and \eqref{e:H_I-between} that $L_R \cap \surfacelift{N,\stopk} = \upperface_R \subseteq \surfacelift{T,\stopk}$. In particular, 
	\begin{align}
		A^* \coloneqq \surfacelift{N,\stopk} \cap b_Q = \bigcup_{R \in \calQ} \upperface_R^* \subseteq \surfacelift{T,\stopk}.
	\end{align}
	It now only remains to check the lower regularity and diameter estimates for $A^*$. Since they are bi-Lipschitz equivalent, it suffices to check the same properties for the set 
	\begin{align}
		A \coloneqq \bigcup_{R \in \calQ} \upperface_R.
	\end{align} 
	We will check the hypotheses of Lemma \ref{l:gluing-ADR} with \begin{align}
		X = \Gamma \cap \pi^{-1}(I_0), \quad \{Y_i\} = \{\upperface_R\} \quad \mbox{ and } \quad \{B_i\} = \{B_R\}. 
	\end{align}
	First, the upper regularity of $X$ is immediately since $\Gamma$ is a graph and the upper regularity of each $Y_i$ follows from since they are $d$-dimensional cubes. Lemma \ref{l:gluing-ADR} \eqref{prop:cover} and \eqref{prop:close} follow directly from Lemma \ref{l:I_R-sigma} and the definition of $\mathcal{Q}$ with constant $C_1 = 1$. The bounded overlap of $\{5B_R\}$ (Lemma \ref{l:gluing-ADR} \eqref{prop:bounded-overlap}) follows directly from Lemma \ref{l:cubes-bounded overlap}. 
\end{proof}

Proposition \ref{p:b_Q} \eqref{prop:b_Q-gluing} is a consequence of Lemma \ref{l:go-up} and Lemma \ref{l:go-down} below. 

\begin{lem}\label{l:go-up}
	For each $x \in b_Q$ and $y \in \R^{n+1}$ such that $\pi_1(y) \geq \taustop d_{N-1,\stopk}(\pi_n(y))$ we have 
	\begin{align}\label{e:glue-up}
			\inf_{a \in \upperface_Q^*} |x-a| + |a-y| \lesssim |x-y|. 
	\end{align}
\end{lem}

\begin{proof}
	If $x \in C_Q$, then this already follows from \eqref{e:pears-1}, so let us suppose that $x \in L_R$ for some $R \in \calQ$. Since $\pi_n(x) \in (\innercube/2)B_Q$, if $\pi_n(y) \notin \innercube B_Q$ then $|x-y| \gtrsim \ell(Q) \sim \diam(b_Q)$ and the estimate follows from Lemma \ref{l:diameter-PI}. Suppose now that $x,y \in \innercube B_Q$. In this case, \eqref{e:H_I-between} and Lemma \ref{l:well-separated} imply 
	\begin{align}
		\pi_1(y) - \pi_1(x) \geq \taustop d_{N-1,\stopk}(\pi_n(y)) - 2\taustop d_{N,\stopk}(\pi_n(x)) \gtrsim \ell(Q) \sim \diam(b_Q).
	\end{align}
	The estimate \eqref{e:glue-up} now follows by Lemma \ref{l:diameter-PI}.
\end{proof}

\begin{lem}\label{l:go-down}
	For each $x \in b_Q$ and $y \in \R^{n+1}$ such that $\pi_1(y) \leq \taustop d_{N,\stopk}(\pi_n(y))$ there exists $R \in \calQ$ such that  
	\begin{align}\label{e:glue-down}
		\inf_{a \in \upperface_R^*} |x-a| + |a-y| \lesssim |x-y|. 
	\end{align}
\end{lem}

\begin{proof}
	If $x \in L_R$ for some $R \in \calQ$ then this already follows from \eqref{e:lemons-3}. Suppose that $x \in C_Q$. If $x \in \surfacelift{N-1,\stopk}$, then $\dist(x,\surfacelift{N,k}) \gtrsim \ell(Q) \sim \diam(b_Q)$ by Lemma \ref{l:well-separated}, and the estimate follows from Lemma \ref{l:diameter-PI}. Thus, we will assume also that $x \notin \surfacelift{N-1,\stopk}$. Define a point
	\begin{align}
		w \in C_Q \cap \text{Graph}_{\R^n}(2\taustop d_{N,\stopk})
	\end{align}
	as follows. If $x \in \text{Graph}_{\R^n}(2\taustop d_{N,\stopk})$, let $w = x$. Otherwise, let $x' \in \upperboundary_Q$ and $t \in [0,1]$ such that $x = H(x',t)$, and set $w = H(x',1) \in G_Q^*$. In either case, using the fact that $H$ is bi-Lipschitz, we have
	\begin{align}
		|w -x| \lesssim 1 - t. 
	\end{align}
	Let $R \in \calQ$ such that $\graphcube_R^*$ is closest to $w$, and let $y \in \graphboundary_R^*$. In this way, it follows from \eqref{e:G-close} that 
	\begin{align}
		|w - y | \lesssim \ell_R. 
	\end{align}
	Let $y' \in \graphboundary_R$ such that $y = H_R(y',0)$ and let $z = H_R(y',\ell_R)$. Since $H_R$ is bi-Lipschitz, we have 
	\begin{align}
		|y-z| \lesssim \ell_R. 
	\end{align}
	Combing the previous estimates, we conclude that 
	\begin{align}\label{e:x-z}
		|x-z| \lesssim 1-t + \ell_R. 
	\end{align}
	On the other hand, by \eqref{e:H-between-2*}, 
	\begin{align}\label{e:pi_1}
		\hspace{2em} |\pi_1(x-z)| &\geq 2\taustop d_{N,\stopk}(\pi_n(x)) + c(1-t) - \taustop d_{N,\stopk}(\pi_n(z)) \\
		&\geq \taustop d_{N,\stopk}(\pi_n(x)) - \taustop d_{N,\stopk}(\pi_n(z)) + c(1 - t ) + d_{N,\stopk}(\pi_n(x)). 
	\end{align}
	Lemma \ref{l:L_I-bilip} implies that $\pi_n(z) \in (\innercube/4)B_R$. Thus, if $1 - t \leq \ell_R$, then it follows from \eqref{e:x-z} that $\pi_n(x) \in CB_R$ for some $C > 1$, and from Lemma \ref{l:stop-size-cutoff} that $d_{N,\stopk}(\pi_n(x)) \sim \ell(R) \sim \ell_R$. If $1 - t > \ell_R$, then $1 -t \gtrsim 1 -t + \ell_R$. In either case, we have 
	\begin{align}
		1-t + d_{N,\stopk}(\pi_n(x)) \gtrsim 1-t+\ell_R .
	\end{align}
	Thus, carrying on from \eqref{e:pi_1} and applying \eqref{e:x-z}, we get 
	\begin{align}
		|\pi_1(x-z)| \geq\taustop d_{N,\stopk}(\pi_n(x)) - \taustop d_{N,\stopk}(\pi_n(z)) + c|x-z|.
	\end{align}
	The estimate \eqref{e:glue-down} now follows directly from Lemma \ref{l:checking-the-gluing-condition}.
\end{proof}

 \section{Local Poincar\'e inequalities for two consecutive layers}\label{s:conse-layers}

As mentioned at the beginning of Section \ref{s:skeleta-per-cube}, the set $\envelope$ from Theorem \ref{t:main} will be obtained as a limit of (a subsequence of) a sequence of sets $\{Z(k)\}_{k \in \N}$. We now have all the ingredients need to define this sequence. Indeed, for $k \in \N$, the set $Z(k)$ is defined as 
\begin{align}
	\envelope(k) \coloneqq \bigcup_{N \in \N}\ \left( \bigcup_{S \in \calT^N(k)}\surfacelift{S,k} \cup \bigcup_{Q \in \Stopprime(N,k)} b_{Q,k}\right). 
\end{align}

We will show in the next section that each $\envelope(k)$ is Ahlfors regular and supports a Poincar\'e inequality. In this section, we prove two useful preliminary results in this direction. First, we prove upper regularity of the union of a $\surfacelift{S,k}$ together with all its bridges (Proposition \ref{p:UR}). Second, we show that the union of the sets and bridge corresponding to nearby stopping-time regions in consecutive layers are contained in a subset of $\envelope(k)$ which supports a Poincar\'e inequality (Proposition \ref{p:PI-in-layers}). 

\begin{prop}\label{p:UR}
	For each $S \in \calT^N(k)$, the set 
	\begin{align}
		F_S \coloneqq \surfacelift{S,k} \cup \bigcup_{\substack{Q \in \Stopprime(N,k) \\ Q \subseteq Q(S)}} b_{Q,k}
	\end{align}
	is upper Ahlfors $d$-regular with diameter at most a constant multiple of $\ell(Q(S))$ and satisfies $\pi_n(F_S) \subseteq 2\sqrt{n}\lambdadyadic B_{Q(S)}$. 
\end{prop}

\begin{proof}
	The diameter estimate and the inclusion $\pi_n(F_S) \subseteq 2\sqrt{n}\lambdadyadic B_{Q(S)}$ follow directly from Lemma \ref{l:contained-in-ball} and Proposition \ref{p:b_Q} \eqref{prop:b_Q-subset}. By Lemma \ref{l:contained-in-ball} we need only check need only check up regularity of 
	\begin{align}
		\bigcup_{\substack{Q \in \Stopprime(N,k) \\ Q \subseteq Q(S)}} b_{Q,k}
	\end{align}
	to conclude upper regularity of $F_S$. For this, we we will check the hypotheses of Lemma \ref{l:gluing-ADR} with data 
	\begin{align}
		X = \surfacelift{N,k}, \quad \{Y_i\} = \{b_Q\}_{Q \in \Stopprime(N,k)} \quad \mbox{ and } \quad \{B_i\} = \{\mathfrak{B}_Q\}_{Q \in \Stopprime(N,k)}
	\end{align}
	for some collection $\{\mathfrak{B}_Q\}_{Q \in \Stopprime(N,k)}$ to be define momentarily. To begin with, the regularity properties with some uniform constant $C_0$ follow immediately from Lemma \ref{l:contained-in-ball} and Proposition \ref{p:b_Q} \eqref{prop:b_Q-ADR-PI}. To define the collection of balls, for each $Q \in \Stopprime(N,k)$ choose a point $\zeta_Q \in b_Q \cap \surfacelift{N,k}$ and consider
	\begin{align}
		\mathfrak{B}_Q \coloneqq B(\zeta_Q,(\innercube/2)\ell(Q)). 
	\end{align}
	It follows from Proposition \ref{p:b_Q} \eqref{prop:b_Q-ADR-PI} and \eqref{prop:b_Q-subset} that Lemma \ref{l:gluing-ADR} \eqref{prop:cover} holds with absolute constant $C_1 \sim 1$. It is easy to check that $\pi_n(5C_1\mathfrak{B}_Q) \subseteq 6C_1B_Q$. Thus, by a volume argument in $E$ we conclude that $\{5C_1\mathfrak{B}_Q\}_{Q \in \Stopprime(N,k)}$ has bounded overlap with some constant $C_2 \sim 1$. Finally, Lemma \ref{l:gluing-ADR} is trivial since we have chosen the balls $\mathfrak{B}_Q$ to be centred on $\surfacelift{N,k}$. 
\end{proof}

\begin{prop}\label{p:PI-in-layers}
	Let $k \in \N$, $L \geq 1$ and $S_L,T_L \in \calT^L(k)$ such that $\Afirst B_{Q(S_L)} \cap \Afirst B_{Q(T_L)} \neq \emptyset$. For any $Q_L \in \min(S_L) \cap \Stopprime(L,k)$ and $R_L \in \min(T_L) \cap \Stopprime(L,k)$ there is a set $X \subseteq \envelope(k)$ which is Ahlfors $d$-regular, supports a $(1,1)$-Poincar\'e inequality and contains
	\begin{align}\label{e:PI-in-layers}
		\bigcup_{i=L-1}^L \left( \surfacelift{S_i,k} \cup \surfacelift{T_i,k} \cup   b_{Q_i,k} \cup b_{R_{i},k}  \right),
	\end{align} 
	where $Q_{L-1},R_{L-1} \in \Stop(L-1)$ and $S_{L-1} , T_{L-1} \in \calT^{L-1}(k)$ are such that 
	\begin{align}
		Q(S_L) \subseteq Q_{L-1} \subseteq Q(S_{L-1}) \quad \mbox{ and } \quad Q(T_L) \subseteq R_{L-1} \subseteq Q(T_{L-1}).
	\end{align}
\end{prop}

We will often apply Proposition \ref{p:PI-in-layers} via the following easy corollary. 

\begin{cor}\label{l:bridge-up-and-down}
	Let $k \in \N$, $L \geq 1$, $Q \in \Stop(L-1)$ and $R \in \Stopprime(L,k)$ such that $R \subseteq Q$. Then, there exists a set $X \subseteq \envelope(k)$ which is Ahlfors $d$-regular, supports a $(1,1)$-Poincar\'e inequality, and contains $b_{Q,k} \cup b_{R,k}$.
\end{cor}

\begin{proof}[Proof of Corollary \ref{l:bridge-up-and-down} assuming Proposition \ref{p:PI-in-layers}]
	Let $S = S_L = T_L \in \calT^L(k)$ be the stopping-time regions such that $R \subseteq Q(S) \subseteq Q$. The set $X$ is obtained by applying Proposition \ref{p:PI-in-layers} with this data (note the the intersection property is trivial). 
\end{proof}

We begin the proof of Proposition \ref{p:PI-in-layers} in earnest. Let $k \in \N$, $L \geq 1$ and $S_L,T_L \in \calT^L(k)$ as in the statement of Proposition \ref{p:PI-in-layers}, and fix some arbitrary $Q_L \in \min(S_L) \cap \Stopprime(L,k)$ and $R_L \in \min(T_L) \cap \Stopprime(L,k)$. To ease notation, we will drop the dependency of $k$ on our notation for the subset. For example, we will write 
\begin{align}
	\surface{S} = \surface{S,k}, \ \surfacelift{S} = \surfacelift{S,k}  \mbox{ and } b_{Q} = b_{Q,k}. 
\end{align}

The first thing we show is that $Q(S_{L-1})$ and $Q(T_{L-1})$ are sufficiently close (in terms of Proposition \ref{p:Gamma}).

\begin{lem}\label{e:big-intersection-above}
	We have 
	\begin{align}
		2 B_{Q(S_{L-1})} \cap 2 B_{Q(T_{L-1})} \neq \emptyset.
	\end{align}
	In particular, there exists $\Gamma_{L-1} \subseteq \surface{L-1} \cap \Asecond B_{Q(S_{L-1})}$ with 
	\begin{align}
		\diam(\Gamma_{L-1}) \lesssim \min\{\ell(Q(S_{L-1})) , \ell(Q(T_{L-1}))\}
	\end{align}
	which is Ahlfors $d$-regular, supports a $(1,1)$-Poincar\'e inequality and satisfies 
	\begin{align}
		\surface{S_{L-1}} \cup  \surface{T_{L-1}} \cup  \subseteq \Gamma_{L-1}. 
	\end{align}
\end{lem}

\begin{proof}
	The existence of the set $\Gamma_{L-1} \subseteq \surface{L-1} \cap \Asecond B_{Q(S_{L-1})}$ (which implies the diameter estimate by Lemma \ref{l:comp-size}) satisfying the above properties follows directly from Proposition \ref{p:Gamma} once we prove the stated intersection property. Since $Q_{L-1} \in \Stopprime(L-1,k)$ and $Q(S_{L})  \subseteq Q_{L-1} \subseteq Q(S_{L-1})$, it follows from Lemma \ref{l:N-and-k-com} that
	\begin{align}
		\ell(Q(S_{L})) \leq \ell(Q_{L-1}) \lesssim \taustop \ell(Q(S_{L-1})). 
	\end{align}
	Thus, by the triangle inequality and the fact that $Q(S_L) \subseteq Q(S_{L-1})$, we have
	\begin{align}
		\Afirst B_{Q(S_{L})} \subseteq B(x_{Q(S_{L-1})}, \ell(Q(S_{L-1})) + \Afirst \ell(Q(S_{L}))) \subseteq (1+C\taustop)B_{Q(S_{L-1})}
	\end{align}
	for some $C = C(\Afirst)$. For $\taustop$ sufficiently small, depending on $\Afirst$, this implies that 
	\begin{align}
		\Afirst B_{Q(S_{L})} \subseteq 2B_{Q(S_{L-1})}. 
	\end{align}
	Similarly, 
	\begin{align}
		\Afirst B_{Q(T_{L})} \subseteq 2B_{Q(T_{L-1})} .
	\end{align}
	Since $\Afirst B_{Q(S_{L})}$ and $\Afirst B_{Q(T_{L})}$ have non-empty intersection, this proves the lemma. 
\end{proof}

By Lemma \ref{e:big-intersection-above}, we are able to apply Proposition \ref{p:Gamma} to construct an Ahlfors regular PI space which covers the set corresponding to level $L-1$. More precisely, we have the following.

\begin{lem}\label{l:X_{L-1}}
	If $\Gamma_{L-1}^*$ denotes the graph of $\taustop d_{L-1,k}$ over $\Gamma_{L-1}$, then the set 
	\begin{align}
		X_{L-1} \coloneqq \Gamma_{L-1}^* \cup b_{Q_{L-1}} \cup b_{R_{L-1}} 
	\end{align}
	is Ahlfors $d$-regular, supports a $(1,1)$-Poincar\'e inequality and has diameter at most some constant multiple of $\min\{\ell(Q(S_{L-1})) , \ell(Q(T_{L-1}))\}$. 
\end{lem}

\begin{proof}
	To simplify notation, let us drop the dependency on $L$. For example, we will denote $\Gamma = \Gamma_{L-1}$ and $b_{Q} = b_{Q_{L-1}}$. The Ahlfors regularity of $X_{L-1}$ follows directly from Lemma \ref{l:ADR-gluing} since each $b_{Q'}$ and $\Gamma_{L-1}^*$ is Ahlfors $d$-regular. We begin with the diameter estimate. Since $\Gamma^*$ is the graph a Lipschitz function with small constant over $\Gamma$, it follows from Lemma \ref{e:big-intersection-above} that $\diam(\Gamma^*) \lesssim \min\{\ell(Q(S_{L-1})) , \ell(Q(T_{L-1}))\}$. By Proposition \ref{p:b_Q} \eqref{prop:b_Q-ADR-PI}, we have $\diam(b_{Q}) \lesssim \ell(Q) \leq \ell(Q(S))$ and, similarly, $\diam(b_{R}) \lesssim \ell(Q(T))$. Since $Q(S)$ and $Q(T)$ have comparable sizes (Lemma \ref{l:comp-size}) and $X_{L-1}$ is connected, these estimates imply the desired diameter estimate.

	We now prove the Poincar\'e inequality. Since $\Gamma$ supports a $(1,1)$-Poincar\'e inequality, the same is true of $\Gamma^*$ with comparable constants (since they are bi-Lipschitz equivalent with constant close to 1). Recall from Proposition \ref{p:b_Q} \eqref{prop:b_Q-ADR-PI} that both $b_{Q}$ and $b_{R}$ support a $(1,1)$-Poincar\'e inequality. We start by showing that 
	\begin{align}\label{e:one-gluing}
		\tilde{X} \coloneqq \Gamma^* \cup b_Q \mbox{ supports a } (1,1)\mbox{-Poincar\'e inequality.} 
	\end{align}
	Indeed, since $\Gamma^* \supseteq \surfacelift{S}$, it follows from Proposition \ref{p:b_Q} \eqref{prop:b_Q-ADR-PI} and \eqref{prop:b_Q-upper} that $\Gamma^* \cap b_Q$ is lower Ahlfors $d$-regular on all balls with radius at most $\ell(Q) \sim \diam(b_Q)$. Moreover, by Proposition \ref{p:b_Q} \eqref{prop:b_Q-upper} and \eqref{prop:b_Q-gluing}, 
	\begin{align}
		\inf_{a \in \Gamma^* \cap b_Q} |x-a| + |a-y| \lesssim |x-y|
	\end{align}
	for all $x \in b_Q$ and $y \in \Gamma^*$. Thus, \eqref{e:one-gluing} follows immediately from Lemma \ref{l:gluing}. 
	
	We now check that 
	\begin{align}
		X = \tilde{X} \cup b_R \mbox{ supports a $(1,1)$-Poincar\'e inequality.} 
	\end{align}
	Again, we will use Lemma \ref{l:gluing}. We begin by noting that, by Proposition \ref{p:b_Q} \eqref{prop:b_Q-subset}, the bridges satisfy 
	\begin{align}\label{e:bridges-are-far}
		\dist(b_Q , b_R) \gtrsim  \ell(R) \sim \diam(b_R). 
	\end{align} 
	In particular, the bridges are disjoint and (for the same reasons as above) this implies that $\tilde{X} \cap b_R$ is lower Ahlfors $d$-regular on all balls with radius at most $\ell(R) \sim \diam(b_R)$. The only thing left to show is that 
	\begin{align}
		\inf_{a \in \Gamma^*} |x-a| + |a-y| \lesssim |x-y|
	\end{align}
	for all $x \in b_R$ and $y \in \tilde{X}$. As in the proof of \eqref{e:one-gluing}, if $y \in \Gamma^*$ then this follows directly from Proposition \ref{p:b_Q} \eqref{prop:b_Q-upper} and \eqref{prop:b_Q-gluing}. If $y \in b_Q$, this follows from Lemma \ref{l:diameter-PI}. Indeed, by \eqref{e:bridges-are-far}, $|x-y| \gtrsim  \diam(b_R)$. 
\end{proof}

We now want to find a set which covers the sets corresponding to level $L$ and which interacts well (in terms of a Poincar\'e inequality) with the set $\Gamma_{L-1}^*$. The construction of this set (and the proof of Proposition \ref{p:PI-in-layers}) splits into two cases. Recall from Proposition \ref{p:Gamma} that there are two possibilities for the relationship between $\Gamma_{S_L}$ and $\Gamma_{T_L}$. Either they are contained in some set which supports a Poincar\'e inequality (within the same layer) or they are well-separated. We deal now with the first case.   \\

\noindent\textbf{Case 1.} Suppose that there exists a set $\Gamma_{L} \subseteq \surface{L} \cap \Asecond B_{Q(S_L)}$ which contains $\Gamma_{S_L} \cup \Gamma_{T_L}$, is Ahlfors $d$-regular and supports a $(1,1)$-Poincar\'e inequality. \\

We start by noting the following, whose proof is the same (up to notational changes) as the proof of Lemma \ref{l:X_{L-1}}. 

\begin{lem}\label{l:X_{L}}
	If $\Gamma_L^*$ denotes the graph of $\tau d_{L,k}$ over $\Gamma_L$, then the set 
	\begin{align}
		X_{L} \coloneqq \Gamma_L^* \cup b_{Q_{L}} \cup b_{R_{L}} 
	\end{align}
	is Ahlfors $d$-regular, supports a $(1,1)$-Poincar\'e inequality and satisfies $\diam(X_L) \lesssim \min\{\ell(Q(S_L)),\ell(Q(T_L))\}$.  
\end{lem}

\begin{proof}[Proof of Proposition \ref{p:PI-in-layers} in Case 1.]
	We set 
	\begin{align}
		X \coloneqq X_{L-1} \cup X_L.
	\end{align}
	The inclusion \eqref{e:PI-in-layers} is immediate from the definition of $X_{L-1}$ and $X_L$. The Ahlfors regularity of $X$ follows from Lemma \ref{l:ADR-gluing} since each of $X_{L-1}$ and $X_L$ are Ahlfors regular with non-trivial intersection. To check the Poincar\'e inequality we will (as usual) apply Lemma \ref{l:gluing}. We show first that 
	\begin{align}\label{e:X_L-LR}
		X_{L-1} \cap X_L \mbox{ is lower Ahlfors $d$-regular}. 
	\end{align}	
	Indeed, recall from Proposition \ref{p:Gamma} and Lemma \ref{l:X_{L}} that 
	 \begin{align}
		\Gamma_L^* \supseteq \Gamma^*_{S_L} \supseteq \surfacelift{U_{S_L}}.
	\end{align}
	Thus, it follows from Proposition \ref{p:b_Q} \eqref{prop:b_Q-ADR-PI} and \eqref{prop:b_Q-lower} that $X_L \cap b_{Q_{L-1}}$ is lower Ahlfors $d$-regular on all balls with radius at most $\ell(Q_{L-1}) \sim \ell(Q(S_L))$. By similar reasoning, $X_L \cap b_{R_{L-1}}$ is lower Ahlfors $d$-regular on all balls with radius at most $\ell(R_{L-1}) \sim \ell(Q(T_L))$. Since $\diam(X_L) \lesssim \min\{\ell(Q_{L-1}) , \ell(R_{L-1})\}$ (by Lemma \ref{l:X_{L}}) and $X_{L-1} \cap X_L$ is contained in $b_{Q_{L-1}} \cup b_{R_{L-1}}$, the lower regularity of this set on balls with radius at most $\diam(X_L)$ is implied by the lower regularity estimates we established above. 
	
	Now we check that 
	\begin{align}\label{e:glue-X_L}
		\inf_{a \in X_{L-1} \cap X_L} |x-a| + |a-y|  \lesssim |x-y|
	\end{align}
	for all $x \in X_{L-1}$ and $y \in X_L$. If $x \in b_{Q_{L-1}}$ or $x \in b_{R_{L-1}}$, then this estimate follows from Proposition \ref{p:b_Q} \eqref{prop:b_Q-gluing} since, by Proposition \ref{p:b_Q} \eqref{prop:b_Q-subset} and the definition of $\Gamma_L^*$, we have $\pi_1(y) \leq \taustop d_{L,k}(\pi_n(x))$. Suppose instead that $x \in \Gamma_{L-1}^*$. Since 
	\begin{align}\label{e:n-y}
		\pi_n(y) \in \Asecond B_{Q(S_L)} \subseteq \Asecond B_{Q_{L-1}},
	\end{align}
	if $\pi_n(x) \notin 2\Asecond B_{Q_{L-1}}$, then we have $|x-y| \gtrsim \ell(Q(S_L)) \sim \diam(X_L)$ and the estimate follows from Lemma \ref{l:diameter-PI}. Suppose instead that $\pi_n(x) \in 2\Asecond B_{Q_{L-1}}$. Then, recalling \eqref{e:n-y} and applying Lemma \ref{l:well-separated}, we have 
	\begin{align}
		\pi_1(x) - \pi_1(y) &\geq \taustop d_{L-1,k}(\pi_n(x)) - \taustop d_{L,k}(\pi_n(y)) \\
		&\gtrsim  \ell(Q_{L-1}) \sim \diam(X_L),
	\end{align}
	and the estimate again follows from Lemma \ref{l:diameter-PI}. Properties \eqref{e:X_L-LR} and \eqref{e:glue-X_L} imply that $X_{L-1}$ and $X_L$ satisfies the hypotheses of Lemma \ref{l:gluing}. Thus, we reach our desired conclusion. 
\end{proof}

\noindent\textbf{Case 2.} Suppose that there is no set satisfying the hypotheses of Proposition \ref{p:Gamma}. In which case, 
\begin{align}\label{e:S,T-far}
	\dist(\Gamma_{S_L}, \Gamma_{T_L}) \gtrsim \ell(Q(S)) + \ell(Q(T)). 
\end{align}

We define 
\begin{align}
	X_L \coloneqq \Gamma_{S_L}^* \cup \Gamma_{T_L}^*\cup b_{Q_L} \cup b_{R_L}. 
\end{align}
Then, the following lemma implies Proposition \ref{p:PI-in-layers} in this case. 

\begin{lem}
	The set $X \coloneqq X_{L-1} \cup X_L$ is Ahlfors $d$-regular and supports a $(1,1)$-Poincar\'e inequality. 
\end{lem}

\begin{proof}
	The inclusion \eqref{e:PI-in-layers} is immediate from the definitions of $X_{L-1}$ and $X_L$.  It follows from the definitions that $X$ is the union of a sequence of Ahlfors regular spaces with non-trivial intersections. Thus, by Lemma \ref{l:ADR-gluing}, it is Ahlfors regular. Since $X_L$ does not support a Poincar\'e inequality, we cannot directly apply the gluing lemma (Lemma \ref{l:gluing}). As in the proof of Lemma \ref{l:X_{L-1}}, we will glue $X_L$ to $X_{L-1}$ in two steps. We begin by showing that 
	\begin{align}\label{e:PI-case-2}
		\tilde{X} \coloneqq X_{L-1} \cup \Gamma_{S_L}^* \cup b_{Q_L} \mbox{ supports a $(1,1)$-Poincar\'e inequality}. 
	\end{align}
	
	We start with the lower regularity of $X_{L-1} \cap (\Gamma_{S_L}^* \cup b_Q)$. Note to begin with that 
	\begin{align}
		X_{L-1} \cap (\Gamma_{S_L}^* \cup b_{Q_L}) \subseteq \surfacelift{L}. 
	\end{align}
	Since the only parts of $X_{L-1}$ which have non-empty intersection with $\surfacelift{L}$ are $b_{Q_{L-1}}$ and $b_{R_{L-1}}$, and since (by Lemma \ref{p:Gamma} and Proposition \ref{p:b_Q} \eqref{prop:b_Q-upper}) we have $b_{Q_L} \cap \surfacelift{L} \subseteq \Gamma_{S_L}^*$, it follows that 
	\begin{align}
		X_{L-1} \cap (\Gamma_{S_L}^* \cup b_{Q_L}) \subseteq \Gamma_{S_L}^* \cap (b_{Q_{L-1}} \cup b_{R_{L-1}}). 
	\end{align}
	By Proposition \ref{p:Gamma} and Proposition \ref{p:b_Q} \eqref{prop:b_Q-lower}, we have $b_{R_{L-1}} \cap \surfacelift{L} \subseteq \surfacelift{U_{T_L}} \subseteq \Gamma_{T_L}^*$ which, since $\Gamma_{S_L}^*$ and $\Gamma_{T_L}^*$ are well-separated, implies that 
	\begin{align}
		X_{L-1} \cap (\Gamma_{S_L}^* \cup b_{Q_L}) = \Gamma_{S_L}^* \cap b_{Q_{L-1}}. 
	\end{align}
	The lower regularity now follows directly from Proposition \ref{p:b_Q} \eqref{prop:b_Q-lower}. To establish \eqref{e:PI-case-2}, the only other thing one needs to check is that 
	\begin{align}\label{e:PI-tilde-X}
		\inf_{a \in X_{L-1} \cap (\Gamma_{S_L}^* \cup b_{Q_L})} |x-a| + |a-y|  \lesssim |x-y|
	\end{align}
	for all $x \in X_{L-1}$ and $y \in \Gamma_{S_L}^* \cup b_{Q_L}$. The proof of \eqref{e:PI-tilde-X} (up to notational changes) identical to the proof of \eqref{e:glue-X_L}, and we will not repeat the argument. 
	
	We now finish by checking that 
	\begin{align}
		\mbox{$X = \tilde{X} \cup \Gamma_{T_L}^* \cup b_{R_L}$ supports a $(1,1)$-Poincar\'e inequality.} 
	\end{align}
	By the same argument as above, we have that 
	\begin{align}
		X_{L-1} \cap (\Gamma_{T_L}^* \cup b_{R_L}) = \Gamma_{T_L}^* \cap b_{R_{L-1}} .
	\end{align}
	By Proposition \ref{p:b_Q} \eqref{prop:b_Q-lower}, this set is lower Ahlfors $d$-regular. Thus, the only thing left to show is that 
	\begin{align}
		\inf_{a \in X_{L-1} \cap (\Gamma_{S_L}^* \cup b_{Q_L})} |x-a| + |a-y|  \lesssim |x-y|
	\end{align}
	for all $x \in \tilde{X}$ and $y \in \Gamma_{T_L}^* \cup b_{R_L}$. If $x \in \tilde{X}$ then the same proof as in \eqref{e:glue-X_L} gives the desired estimate. So, let us suppose that $x \in \Gamma_{S_L}^* \cup b_{Q_L}$. In this case, it follows from \eqref{e:S,T-far} that $|x-y| \gtrsim \ell(Q(S_{L})) \sim \diam(\Gamma_{S_L} \cup b_{Q_L})$, and the estimate follows from Lemma \ref{l:diameter-PI}. 
\end{proof}

 \section{Global Poincar\'e inequalities across all layers}\label{s:global}

We have now reached the point where we can prove the desired properties of the set $\envelope(k)$ (Ahlfors regularity, Poincare inequality, approximation of $E$). The precise statement that we will prove is the following. 

\begin{prop}\label{p:Z(k_0)}
	For each $k \in \N$, the set $Z(k)$ is closed, contained in the closed ball $\overline{B}(0,15\sqrt{n}\lambdadyadic)$ in $\R^{n+1}$, Ahlfors $d$-regular and supports a $(1,1)$-Poincar\'e inequality with constant that depend only on $n,d$ and the uniformly rectifiability constants i.e., independent of $k$. Moreover, we have 
	\begin{align}\label{e:E-Z(k_0)}
		\dist(x,Z(k)) \lesssim \rho^{k} 
	\end{align}
	for each $x \in E$. 
\end{prop}

Before proving Proposition \ref{p:Z(k_0)}, let us see why it implies Theorem \ref{t:main}. 

\begin{proof}[Proof of Theorem \ref{t:main} assuming Proposition \ref{p:Z(k_0)}]
	Since the space of non-empty compact subsets of $\overline{B}(0,15\sqrt{n}\lambdadyadic)$ is compact with respect to Hausdorff convergence, there is a subsequence $\{Z(k_i)\}_{i\in\N}$ converging to a compact set $Z \subseteq \overline{B}(0,15\sqrt{n}\lambdadyadic)$. Since each $Z(k_i)$ is Ahlfors $d$-regular and supports a $(1,1)$-Poincar\'e inequality with uniform constants, the same is true of $Z$. Finally, \eqref{e:E-Z(k_0)} implies $E \subseteq Z$.
\end{proof}

Fix some $k \in \N$. We now prove Proposition \ref{p:Z(k_0)}.  

\begin{lem}
	The set $Z(k)$ is closed and contained in the closed ball $\overline{B}(0,15\sqrt{n}\lambdadyadic)$ in $\R^{n+1}$. 
\end{lem}

\begin{proof}
	It is clear that $\envelope(k)$ is closed since it is the union of finitely many closed sets. To prove the inclusion, fix some the moment some $N \in \N$ and $S \in \calT^N(k)$. By Proposition \ref{p:UR}, we know that $\pi_n(F_S) \subseteq 2\sqrt{n}\lambda B_{Q(S)}$. Since, $0 \leq \pi_1(x) \leq \taustop d_{N,k}(\pi_n(x))$ for each $x \in F_S$, this implies that $F_S \subseteq 3\sqrt{n}\lambda B_{Q(S)}$ (as a ball in $\R^{n+1}$). Recalling \eqref{e:largest0-ball}, this implies 
	\begin{align}
		F_S \subseteq 3\sqrt{n}\lambdadyadic B_{Q(S)} \subseteq 3\sqrt{n}\lambdadyadic B_{Q_0} \subseteq B(0,15\sqrt{n}\lambdadyadic ). 
	\end{align}
	Since $\envelope(k)$ is just the union over the $F_S$, this finishes the proof. 
\end{proof}

\begin{lem}\label{l:lower-reg-Z}
	The set $\envelope(k)$ is lower Ahlfors $d$-regular with regularity constant depending only on $n,d$ and the uniform rectifiability constants for $E$.
\end{lem}

\begin{proof}
	Fix $x \in \envelope(k)$ and $r \in (0,\diam(\envelope(k)))$. Moreover, suppose that
	\begin{align}
		x \in \surfacelift{S,k} \cup b_{Q,k} 
	\end{align}
	for some $N \in \N$, $S \in \calT^N(k)$ and $Q' \in \Stopprime(N,k)$. We can suppose without loss of generality that $Q' \subseteq Q(S)$. Let $Q^* \in \Stop(N-1)$ such that $Q(S) \subseteq Q^*$. 
	
	Let $\lambda > 1$ (to be chosen sufficiently large). If $r \leq \lambda \ell(Q^*) \sim_{\lambda} \diam(\Gamma_{S,k}^* \cup b_{Q,k})$, then the lower regularity estimates follow directly from the Ahlfors regularity of this space (recall that $\Gamma_{S,k}^*$ and $b_{Q,k}$ are Ahlfors regular by Propositions \ref{p:Gamma} and \ref{p:b_Q}, respectively). So, we suppose instead that 
	\begin{align}\label{e:choice-r}
		r >  10\lambda \ell(Q^*).
	\end{align}
	Let $Q^{**}$ be the largest cube in the set 
	\begin{align}
		\left\{ R \in \bigcup_{L =0}^{N-1} \Stop(L) \colon  R \supseteq Q \mbox{ and } \lambda \ell(R) \leq \pi_1(x) + r/10 \right\}
	\end{align}
	(which is non-empty by \eqref{e:choice-r}). Let $L \in \{0,\dots,N-1\}$ and $T \in \calT^L(k)$ such that $Q^{**} \in T$. By maximality, the parent $R$ of $Q(T)$ satisfies $\lambda \ell(R) > \pi_1(x) + r/10$. In particular, 
	\begin{align}\label{e:diam-est}
		\diam(\Gamma_{T,k}^*) \sim \diam(\surfacelift{T,k}) \sim \ell(Q(T)) \sim \ell(R) \gtrsim r. 
	\end{align} Let $x^* \in \surfacelift{T,k}$ such that 
	\begin{align}\label{e:x^*}
		|\pi_n(x_{Q^{**}}) - \pi_n(x^*)| \leq 2\ell(Q^{**}),
	\end{align}
	which exists by Lemma \ref{l:close-point}. We claim that this point satisfies 
	\begin{align}\label{e:inclusion}
		B(x^* , r/2) \subseteq B(x,r). 
	\end{align} 
	Since $\Gamma^*_{T,k}$ is Ahlfors regular (Proposition \ref{p:Gamma}) with diameter at least some constant multiple of $r$ (by \eqref{e:diam-est}), this is sufficient for the lower regularity estimates. We commence with the proof of \eqref{e:inclusion} in earnest. First, by Lemma \ref{l:contained-in-ball} and Proposition \ref{p:b_Q} \eqref{prop:b_Q-subset} we know that 
	\begin{align}\label{e:proj}
		\pi_n(\surfacelift{S,k} \cup b_{Q,k}) \subseteq 2\sqrt{n}\lambdadyadic B_{Q(S)} \subseteq 2\sqrt{n}\lambdadyadic B_{Q^{**}} 
	\end{align}
	and $\pi_1(x) \leq \taustop d_{N,k}(\pi_n(x))$. Thus, by \eqref{e:stop-size-cutoff-2} and Lemma \ref{l:monoton-in-N}, we have  
	\begin{align}
		\begin{split}\label{e:pi_1(x)}
			0 \leq \pi_1(x) \leq \lambda  \ell(Q^{**}) \leq r/10
		\end{split}
	\end{align}
	for $\lambda$ sufficiently large. On the other hand, it follows from \eqref{e:x^*} and Lemma \ref{l:N-and-k-com} that we have $\pi_1(x^*) \sim \taustop d_{L,k}(\pi_n(x^*)) \sim \ell(Q^{**})$. Thus, for $\lambda$ sufficiently large, 
	\begin{align}\label{e:pi_1(x^*)}
		0 \leq \pi_1(x^*) \leq \lambda  \ell(Q^{**}) \leq \pi_1(x) + r/10 .
	\end{align}
	Inequalities \eqref{e:pi_1(x)} and \eqref{e:pi_1(x^*)} give 
	\begin{align}\label{e:pi-1}
		|\pi_1(x) - \pi_1(x^*)| \leq r/10. 
	\end{align}
	Furthermore, by \eqref{e:x^*} and \eqref{e:proj},
	\begin{align}
		\begin{split}\label{e:pi-n}
			|\pi_n(x) - \pi_n(x^*)| &\leq |\pi_n(x) - \pi_n(x_{Q^**})| + |\pi_n(x_{Q^{**}}) - \pi_n(x^*)| \\
			&\leq (2\sqrt{n}\lambdadyadic + 2) \ell(Q^{**}) \leq r/10
		\end{split}
	\end{align}
	after choosing $\lambda$ sufficiently large. Inequalities \eqref{e:pi-1} and \eqref{e:pi-n} imply $|x - x^*| \leq r/5$. The inclusion \eqref{e:inclusion} is an easy consequence of this estimate. 
\end{proof}

\begin{lem}
	The set $\envelope(k)$ is upper Ahlfors $d$-regular with regularity constant depending only on $n,d$ and the uniform rectifiability constants for $E$.
\end{lem}

\begin{proof}
	Fix $x \in \envelope(k)$ and $r \in (0,\diam(\envelope(k)))$, and let $N \in \N$, $S \in \calT^N(k)$ and $Q' \in \Stopprime(N,k)$ as in the proof of Lemma \ref{l:lower-reg-Z}. For each $T \in \calT(k)$, consider the set 
	\begin{align}
		F_T \coloneqq \surfacelift{T,k} \cup \bigcup_{\substack{Q \in \Stopprime(N,k) \\ Q \subseteq Q(T)}} b_{Q,k}.
	\end{align}
	Then, let 
	\begin{align}
		\calT_r(k) \coloneqq \{T \in \calT(k) \colon F_T \cap B(x,r) \neq \emptyset\}.
	\end{align}
	By Proposition \ref{p:UR}, each $F_T$ is upper Ahlfors $d$-regular with diameter at most some constant multiple of $\ell(Q(T))$. In particular, 
	\begin{align}\label{e:upper-Z}
		\calH^d(\envelope(k) \cap B(x,r)) &\leq \sum_{T \in \calT_r(k)} \calH^d(F_T \cap B(x,r)) \\
		&\lesssim \sum_{T \in \calT_r(k)} \min\{\ell(Q(T)),r\}^d.
	\end{align}
	To deal with this final sum, we split it into two pieces. Indeed, let $\lambda > 1$ (to be chosen sufficiently large) and consider the sets 
	\begin{align}
		\calT_{\text{small}}(k) &\coloneqq \{T \in \calT_r(k) \colon  \ell(Q(T)) \leq \lambda r \}; \\
		\calT_{\text{large}}(k) &\coloneqq \{T \in \calT_r(k) \colon \ell(Q(T)) > \lambda r \}.
	\end{align}
	To deal with the sum over cubes in $\calT^N_{\text{small}}(k)$, first recall from Proposition \ref{p:UR} that $\pi_n(F_T) \subseteq 2\sqrt{n}\lambdadyadic B_{Q(T)}$. Thus, $2\sqrt{n}\lambdadyadic B_{Q(T)} \cap \pi_n(B(x,r)) \neq \emptyset$. For each $T \in \calT^N_{\text{large}}(k)$, this means that $B_{Q(T)} \subseteq B(\pi_n(x),R)$ for some $R$ depending on $\lambda$. Thus, since $\calT$ satisfies a Carleson packing condition (Lemma \ref{l:Carleson}), we have  
	\begin{align}\label{e:small}
		\sum_{T \in \calT_{\text{small}}(k)} \ell(Q(T))^d \lesssim r^d.
	\end{align}
	The upper regularity of $\envelope(k)$ in $B(x,r)$ follows from \eqref{e:upper-Z} and \eqref{e:small} once we show that 
	\begin{align}\label{e:large}
		\sum_{T \in \calT_{\text{large}}(k)} r^d \lesssim r^d. 
	\end{align}
	To prove \eqref{e:large}, consider first the set $\calN$ of integers $N \in \N$ such that 
	\begin{align}
		\calT^N_{\text{large}}(k) \coloneqq \calT^N(k) \cap \calT_{\text{large}}(k) \neq \emptyset. 
	\end{align}
	Since $2\sqrt{n}\lambdadyadic B_{Q(T)} \cap \pi_n(B(x,r)) \neq \emptyset$ and $\ell(Q) > \lambda r$ for each $T \in \calT_{\text{large}}(k)$, it follows from Lemma \ref{l:comp-size} and a simple volume argument that 
	\begin{align}\label{e:T-N}
		|\calT^N_{\text{large}}(k)| \lesssim_{d,n,C_0} 1 \quad  \mbox{ for all } N \in \calN. 
	\end{align}
	We claim that 
	\begin{align}\label{e:N}
		|\calN| \leq 2,
	\end{align}
	which, along with \eqref{e:T-N}, is enough to imply \eqref{e:large}. Let $N$ be the largest integer in $\calN$, $T \in \calT^N_{\text{large}}(k)$ and $Q \in \Stopprime(N-1,k)$ such that $Q(T) \subseteq Q$. Moreover, fix a point
	\begin{align}
		z \in F_{T} \cap B(x,r). 
	\end{align}
	Note that, since $\pi_n(F_T) \subseteq 2\sqrt{n}\lambdadyadic B_{Q(T)}$ and $\ell(Q(T)) > \lambda r$, we have 
	\begin{align}\label{e:TTTT}
		\pi_n(B(x,r)) \subseteq 6\sqrt{n}\lambdadyadic B_{Q(T)} \subseteq 6\sqrt{n} \lambdadyadic B_{Q} . 
	\end{align}
	Suppose towards a contradiction that there exists $M \leq N-2$, $S \in \calT^M_{\text{large}}(k)$ and $w  \in F_S \cap B(x,r)$. By \eqref{e:TTTT}, $\pi_n(w),\pi_n(z) \in 6\sqrt{n}\lambda B_{Q}$. Then, since $\pi_1(z) \leq \taustop d_{N,k}(\pi_n(z))$ and $\pi_1(w) \geq \taustop d_{M+1}(\pi_n(w))$ (and $M+1 \leq N-2$), it follows from Lemma \ref{l:well-separated} that 
	\begin{align}
		\pi_1(w) - \pi_1(z) \geq c \ell(Q) \geq c\lambda r > r
	\end{align}
	for $\lambda$ sufficiently large, which is a contradiction. We have just argued that $\calN \subseteq \{N,N-1\}$, which clearly implies \eqref{e:N}. 
\end{proof}

\begin{lem}
	Inequality \eqref{e:E-Z(k_0)} holds. 
\end{lem}

\begin{proof}
	Let $x \in E$ and let $Q \in \cubes_{k}$ such that $x \in Q$. By construction, there exists $N \in \N$ and $S \in \calT^N(k)$ such that $Q \in \min(S(k))$. Thus, by Lemma \ref{l:close-point}, there exists $z \in \surface{S,k}$ such that $|x-z| \leq 2\ell(Q)$. Moreover, this point satisfies $\taustop d_{N,k}(z) \lesssim \ell(Q)$ (see Lemma \ref{l:stop-size-cutoff}). Letting $z^* \coloneqq (z,\tau d_{N,k}(z)) \in \surfacelift{S}$, we now have 
	\begin{align}
		|x - z^*| \leq |x-z| + |z-z^*| \lesssim \ell(Q) \lesssim \rho^k,
	\end{align}
	as required. 
\end{proof}

The only thing left to show is that $Z(k)$ supports a $(1,1)$-Poincar\'e inequality. Fix some $x,y \in Z(k)$ and $u \colon Z(k) \to \R$ Lipschitz. Our goal is to prove \eqref{e:PI-max} with $C_1$ independent of $x,y$ and $u$. Let $N \geq 0$, $S \in \calT^N(k)$ and $Q_N \in \Stopprime(N,k)$ such that $Q_N \subseteq Q(S)$ and 
\begin{align}\label{e:x-in-S}
	x \in \surfacelift{S} \cup b_{Q_N},
\end{align}
and let $M \geq 0$, $T \in \calT^N(k)$ and $R_M \in \Stopprime(M,k)$ such that $R_M \subseteq Q(T)$ and 
\begin{align}\label{y-in-T}
	y \in \surfacelift{T} \cup b_{R_M}.
\end{align}
We can assume, without loss of generality, that 
\begin{align}
	M \leq N  .
\end{align}
Note that 
\begin{align}\label{e:heights}
	\pi_1(x) \leq \taustop d_{N,k}(\pi_n(x)) \quad \mbox{ and } \quad \pi_1(y) \geq \taustop d_{M+1,k}(\pi_n(y)).
\end{align}
Let $Q_0 \in \cubes$ with smallest cube such that $\pi_n(x) , \pi_n(y) \in 3\sqrt{n}\lambdadyadic B_{Q_0}$ (recall the definition of $\lambdadyadic$ from \eqref{e:J-contains-ball}). Let $K \geq 0$ and $U \in \calT^K$ such that $Q_0 \in U$. 

The proof of \eqref{e:PI-max} will be split into a number of sub-cases. In two of these cases, we check \eqref{e:PI-max} by constructing a \textit{pencil of balls} from $x$ to $y$ satisfying the hypotheses of Lemma \ref{l:local-to-global}. The balls will a subset of the following set.  

\begin{lem}\label{l:stop-balls}
	For each $Q \in \Stop(N,k)$ there is a ball $\mathfrak{B}_Q$ centred on $b_Q$ with radius comparable to $\ell(Q)$ and such that 
	\begin{align}
		\envelope(k) \cap \mathfrak{B}_Q = b_Q \cap \mathfrak{B}_Q.  
	\end{align}
\end{lem}

The sets in Lemma \ref{l:local-to-global} (which contain these balls) will be constructed by Proposition \ref{p:PI-in-layers}. The following will be useful in verifying the hypotheses of that lemma. 

\begin{lem}\label{l:tops-are-nearby}
	Let $L \leq K$ and $V,W \in \calT^L(k)$ such that $Q_N \subseteq Q(V)$ and $R_M \subseteq Q(W)$. Then, $\Afirst B_{Q(V)} \cap \Afirst B_{Q(W)} \neq \emptyset$. 
\end{lem}

\begin{proof}
	Recall from Lemma \ref{l:contained-in-ball} that $\pi_n(x) \in 2\sqrt{n}\lambdadyadic B_{Q(S)}$. By definition, we have $Q(S) \subseteq Q(V)$. Thus, it follows from the triangle inequality that $\pi_n(x) \in 3\sqrt{n}\lambdadyadic B_{Q(V)}$. In particular, $3\sqrt{n}\lambdadyadic B_{Q(V)} \cap 3\sqrt{n}\lambdadyadic B_{Q_0} \neq \emptyset$. Applying Lemma \ref{l:comp-size} and the triangle inequality, it is easy to check that 
	\begin{align}
		3\sqrt{n}\lambdadyadic B_{Q_0} \subseteq 3(2+\Ccomp)\sqrt{n}\lambdadyadic B_{Q(V)} \subseteq  \Afirst B_{Q(V)}.
	\end{align}
	Since $3\sqrt{n}\lambdadyadic B_{Q(W)} \subseteq A_1B_{Q(W)}$ and $3\sqrt{n}\lambdadyadic B_{Q(W)} \cap 3\sqrt{n}\lambdadyadic B_{Q_0} \neq \emptyset$, the result now follows. 
\end{proof}

We commence with the proof of \eqref{e:PI-max}. \\ 

\noindent\textbf{Case 1.} Suppose that $K \leq M$. \\

We construct a sequence of balls from $x$ to $y$ and subsets of $\envelope(k)$ satisfying the hypotheses of Lemma \ref{l:local-to-global}. We start with the balls. Let $\{Q_i\}_{i=K}^{N-1} \subseteq \cubes$ be such that $Q_i \in \Stop(i,k)$ satisfies $Q(S) \subseteq Q_i$, and let $\{R_j\}_{i=K}^{M-1} \subseteq \cubes$ be such that $R_j \in \Stop(j,k)$ satisfies $Q(S) \subseteq R_j$.

\begin{lem}\label{l:good-balls-case-1}
	The collection of balls 
	\begin{align}
		\left\{\mathfrak{B}_Q \colon Q \in \{Q_i\}_{i=K}^{N-1} \cup \{R_j\}_{j=K}^{M-1}\right\}
	\end{align}
	satisfy Lemma \ref{l:local-to-global} \eqref{l:local-to-global-ends} and \eqref{l:local-to-global-r}. That is, there exists $A \geq 1$ such that $x \in A\mathfrak{B}_{Q_i}$ and $y \in A \mathfrak{B}_{R_j}$ for each $i \in \{K,\dots,N-1\}$ and $j \in \{K,\dots,M-1\}$, and 
	\begin{align}\label{e:sum-radii-bounded}
		\sum_{j=K}^{N-1} r_{\mathfrak{B}_{Q_i}}  + \sum_{j=K}^{M-1} r_{\mathfrak{B}_{R_j}}  \lesssim |x-y|. 
	\end{align} 
\end{lem}

\begin{proof}
	We only deal with properties related to $\{Q_i\}$ and $x$. The properties for $\{R_i\}$ and $y$ are identical up to notational changes. Let us begin by showing that 
	\begin{align}\label{e:x-in-A}
		\mbox{$x \in A\mathfrak{B}_{Q_i}$ for some suitable choice of $A$. }
	\end{align}
	Fix some $i \in \{K,\dots,N-1\}$. For brevity, let us define $x_i \coloneqq x_{B_{Q_i}}$. First, by Lemma \ref{l:contained-in-ball}, we have $\pi_{n}(x) \in 2\sqrt{n}\lambdadyadic B_{Q(S)} \subseteq 3\sqrt{n}\lambdadyadic B_{Q_i}$. On the other hand, since $\mathfrak{B}_{Q_i} \subseteq b_{Q_i}$, we have by Proposition \ref{p:b_Q} \eqref{prop:b_Q-subset} that $\pi_n(x_i) \in (\innercube/2)B_{Q_i}$. In particular, 
	\begin{align}
		|	\pi_n(x-x_i) | \lesssim \ell(Q_i). 
	\end{align}
	Since $\pi_1(x) \leq \taustop d_{N,k}(x)$, \eqref{e:stop-size-cutoff-2} and Lemma \ref{l:monoton-in-N} implies $\pi_1(x) \lesssim \ell(Q_i)$. On the other hand, by Proposition \ref{p:b_Q} \eqref{prop:b_Q-subset} and Lemma \ref{l:N-and-k-com}, we have $\pi_1(x_i) \lesssim \taustop d_{N,k}(\pi_n(x_i)) \lesssim \ell(Q_i)$. In particular, 
	\begin{align}
		|\pi_1(x-x_i) | \lesssim \ell(Q_i). 
	\end{align}
	Recalling from Lemma \ref{l:stop-balls} that $r_{B_{Q_i}} \gtrsim \ell(Q_i)$, the above inequalities imply \eqref{e:x-in-A}.

	We now turn our attention to \eqref{e:sum-radii-bounded}. Recall from Lemma \ref{l:stop-size-cutoff} that $\ell(Q_i) \sim \taustop d_{i,k}(Q_i)$ for each $i$. Since $Q_{i+1} \subseteq Q_i$, it follows from Lemma \ref{l:N-and-k-com} that $\ell(Q_{i+1}) \lesssim \taustop \ell(Q_i)$. In particular, by summing over a geometric series, it suffices to show that 
	\begin{align}
		\ell(Q_K) \lesssim |x-y|. 
	\end{align}
	This estimate follows by minimality. Indeed, let $Q \in \text{Child}(Q_0)$ such that $\pi_n(x) \in 2\sqrt{n}\lambdadyadic B_Q$ (note that such a cube exists since $\pi_n(x) \in 2\sqrt{n}\lambdadyadic B_{Q(S)}$ by Lemma \ref{l:contained-in-ball}). By minimality, we know that $\pi_n(y) \not\in 3\sqrt{n}\lambdadyadic B_Q$. Since $\ell(Q_K) \lesssim \ell(Q)$, the estimate now easily follows. 
\end{proof}

We now construct the sets satisfying Lemma \ref{l:local-to-global} \eqref{l:local-to-global-PI}--\eqref{l:local-to-global-X}. These are constructed by applying Proposition \ref{p:PI-in-layers}. The properties Lemma \ref{l:local-to-global} \eqref{l:local-to-global-PI} and \eqref{l:local-to-global-X} are explicit in the statement below, we note that \eqref{e:bridge-in-X} and \eqref{e:bridge-in-X_K-1} are sufficient for Lemma \ref{l:local-to-global} \eqref{l:local-to-global-balls} by Lemma \ref{l:stop-balls}. 

\begin{lem}\label{l:space-case-1}
	There are sets $\{X_i\}_{i=K-1}^{N-1} , \{Y_j\}_{j=K-1}^{M-1} \subseteq \envelope(k)$, each of which is Ahlfors $d$-regular and supports a $(1,1)$-Poincar\'e inequality, such that $x \in X_{N-1}$, $y \in Y_{M-1}$ and 
	\begin{align}\label{e:bridge-in-X}
		b_{Q_i} \cup b_{Q_{i+1}} \subseteq X_{i} \quad \mbox{ and } \quad b_{R_j} \cup b_{R_{j+1}} \subseteq Y_j 
	\end{align}
	for each $i \in \{K,\dots,N-1\}$ and $j \in \{K,\dots,M-1\}$, and 
	\begin{align}\label{e:bridge-in-X_K-1}
		b_{Q_K} \cup  b_{R_K} \subseteq X_{K-1}.
	\end{align}
\end{lem}

\begin{proof}
	Except for the end points $i \in \{K-1, N-1\}$ and $j \in \{K-1, M-1\}$, we let $X_i$ and $Y_j$ be the sets constructed by Corollary \ref{l:bridge-up-and-down} with the cubes $Q_i,Q_{i+1}$ and $R_j,R_{j+1}$, respectively. In this range of indices, \eqref{e:bridge-in-X} is then immediate. For $i = N-1$ (resp. $j = M-1$), we let $X_{N-1}$ (resp. $Y_{M-1}$) be the set construct by applying Proposition \ref{p:PI-in-layers} with data $S,U_S$ and $Q_M$ (resp. $T,U_T$ and $R_M$). Then, 
	\begin{align}
		\surfacelift{S} \cup b_{Q_N} \subseteq X_{N-1} \quad \mbox{ and } \quad \surfacelift{T} \cup b_{R_{M}} \subseteq Y_{M-1}, 
	\end{align}
	which implies \eqref{e:bridge-in-X} for this choice of indices and also implies $x \in X_{N-1}$ and $y \in Y_{M-1}$ (recall \eqref{e:x-in-S} and \eqref{y-in-T}). For \eqref{e:bridge-in-X_K-1}, we also apply Proposition \ref{p:PI-in-layers}. Indeed, let $V,W \in \calT^K(k)$ such that $Q_K \in V$ and $R_K \in W$. It follows from Lemma \ref{l:tops-are-nearby} that $\Afirst B_{Q(V)} \cap \Afirst B_{Q(W)} \neq \emptyset$, and now from Proposition \ref{p:PI-in-layers} that such a set exists. 
\end{proof}

\noindent\textbf{Case 2.} Suppose that $N \geq M+2$ and $K \geq M +1$. \\

As in Case 1, we construct a sequence of balls from $x$ to $y$ and subsets of $\envelope(k)$ satisfying the hypotheses of Lemma \ref{l:local-to-global}. Let 
\begin{align}
	I = \{M,\dots,N-1\} \quad \mbox{ and } \quad I' = I \setminus \{M\}.
\end{align}
For each $i \in I$, let $Q_i \in \Stop(i)$ such that $Q(S) \subseteq Q_i$. 

\begin{lem}\label{l:geometric-2}
	The collection of balls $\{ \mathfrak{B}_{Q_i} \colon i \in I'\}$ satisfy Lemma \ref{l:local-to-global} \eqref{l:local-to-global-ends} and \eqref{l:local-to-global-r} i.e., there is $A \geq 1$ such that $x \in A\mathfrak{B}_{Q_i}$ for each $i \in I'$ and 
	\begin{align}\label{e:I'}
		\sum_{i \in I'} \ell(Q_i) \lesssim |x-y|. 
	\end{align}
\end{lem}

\begin{proof}
	Since $I' = \emptyset$ in the case that $N = M+1$, we can assume that $N \geq M+2$. The proof that $x \in A\mathfrak{B}_{Q_i}$ is the same as in the proof of Lemma \ref{l:good-balls-case-1} and we will not repeat it. Furthermore, as in the proof of Lemma \ref{l:good-balls-case-1}, \eqref{e:I'} follows once we show that 
	\begin{align}\label{e:Q_{M+1}}
		\ell(Q_{M+1}) \lesssim |x-y|. 
	\end{align}
	It suffice to prove \eqref{e:Q_{M+1}} in the case that $|\pi_n(x) - \pi_n(y)| \leq \ell(Q_{M+1})$ (otherwise there is nothing to show). In this case, we have $\pi_n(x),\pi_n(y) \in 3\sqrt{n}\lambdadyadic B_{Q_{M+1}}$. Thus, by Lemma \ref{l:well-separated} and \eqref{e:heights}, we have 
	\begin{align}
		|x-y| \geq \pi_1(y) - \pi_1(x) \gtrsim \ell(Q_{M+1}),
	\end{align}
	as required. 
\end{proof}

The properties Lemma \ref{l:local-to-global} \eqref{l:local-to-global-PI}--\eqref{l:local-to-global-X} follow directly from the next two lemmas.  

\begin{lem}
	There are sets $\{X_i\}_{i \in I'} \subseteq Z(k)$, each of which is Ahlfors $d$-regular and supports a $(1,1)$-Poincar\'e inequality, such that $x \in X_{N-1}$ and 
	\begin{align}\label{e:bridg-1}
		b_{Q_i} \cup b_{Q_{i+1}} \subseteq X_i
	\end{align}
	for each $i \in I'$.
\end{lem}

\begin{proof}
	Except for the endpoint $i = N-1$, we let $X_i$ be the set constructed by Corollary \ref{l:bridge-up-and-down} with the cubes $Q_i$ and $Q_{i+1}$. In this range of indices, \eqref{e:bridg-1} is then immediate. For $i = N-1$, we let $X_{N-1}$ be the set constructed by applying Proposition \ref{p:PI-in-layers} with the data $S,U_S$ and $Q_N$. The remaining properties hold in the same way as in Lemma \ref{l:space-case-1}. 
\end{proof}

\begin{lem}
	We have $|y - x_{\mathfrak{B}_{Q_{M+1}}}| \lesssim |x-y|$. Moreover, there exists a space $X \subseteq \envelope(k)$ which is Ahlfors $d$-regular supports a $(1,1)$-Poincar\'e inequality and satisfies 
	\begin{align}
		\{y\} \cup b_{Q_{M+1}} \subseteq X. 
	\end{align}
\end{lem}

\begin{proof}
	To simplify our notation, we will denote $z \coloneqq x_{\mathfrak{B}_{Q_{M+1}}}$. We start by estimating
	\begin{align}
		|y-z| &\leq |y-x| + |\pi_1(x) - \pi_1(z) | + |\pi_n(x) - \pi_n(z) | \\
		&\leq |y-x| + |\pi_1(z)| + |\pi_n(x) - \pi_n(z) | .
	\end{align}
	Since $Q(S) \subseteq Q_{M+1}$, it follows immediately that  $|\pi_n(x) - \pi_n(z) | \lesssim \ell(Q_{M+1})$. Also, since $z \in b_{Q_{M+1}}$, it follows from Proposition \ref{p:b_Q} \eqref{prop:b_Q-subset} and Lemma \ref{l:stop-size-cutoff} that 
	\begin{align}
		\pi_1(z) \leq \taustop d_{M+1,k}(\pi_n(z)) \sim \ell(Q_{M+1}). 
	\end{align}
	As we proved in Lemma \ref{l:geometric-2}, we have $\ell(Q_{M+1}) \lesssim |x-y|$, which finishes the proof of the first assertion. 
	
	The set $X$ is constructed as follows. Let $V \in \calT^{M+1}(k)$ such that $Q_{M+1} \in V$ and let $W \in \calT^{M+1}(k)$ such that $Q(W) \subseteq R_M$. By Lemma \ref{l:tops-are-nearby}, we have that $\Afirst B_{Q(V)} \cap \Afirst B_{Q(W)} \neq \emptyset$. Thus, we let $X$ be the set constructed by applying Proposition \ref{p:PI-in-layers} with this data. It is easy to see that this set satisfies the desired conclusion. 
\end{proof}

\noindent\textbf{Case 3.} Suppose that $N = M+1$ or $N = M$. \\ 

In this final case, we simply find an Ahlfors regular PI space in $\envelope(k)$ which contains the endpoints $x$ and $y$. 

\begin{lem}
	There is a set $X \subseteq \envelope(k)$ which is Ahlfors $d$-regular and supports a $(1,1)$-Poincar\'e inequality and satisfies $x,y \in X$. In particular, \eqref{e:PI-max} holds.  
\end{lem} 

\begin{proof}
	If $N=M+1$, let $W \in \calT^N(k)$ such that $Q(W) \subseteq R_{M}$ and $\pi(y) \in 3\sqrt{n}\lambdadyadic B_{Q(W)}$ (such a cube exists by the minimality of $Q_0$). If $N = M$, let $W = T$. In either case, we have by Lemma \ref{l:tops-are-nearby} that 
	\begin{align}
		\Afirst B_{Q(S)} \cap \Afirst B_{Q(W)} \neq \emptyset. 
	\end{align}
	Let $X$ be the set constructed by applying Proposition \ref{p:PI-in-layers} with this data. The conclusion is then immediate.
\end{proof}

\appendix

\section{Bi-Lipschitz maps associated to the bridges}\label{a:useful}

In this section, we fill in the details for the proofs of Lemma \ref{l:useful}. We start by explicitly defining the map $H$ that we will work with. 
By scaling, we will suppose that $\ell(I) = 1.$ Define a map $\tilde{h} \colon [-1,2] \to \R^n$ (piecewise) as follows. First, let $\tilde{h} \colon I \times [-1,0] \to Q$ be an isometry with $\tilde{h}(x,0) = h(x,0) =x$. For $x \in I$ and $t \in [0,1]$ let 
\begin{align}
	\tilde{h}(x,t) =
	\begin{cases}
		h(x,2t) & \mbox{ if } t \in [0,1/2]; \\
		h(x,1) + 2(t-1/2)u(h(x,1)) & \mbox{ if } t \in [1/2,1].
	\end{cases}
\end{align} 
Finally, let $h' \colon I \times [1,2] \to R$ be a surjective map such that $h'(x,1) = h(x,1)$ and  
\begin{align}
	|h'(x,t) - h'(y,s)| = \frac{\ell(R)}{\ell(Q)} |(x,t) - (y,s)| \mbox{ for all } (x,t) , (y,s) \in I \times [1,2],
\end{align}
and set 
\begin{align}
	\tilde{h} \coloneqq u \circ h'. 
\end{align}

We have that $\tilde{h}$ is continuous since it is continuous on each component and agrees at the boundaries. Also, $\tilde{h}$ is an isometry (1-bi-Lipschitz) on $I \times [-1,0]$ and $\ell(R)/\ell(Q)$-bi-Lipschitz on $I \times [1,2]$. It follows from the definitions and \eqref{e:all-in-ball} that  
\begin{align}\label{e:tilde-h}
	\tilde{h}(I \times [-1,0]) = Q , \  \tilde{h}(I \times [1,2]) = u(R) = K \mbox{ and } \tilde{h}(I \times [-1,2]) \subseteq B. 
\end{align}
Finally, one can also check, using \eqref{e:useful-1}--\eqref{e:useful-3} and the fact that $u$ is $1$-Lipschitz, that if $x,y \in I$ and $t,s \in [0,1]$ then 
\begin{align}
	\delta | x-y| \leq | \tilde{h}(x,t) - \tilde{h}(y,t) | &\leq 2|x-y|; \label{e:useful-2'} \\
	| \tilde{h}(x,1) - \tilde{h}(y,1) | &\leq 2C_1  \delta |x-y|;  \label{e:useful-2.5'} \\
	|\tilde{h}(x,t) - \tilde{h}(x,s)| &\leq C_1|t-s|. \label{e:useful-3'}
\end{align}

We define $H \colon I \times [-1,2] \to \R^{n+1}$ by 
\begin{align}\label{e:def-H}
	H(x,t) \coloneqq 
	\begin{cases}
		\left(\tilde{h}(x,t) , g_1(\tilde{h}(x,t))\right) & \mbox{ for } t \in [-1,0]; \\
		\left(\tilde{h}(x,t),\left(1-t \right) g_1(\tilde{h}(x,0)) +  t g_2(\tilde{h}(x,1))\right) & \mbox{ for } t \in [0,1] ; \\
		\left(\tilde{h}(x,t) , g_2(\tilde{h}(x,t))\right) &\mbox{ for } t \in [1,2]. 
	\end{cases}
\end{align}

\begin{lem}
	The map $H$ takes values in $B \times \R$ and satisfies \eqref{e:useful-4}--\eqref{e:useful-5.6}. 
\end{lem}

\begin{proof}
	That $H$ takes values in $B \times \R$ and satisfies \eqref{e:useful-4} and \eqref{e:useful-4.5} follows directly from \eqref{e:tilde-h}. Let us check \eqref{e:useful-5.6}. Fix $x \in I$ and $t \in [0,1].$ Recalling that $\tilde{h}(x,0) \in I$ and $\tilde{h}(x,1) \in u(J)$, by applying \eqref{e:useful-1} we have 
	\begin{align}\label{e:old-6}
		g_1(\tilde{h}(x,0)) - \pi_{\R}(H(x,t)) \geq t \gamma.
	\end{align}
	Similarly, 
	\begin{align}\label{e:old-7}
		\pi_{\R}(H(x,t)) - g_2(\tilde{h}(x,1)) = (1-t)\left(  g_1(x) - g_2(\tilde{h}(x,1) ) \right) \geq (1-t)\gamma.
	\end{align}
	Applying \eqref{e:useful-3'}, \eqref{e:old-6}, and using the fact that $g_1$ is $\eta$-Lipschitz, if $\eta = \eta(C_1,\gamma)$ is chosen small enough, we have 
	\begin{align}\label{e:new-6}
		\pi_{\R}(H(x,t)) &\leq g_1(\tilde{h}(x,0)) - t \gamma \\
		&\leq g_1(\tilde{h}(x,t)) + |g_1(\tilde{h}(x,t)) - g_1(\tilde{h}(x,0))| - t \gamma \\
		&\leq g_1(\tilde{h}(x,t)) + (C_1\eta - \gamma) t \leq g_1(\tilde{h}(x,t))  - \gamma t/2.
	\end{align}
	Similarly, applying \eqref{e:useful-3'}, \eqref{e:old-7}, and using the fact that $g_2$ is $\theta$-Lipschitz, if $\theta = \theta(C_1,\gamma)$ is chosen small enough, we have 
	\begin{align}\label{e:new-7}
		\pi_{\R}(H(x,t)) \geq g_2(\tilde{h}(x,t)) + (1-t)\gamma/2. 
	\end{align}
	Inequalities \eqref{e:new-6} and \eqref{e:new-7} imply \eqref{e:useful-5.6}.
\end{proof}

Next, we shall check \eqref{e:banana-4} and \eqref{e:banana-5}, which will be important for checking that $H$ is bi-Lipschitz. 

\begin{lem}
	The estimates \eqref{e:banana-4} and \eqref{e:banana-5} hold.
\end{lem}

\begin{proof}
	The proofs of \eqref{e:banana-4} and \eqref{e:banana-5} are similar, and so we will only show \eqref{e:banana-5}. Fix $x \in X$. If $x \in X_{\rm down}$ then \eqref{e:banana-5} is trivial (just take $a = x$). Suppose instead that $x \in H(I \times [-1,1])$. Let $(w,t) \in I \times [-1,1]$ such that $x = H(w,t).$ Then, let $a_2 = H(w,1).$ The proof split into two cases. 
	
	Suppose first that $t \in [0,1]$. Since $H$ is Lipschitz, we have
	\begin{align}\label{e:peach-2}
		| x - a_2| \lesssim 1 - t.
	\end{align}
	On the other hand, it follows from \eqref{e:useful-5.6} that 
	\begin{align}
		|\pi_1(x)-\pi_1(a_2)| &\geq \lambda g_2(\pi_n(x)) + c\gamma \cdot (1-t) - \lambda g_2(\pi_n(a)) \\
		&\geq c\gamma |x-a_2| + \lambda g_2(\pi_n(x)) - \lambda g_2(\pi_n(y)). 
	\end{align}
	The estimate \eqref{e:banana-5} now follows directly from Lemma \ref{l:checking-the-gluing-condition} so long as $\theta = \theta(\gamma)$ is sufficiently small. 
	
	Suppose now that $t \in [-1,0]$. Let $x' = H(w,0)$, so that 
	\begin{align}\label{e:x=x'}
		|x-x'| \lesssim 1 \sim |x-a_2|. 
	\end{align}
	Repeating the estimates above for the point $x'$, we have 
	\begin{align}\label{e:x'-est}
		|\pi_1(x')-\pi_1(a_2)| \geq  \lambda g_2(\pi_n(x')) +c\gamma - \lambda g_2(\pi_n(y)). 
	\end{align} 
	Since $x$ and $x'$ lie in the graph of $g_1$, we can use \eqref{e:x=x'}, \eqref{e:x'-est} and the fact that $g_1$ is $\eta$-Lipschitz and $\lambda g_2$ is $\lambda \theta$-Lipschitz, to get 
	\begin{align}
		|\pi_1(x) - \pi_1(a_2)| &\geq |\pi_1(x') - \pi_1(a_2)| - C\eta \\
		&\geq (c\gamma - C\eta - C\lambda \theta) +  \lambda g_2(\pi_n(x)) - \lambda g_2(\pi_n(y)) 
	\end{align}
	Assuming $\eta = \eta(\gamma)$ and $\theta = \theta(\gamma)$ are small enough, this implies 
	\begin{align}
		|\pi_1(x) - \pi_1(a_2)|	&\geq c\gamma |x-a_2|  +  \lambda g_2(\pi_n(x)) - \lambda g_2(\pi_n(y)).
	\end{align}
	The estimate \eqref{e:banana-5} now follows directly from Lemma \ref{l:checking-the-gluing-condition} so long as $\theta = \theta(\gamma)$ is sufficiently small.  
\end{proof}

The only thing left to check is that $H$ is bi-Lipschitz. We start with the following preliminary result. 

\begin{lem}\label{l:piecewise-bi-lip}
	The map $H$ is bi-Lipschitz when restricted to either of $I \times [-1,0]$, $I \times [0,1]$ or $I \times [1,2]$. 
\end{lem}

\begin{proof}
	That $H$ is bi-Lipschitz on when restricted to $I \times [-1,0]$ and $I \times [1,2]$ (since it is Lipschitz in the second component and bi-Lipschitz in the first).   To check the last case, fix $x,y \in I$, $t,s \in [0,1]$. Denote 
	\[w = \tilde{h}(x,1) \mbox{ and } z = \tilde{h}(y,1) \]
	and observe for later that, after applying \eqref{e:useful-2.5'} and the fact that $g_2$ is $\theta$-Lipschitz and $u$ is $1$-Lipschitz, we have 
	\begin{align}\label{e:d-w-z}
		| g_2(w) - g_2(z) | & \leq \theta |w-z|  \leq 2C_1 \theta \delta. 
	\end{align}
	By \eqref{e:useful-2'} and \eqref{e:useful-3'}, we have 
	\begin{align}
		| \tilde{h}(x,t) - \tilde{h}(y,s) | \leq | \tilde{h}(x,t) - \tilde{h}(y,t) | + | \tilde{h}(y,t) - \tilde{h}(y,s) | \lesssim_{C_1} |x-y| + |t-s|. 
	\end{align}
	Furthermore, using the second estimate in \eqref{e:useful-1}, \eqref{e:d-w-z}, and the fact that $g_1$ is $\eta$-Lipschitz, it follows that 
	\begin{align}
		\begin{split}\label{e:upper-H-2}
			&| (1-t) g_1(x) +  g_2(w) - (1-s) g_1(y) -  s g_2(z) | \\
			&\hspace{2em}  \leq (1-t)|g_1(x) - g_1(y)| + t |g_2(w) - g_2(z)|  + |s-t| \left(  |g_1(y)| + |g_2(z)| \right)  \\
			&\hspace{2em} \leq (\eta + 2C_1 \theta \delta )  |x-y| + 2C_1|t-s|. 
		\end{split}
	\end{align}
	Combining both previous sets of inequality, we deduce that $H$ is Lipschitz. To prove that $H$ is bi-Lipschitz, we first choose $\eta= \eta(\gamma,\delta,C_1)$ and $\theta = \theta(\gamma,C_1)$ small enough such that 
	\begin{align}\label{e:alpha}
		\frac{\eta}{\delta} + 2 \theta C_1 \leq \frac{\gamma}{4C_1},
	\end{align}
	then split into two cases. Suppose first that $2C_1|t-s| \geq \delta |x-y|.$ 
	By the first inequality in \eqref{e:useful-1},
	\begin{align}
		| g_1(x) - g_2(w) | \geq \gamma. 
	\end{align}
	Applying this, in conjunction with the first estimate in \eqref{e:useful-1}, \eqref{e:d-w-z}, \eqref{e:alpha}, the fact that $g_1$ is $\eta$-Lipschitz, and our assumption that $2C_1|t-s| \geq \delta |x-y|,$ it follows that 
	\begin{align}
		&| (1-t)g_1(x) + tg_2(w) - (1-s)g_1(y) -  s g_2(z) | \\
		&\hspace{2em} \geq |t-s|  | g_1(x) - g_2(w) | - (1-s) |g_1(x) - g_1(y) |  -  s | g_2(w) - g_J(z) | \\
		&\hspace{2em} \geq \gamma |t-s| - (\eta + 2C_1\theta \delta) |x-y| \geq \gamma |t-s| - 2C_1\left( \frac{\eta}{\delta} + 2\theta C_1 \right)  |t-s | \\
		&\hspace{2em} \gtrsim \gamma |t-s| .
	\end{align} 
	Suppose instead that $2C_1|t-s| <  \delta |x-y|.$ Then, using \eqref{e:useful-2}, we have
	\begin{align}
		|h(x,t) - h(y,s) | &\geq |h(x,t) - h(y,t)| - |h(y,t) - h(y,s)| \geq  \delta |x-y| - C_1|t-s| \\
		&\gtrsim \delta |x-y|.
	\end{align}
	The completes the proof.
\end{proof}

The following finishes the proof of Lemma \ref{l:useful}. 

\begin{lem}
	The map $H$ is bi-Lipschitz on $I \times [-1,2]$. 
\end{lem}

\begin{proof}
	Note that $H$ is globally Lipschitz since it is Lipschitz on each component, so we only need to check the lower bounds. Fix $x,y \in I$ and $s,t \in [-1,2]$. For brevity, let us denote $w = H(x,t)$ and $z = H(y,s)$. \\
	
	\noindent\textbf{Case 1.} Suppose that $t \in [-1,0]$ and $s \in [0,1]$.\\
	
	Consider the point $a_1 = H(y,0)$. By \eqref{e:banana-4}, we know that $|w-z| \gtrsim |w-a| + |a-z|$. This gives the lower bound $|H(x,t) - H(y,s)| \gtrsim |x-y| + |t-s|$ since $(y,0) \in I \times \{0\}$ and $H$ is bi-Lipschitz on each of $I \times [-1,0]$ and $I \times [0,1]$ by Lemma \ref{l:piecewise-bi-lip}.  \\
	
	\noindent\textbf{Case 2.} Suppose that $t \in [-1,0]$ and $s \in [1,2]$. \\
	
	Consider the points $w' = H(x,0) \in g_1(I)$ and $z' = H(y,1) \in g_2(u(J))$. By \eqref{e:useful-1}, these points satisfy $|w'-z'| \geq \gamma \ell(I)$. Since $g_1$ is $\eta$-Lipschitz and $H$ is Lipschitz (so that $|w-w'| \lesssim \ell(I)$), we have 
	\begin{align}
		|\pi_1(w-w')| = |g_1(\pi_n(w)) - g_1(\pi_n(w'))| \lesssim \eta \ell(I). 
	\end{align}
	Similarly, $|\pi_1(z-z')| \lesssim \theta \ell(I)$. Putting these estimates together, we get 
	\begin{align}
		|w-z| \geq |\pi_1(w - z)| &\geq |\pi_1(w' - z')| - |\pi_1(w-w')| - |\pi_1(z-z')| \\
		&\geq \gamma \ell(I) - C\eta \ell(I) - C\theta \ell(I) \gtrsim \gamma \ell(I) \\
		&\gtrsim \gamma \left( |x-y| + |t-s| \right). 
	\end{align}
	\hspace{1em}
	
	\noindent\textbf{Case 3.} Suppose that $t \in [0,1]$ and $s \in [1,2]$. \\
	
	This has is similar to Case 1, and the proof passes via Lemma \ref{l:checking-the-gluing-condition}. In this case, one considers the point $a = H(x,1)$. We will leave the details to the reader. 
\end{proof}

\bibliography{Ref.bib}
\bibliographystyle{alpha}
	
\end{document}